\documentclass[11pt,oneside]{amsart}

\usepackage{amsmath}

\usepackage{tikz-cd}

\usepackage{fix-cm}

\DeclareMathAlphabet{\mathcal}{OMS}{cmsy}{m}{n}

\usepackage[arrow,curve,matrix,tips,2cell]{xy}

\usepackage[hmargin=2cm,vmargin=1.5cm]{geometry}

\usepackage{anyfontsize}
\usepackage{mathrsfs}

\usepackage{bm}

\usepackage{bbm}

\usepackage{fancyhdr}
\usepackage[obeyspaces]{url}

\usepackage{tensor}

\usepackage{mathtools}

\usepackage{accents}

\usepackage{amscd}
\usepackage{amssymb}
\usepackage{latexsym}
\usepackage{url}

\usepackage{graphicx}

\usepackage[hidelinks]{hyperref}

\usepackage{enumitem}
\usepackage{leftindex}

\usepackage{tikz}
\usetikzlibrary{automata, positioning, arrows}

\tikzset{
->, 
>=stealth', 
node distance=3cm, 
initial text=$ $, 
}

\newtheorem{theorem}{Theorem}[section]
\newtheorem*{theorem*}{Theorem}
\newtheorem{lemma}[theorem]{Lemma}
\newtheorem*{lemma*}{Lemma}
\newtheorem{corollary}[theorem]{Corollary}
\newtheorem{proposition}[theorem]{Proposition}

\newtheorem{remark}[theorem]{Remark}
\newtheorem{definition}[theorem]{Definition}
\newtheorem*{definition*}{Definition}

\newtheorem{question}[theorem]{Question}
\newtheorem{questions}[theorem]{Questions}
\newtheorem*{question*}{Question}

\newtheorem{example}[theorem]{Example}
\newtheorem{examples}[theorem]{Examples}

\newtheorem{observation}[theorem]{Observation}
\newtheorem{observations}[theorem]{Observations}

\makeatletter
\def\revddots{\mathinner{\mkern1mu\raise\p@
\vbox{\kern7\p@\hbox{.}}\mkern2mu
\raise4\p@\hbox{.}\mkern2mu\raise7\p@\hbox{.}\mkern1mu}}
\makeatother 
\newcommand{\bgl}{\begin{equation}} 
\newcommand{\egl}{\end{equation}}
\newcommand{\bgloz}{\begin{equation*}} 
\newcommand{\egloz}{\end{equation*}}
\newcommand{\bgln}{\begin{eqnarray}} 
\newcommand{\egln}{\end{eqnarray}}
\newcommand{\bglnoz}{\begin{eqnarray*}} 
\newcommand{\eglnoz}{\end{eqnarray*}}
\newcommand{\btheo}{\begin{theorem}}
\newcommand{\etheo}{\end{theorem}}
\newcommand{\btheooz}{\begin{theorem*}}
\newcommand{\etheooz}{\end{theorem*}}

\newcommand{\blemma}{\begin{lemma}}
\newcommand{\elemma}{\end{lemma}}
\newcommand{\blemmaoz}{\begin{lemma*}}
\newcommand{\elemmaoz}{\end{lemma*}}
\newcommand{\bproof}{\begin{proof}}
\newcommand{\eproof}{\end{proof}}
\newcommand{\bbew}{\begin{beweis}}
\newcommand{\ebew}{\end{beweis}}
\newcommand{\bremark}{\begin{remark}}
\newcommand{\eremark}{\end{remark}}
\newcommand{\bdefin}{\begin{definition}}
\newcommand{\edefin}{\end{definition}}
\newcommand{\bdefinoz}{\begin{definition*}}
\newcommand{\edefinoz}{\end{definition*}}
\newcommand{\bex}{\begin{example}}
\newcommand{\eex}{\end{example}}
\newcommand{\bexs}{\begin{examples}}
\newcommand{\eexs}{\end{examples}}

\newcommand{\bobs}{\begin{observation}}
\newcommand{\eobs}{\end{observation}}
\newcommand{\bobss}{\begin{observations}}
\newcommand{\eobss}{\end{observations}}

\newcommand{\bprop}{\begin{proposition}}
\newcommand{\eprop}{\end{proposition}}
\newcommand{\bcor}{\begin{corollary}}
\newcommand{\ecor}{\end{corollary}}
\newcommand{\bfa}{\begin{cases}} 
\newcommand{\efa}{\end{cases}}

\newcommand{\bquestion}{\begin{question}}
\newcommand{\equestion}{\end{question}}
\newcommand{\bquestions}{\begin{questions}}
\newcommand{\equestions}{\end{questions}}
\newcommand{\bquestionoz}{\begin{question*}}
\newcommand{\equestionoz}{\end{question*}}
\newtheorem{introtheorem}{Theorem}

\newtheorem{introcor}[introtheorem]{Corollary}

\newcommand{\cF}{\mathcal F}

\newcommand{\cI}{\mathcal I}

\newcommand{\cL}{\mathcal L}
\newcommand{\cM}{\mathcal M}

\newcommand{\cO}{\mathcal O}
\newcommand{\cP}{\mathcal P}

\newcommand{\cS}{\mathcal S}

\newcommand{\cU}{\mathcal U}

\def\Cz{\mathbb{C}}

\def\Kz{\mathbb{K}}

\def\Nz{\mathbb{N}}

\def\Hz{\mathbb{H}}
\def\Iz{\mathbb{I}}

\def\Qz{\mathbb{Q}}
\def\Rz{\mathbb{R}}
\def\Tz{\mathbb{T}}

\def\Zz{\mathbb{Z}}

\newcommand{\fC}{\mathfrak C}

\newcommand{\mfi}{\mathfrak i}

\newcommand{\mfp}{\mathfrak p}
\newcommand{\mfq}{\mathfrak q}

\newcommand{\an}[1]{``#1''} 
\newcommand{\ti}{\tilde}

\newcommand{\ri}{\rightarrow}

\newcommand{\ma}{\mapsto} 
\newcommand{\onto}{\twoheadrightarrow} 
\newcommand{\into}{\hookrightarrow} 
\newcommand{\Rarr}{\Rightarrow} 
\newcommand{\Larr}{\Leftarrow} 
\newcommand{\LRarr}{\Leftrightarrow} 

\def\SEMI{\mbox{$\times\kern-2pt\vrule height5pt width.6pt \kern3pt $}}

\newcommand{\Spec}{{\rm Spec\,}} 

\newcommand{\id}{{\rm id}}

\renewcommand{\ker}{{\rm ker}\,}

\newcommand{\reg}{^\times} 
\newcommand{\ev}{\operatorname{ev}} 
\newcommand{\lspan}{{\rm span}} 
\newcommand{\clspan}{\overline{\lspan}} 
\newcommand{\abs}[1]{\lvert#1\rvert} 
\newcommand{\norm}[1]{\left\|#1\right\|} 

\newcommand{\dop}{\text{: }} 
\newcommand{\ilim}{\varinjlim} 
\newcommand{\e}{{\rm e}} 
\newcommand{\supp}{{\rm supp}} 
\newcommand{\lge}{\left\{} 
\newcommand{\rge}{\right\}} 
\newcommand{\lsp}{\left\langle} 
\newcommand{\rsp}{\right\rangle} 
\newcommand{\gekl}[1]{\lge #1 \rge} 
\newcommand{\spkl}[1]{\lsp #1 \rsp} 
\newcommand{\menge}[2]{\gekl{ #1 \dop #2 }} 
\newcommand{\bmg}{\bm{g}}
\newcommand{\bmh}{\bm{h}}

\newcommand{\bmx}{\bm{x}}
\newcommand{\bmy}{\bm{y}}

\newcommand{\bmU}{\bm{U}}

\newcommand{\oset}[2]{%
  \mathop{#2}\limits^{\vbox to -1.66ex{%
  \kern -1.4ex\hbox{$#1$}\vss}}}

\newcommand{\B}{\mathscr{B}}
\newcommand{\C}{\mathscr{C}}

\newcommand{\Gn}{{G^{(0)}}}

\newcommand{\acts}{\curvearrowright}
\newcommand{\stca}{\curvearrowleft}

\newcommand{\inte}{{\rm int\,}}

\newcommand{\sing}{{\rm sing}}
\newcommand{\ess}{{\rm ess}}

\newcommand{\im}{{\rm im}\,}

\newcommand{\Cc}{\mathscr{C}_c}

\newcommand{\Gu}{G^{(0)}}
\newcommand{\tiG}{\ti{G}}
\newcommand{\tiGu}{\tiG^{(0)}}

\newcommand{\tif}{\ti{f}}
\newcommand{\tis}{\ti{s}}
\newcommand{\tir}{\ti{r}}

\newcommand{\tiA}{\ti{A}}
\newcommand{\tiI}{\ti{I}}
\newcommand{\tiU}{\ti{U}}
\newcommand{\tiX}{\ti{X}}

\newcommand{\osupp}{\supp^\circ}

\newcommand{\TF}{{\rm TF}}

\newcommand{\SF}{{\rm SF}}

\newcommand{\Isot}{{\rm Iso}}

\newcommand{\congto}{\xrightarrow{\raisebox{-.5ex}[0ex][0ex]{$\sim$}}}

\DeclareFontFamily{U}{mathb}{\hyphenchar\font45}
\DeclareFontShape{U}{mathb}{m}{n}{
      <5> <6> <7> <8> <9> <10>
      <10.95> <12> <14.4> <17.28> <20.74> <24.88>
      mathb10
      }{}
\DeclareSymbolFont{mathb}{U}{mathb}{m}{n}

\DeclareMathSymbol{\sqbullet}{1}{mathb}{"0D}
\usepackage{newtxtext}
\usepackage{newtxmath}

\newcommand{\ord}{{\rm ord}}

\newcommand{\J}{\mathscr{J}}

\newcommand{\oX}{\overline{X}}

\begin{document}

\title{On Hausdorff covers for non-Hausdorff groupoids}

\thispagestyle{fancy}

\author{Kevin Aguyar Brix}
\author{Julian Gonzales}
\author{Jeremy B. Hume}
\author{Xin Li}

\address[K. A. Brix]{Department of Mathematics and Computer Science, University of Southern Denmark, 5230 Odense, Denmark}
\email{kabrix@imada.sdu.dk}

\address{Julian Gonzales, School of Mathematics and Statistics, University of Glasgow, University Place, Glasgow G12 8QQ, United Kingdom}
\email{julian.gonzales@glasgow.ac.uk}

\address{Jeremy B. Hume, School of Mathematics and Statistics, Carleton University, 4302 Herzberg Laboratories, Ottawa, ON, K1S 5B6, Canada}
\email{jeremybhume@gmail.com}

\address{Xin Li, School of Mathematics and Statistics, University of Glasgow, University Place, Glasgow G12 8QQ, United Kingdom}
\email{Xin.Li@glasgow.ac.uk}


\thanks{Brix was supported by a DFF-international postdoc (case number 1025--00004B), a Starting grant from the Swedish Research Council (2023-03315) and the Carlsberg Foundation (grant CF23–1328). Hume was supported by NSERC Discovery Grant RGPIN-2021-03834. This project has received funding from the European Research Council (ERC) under the European Union's Horizon 2020 research and innovation programme (grant agreement No. 817597).}

\begin{abstract}
We develop a new approach to non-Hausdorff {\'e}tale groupoids and their algebras based on Timmermann's construction of Hausdorff covers. As an application, we completely characterise when singular ideals vanish in Steinberg algebras over arbitrary rings. We also completely characterise when $C^*$-algebraic singular ideals have trivial intersection with the non-Hausdorff analogue of subalgebras of continuous, compactly supported functions. This leads to a characterisation when $C^*$-algebraic singular ideals vanish for groupoids satisfying a finiteness condition. Moreover, our approach leads to further sufficient vanishing criteria for singular ideals and reduces questions about simplicity, the ideal intersection property, amenability and nuclearity for non-Hausdorff {\'e}tale groupoids to the Hausdorff case.
\end{abstract}

\maketitle

\setlength{\parindent}{0cm} \setlength{\parskip}{0.5cm}

\section{Introduction}

Topological groupoids arise naturally in a variety of areas such as dynamics, topology, geometry, group theory and $C^*$-algebras, building bridges between all these areas of mathematics. The link between topological groupoids and $C^*$-algebras is established by the construction of convolution algebras \cite{Ren,Con78,Con82}. Many important classes of $C^*$-algebras arise from groupoids (for example AF-algebras, Cuntz-Krieger algebras, graph algebras etc., see also \cite{Li20}), and groupoid models are very helpful for studying structural properties or interesting invariants of the corresponding $C^*$-algebras. For topological groupoids which are Hausdorff (as well as locally compact and {\'e}tale), there are satisfactory results characterising fundamental properties of groupoid $C^*$-algebras such as simplicity, the ideal intersection property or nuclearity in terms of the underlying groupoids (see \cite{Exel,BCFS,KKLRU,AR}), and the notion of Cartan subalgebras allows to reconstruct groupoids from their $C^*$-algebras \cite{Kum,Ren08}. However, several of these results, such as characterisations of simplicity, the ideal intersection property, or reconstruction results using Cartan subalgebras, do not carry over to non-Hausdorff groupoids (see \cite{Exel} and also a related example due to Skandalis in \cite{Ren91}), and new challenges arise in the non-Hausdorff setting because we have to work with compact sets which are not closed and functions which are not continuous. At the same time, important examples of groupoids arising as groupoid of germs, from foliations or self-similar groups typically are non-Hausdorff (see \cite{EP17,Nek18a,Nek18b} for some examples), which provides a strong motivation to systematically develop a better understanding of non-Hausdorff groupoids and their $C^*$-algebras. 

In this context, a new construction called essential groupoid $C^*$-algebras has been introduced for non-Hausdorff groupoids \cite{EP22,KM}, and simplicity as well as the ideal intersection property have been characterised for this new construction in terms of the underlying groupoids \cite{CEPSS,KKLRU}. While these recent developments and their algebraic analogues in the setting of Steinberg algebras have led to progress in our understanding of non-Hausdorff groupoids, it remains an open question when essential groupoid $C^*$-algebras agree with the more classical reduced groupoid $C^*$-algebras, and there is an analogous open question for Steinberg algebras (see \cite{StSz21,StSz23,Yos,GNSV} for progress in the case of groupoids arising from self-similar groups). Moreover, the construction of essential groupoid $C^*$-algebras has been specifically designed to resolve issues arising in the study of simplicity or the ideal intersection property, so that in other contexts, essential groupoid $C^*$-algebras might not be good substitutes for reduced groupoid $C^*$-algebras (because of bad functoriality properties, for example).

The goal of this paper is to develop a new approach to non-Hausdorff groupoids and their algebras, based on the construction of Hausdorff covers as in \cite{Tim} (similar constructions appeared in \cite{KS,Tu}). Given a non-Hausdorff {\'e}tale groupoid $G$ with Hausdorff unit space, its Hausdorff cover $\tiG$ is a Hausdorff {\'e}tale groupoid given by the closure of $G$ in the space of non-empty closed subsets of $G$ with respect to the Fell topology \cite{Fell}. Alternatively, we can think of $\tiG$ as the Gelfand spectrum of the smallest commutative $C^*$-algebra of bounded, complex-valued functions on $G$ (with respect to point-wise multiplication) containing $\Cc(G) \coloneq \lspan(\menge{C_c(U)}{U \subseteq G \text{ open, Hausdorff}})$. Here $C_c(U)$ is the algebra of continuous functions $U \to \Cz$ with compact support contained in $U$, and we extend a function in $C_c(U)$ by zero to view it as a function on $G$. Using the Hausdorff cover, we are able to develop new tools in the non-Hausdorff setting as well as reduce questions about non-Hausdorff groupoids to the Hausdorff case. Here is a summary of our main achievements.
\setlength{\parindent}{0cm} \setlength{\parskip}{0cm}

\begin{itemize}
    \item We establish a complete characterisation when singular ideals (the definition is explained below) vanish in Steinberg algebras over arbitrary rings in terms of a groupoid property. In combination with \cite[Theorem~A]{StSz21}, this leads to a complete characterisation for simplicity of Steinberg algebras over fields, which in particular answers \cite[Question~1]{CEPSS}. 
    \item We also completely characterise when the $C^*$-algebraic singular ideal of a non-Hausdorff {\'e}tale groupoid $G$ has trivial intersection with $\Cc(G)$. This leads to a characterisation when the $C^*$-algebraic singular ideal vanishes, under the assumption that the closure of the unit space $\Gu$ of $G$ has finite fibres over $\Gu$. In combination with existing simplicity characterisations for essential groupoid $C^*$-algebras in \cite{KKLRU}, this yields a characterisation in terms of a groupoid property when the reduced groupoid $C^*$-algebra of $G$ is simple (under further assumptions on $G$), which provides a partial answer to \cite[Question~2]{CEPSS}. Another consequence is that, for a non-Hausdorff ample groupoid $G$ satisfying our finiteness condition, simplicity of the complex Steinberg algebra of $G$ implies simplicity of the reduced groupoid $C^*$-algebra of $G$, which provides a partial answer to \cite[Question~3]{CEPSS}. As a further application, we obtain a conceptual explanation for the results on groupoids associated with contracting self-similar groups in \cite{GNSV}.
    \item For a non-Hausdorff {\'e}tale groupoid $G$ which can be covered by countably many open bisections, we describe singular ideals in terms of the Hausdorff cover and give sufficient conditions for their vanishing, both in the $C^*$-algebraic and the algebraic context. This complements existing results on vanishing of singular ideals (such as those in \cite{CEPSS,NS}) by explaining them in terms of the Hausdorff cover.
    \item We show, again for a non-Hausdorff {\'e}tale groupoid $G$ which can be covered by countably many open bisections, that the ideal intersection property for the essential groupoid $C^*$-algebra $C^*_\ess(G)$ is equivalent to the ideal intersection property of the reduced $C^*$-algebra $C^*_r(\tiG_{\ess})$, where $\tiG_{\ess}$ is the restriction of $\tiG$ to a closed invariant subset of the unit space of $\tiG$. Since $\tiG_{\ess}$ is Hausdorff, this explains why the ideal intersection property for essential groupoid $C^*$-algebras is much more tractable than for reduced groupoid $C^*$-algebras and allows for characterisations as in the Hausdorff case. In particular, we obtain a conceptual explanation for results on simplicity or the ideal intersection property for essential $C^*$-algebras of non-Hausdorff {\'e}tale groupoids in \cite{CEPSS,EP22,KM,KKLRU}.
    \item We show that nuclearity of the reduced groupoid $C^*$-algebra of a non-Hausdorff {\'e}tale groupoid $G$ is equivalent to topological amenability of its Hausdorff cover $\tiG$. If $G$ is $\sigma$-compact, we further show, based on results in \cite{Ren15}, that this is equivalent to amenability of $G$ itself (amenability in the sense of \cite[Definition~2.7]{Ren15}). This completely reduces questions about nuclearity and amenability to the Hausdorff case. 
\end{itemize}
\setlength{\parindent}{0cm} \setlength{\parskip}{0.5cm}

Let us now describe our main results in more detail. Our central objects of study are {\'e}tale groupoids $G$, i.e., small categories with invertible morphisms, equipped with a topology such that all structure maps become local homeomorphisms. We always assume that $G$ is locally compact and that the unit space $\Gn$ is Hausdorff. The space $\Cc(G)$ introduced above is an algebra with respect to convolution, and its completion with respect to the norm from the left regular representation is the reduced groupoid $C^*$-algebra $C^*_r(G)$. The ($C^*$-algebraic) singular ideal $J$ consists of those elements of $C^*_r(G)$ which are only non-zero on a set whose source has empty interior. The essential groupoid $C^*$-algebra $C^*_{\ess}(G)$ is then defined as the quotient $C^*_r(G) / J$. In the purely algebraic context, we consider ample groupoids $G$, which are {\'e}tale groupoids with totally disconnected unit spaces. In this case, given a ring $R$, the Steinberg algebra $RG$ is given by $RG \coloneq \lspan_R(\menge{C_c(U,R)}{U \subseteq G \text{ open, Hausdorff}})$. Here $C_c(U,R)$ is the algebra of locally constant functions $U \to R$ with compact support contained in $U$, and we extend a function in $C_c(U,R)$ by zero to view it as a function on $G$. As is the case for $\Cc(G)$, $RG$ becomes an algebra with respect to convolution, and the algebraic analogue of the singular ideal is denoted by $J_R$ and consists of those $f \in RG$ which vanish outside of a set with empty interior. As mentioned above, it is an open problem to characterise in terms of $G$ when the singular ideal $J$ vanishes, or in the ample case, when $J_R$ vanishes.
\setlength{\parindent}{0.5cm} \setlength{\parskip}{0cm}

The Hausdorff cover $\tiG$ of $G$ is constructed as a space of subsets of $G$, equipped with the Fell topology. It comes with a canonical (non-continuous) inclusion $\iota \colon G \into \tiG$ with dense image. The dual picture mentioned above, describing $\tiG$ as a Gelfand spectrum, yields a canonical inclusion $\mfi \colon \Cc(G) \into C_c(\tiG)$, where $C_c(\tiG)$ denotes the algebra of continuous, complex-valued functions on $\tiG$ with compact support. Moreover, the map $\mfi$ extends continuously to the embedding $C^*_r(\mfi) \colon C^*_r(G) \into C^*_r(\tiG)$. For an ample groupoid $G$, we show that its Hausdorff cover $\tiG$ is again ample, and given a ring $R$, we obtain an analogous inclusion of Steinberg algebras $\mfi_R \colon RG \into R \tiG$.
\setlength{\parindent}{0cm} \setlength{\parskip}{0.5cm}

Our first main result completely characterises when $J \cap \Cc(G)$ vanishes, and when $J_R$ vanishes, for arbitrary rings. Given a ring $R$, we write $\ord(R)$ for the set of all orders of non-zero elements of $R$. Moreover, we set $R_0 \coloneq \Qz$ and $R_t \coloneq \Zz / t \Zz$, where $t$ is an integer with $t > 1$. Furthermore, given $t \in \gekl{0,2,3,\dotsc}$, we say that an {\'e}tale groupoid $G$ satisfies condition ($\cS_t$) if there exist a positive integer $n$ and open bisections $U_i \subseteq G$ for $1 \leq i \leq n$ such that 
\setlength{\parindent}{0cm} \setlength{\parskip}{0cm}

\begin{itemize}
    \item $s(U_i) = s(U_j)$ for all $1 \leq i, j \leq n$,
    \item $U_i \setminus (\bigcup_{j \neq i} U_j)$ is not empty and has empty interior for all $1 \leq i \leq n$,
    \item $\lspan_{R_t}(\menge{b_I}{I \subseteq \gekl{1, \dotsc, n}, \, \# \, I > 1, \, (\bigcap_{i \in I} U_i) \setminus (\bigcup_{j \notin I} U_j) \neq \emptyset}) \neq R_t^n$, where $b_I \coloneq \sum_{i \in I} b_i$ and $b_i$ is the $i$-th standard basis vector in $R_t^n$ (all components of $b_i$ are zero except the $i$-th, which is $1$).

\end{itemize}
\setlength{\parindent}{0cm} \setlength{\parskip}{0cm}
\begin{introtheorem}[see Theorem~\ref{thm:CharSingIdVanish}]
\label{thm:TheVeryFirstOne}
Let $G$ be a non-Hausdorff {\'e}tale groupoid. 
\setlength{\parindent}{0cm} \setlength{\parskip}{0cm}

\begin{enumerate}[label=(\Roman*)]
    \item We have $J \cap \Cc(G) \neq \gekl{0}$ if and only if $G$ satisfies condition ($\cS_0$).
    \item Let $G$ be ample and $R$ a ring. We have $J_R \neq \gekl{0}$ if and only if there exists $t \in \ord(R)$ such that $G$ satisfies condition ($\cS_t$).
\end{enumerate}
\end{introtheorem}
\setlength{\parindent}{0cm} \setlength{\parskip}{0cm}

An immediate consequence of Theorem~\ref{thm:TheVeryFirstOne} is that, for non-Hausdorff ample groupoids $G$, we have $J \cap \Cc(G) = \gekl{0}$ if and only if $J_{\Cz} = \gekl{0}$. Theorem~\ref{thm:TheVeryFirstOne} immediately implies (and thus explains) the result in \cite[Theorem~5.9]{StSz21} that, for a field $\Kz$, whether or not $J_\Kz = \gekl{0}$ depends only on the characteristic of $\Kz$. Moreover, Theorem~\ref{thm:TheVeryFirstOne} combined with \cite[Theorem~A]{StSz21} yields the following complete characterisation for simplicity of Steinberg algebras over fields, which answers \cite[Question~1]{CEPSS}.

\begin{introcor}
\label{cor:intro_1}
Let $G$ be a non-Hausdorff ample groupoid and $\Kz$ a field of characteristic $p$. Then $\Kz G$ is simple if and only if $G$ is minimal, effective and does not satisfy condition ($\cS_p$).
\end{introcor}
\setlength{\parindent}{0cm} \setlength{\parskip}{0.5cm}

Our second main result characterises when the $C^*$-algebraic singular ideal vanishes under a finiteness assumption.
\begin{introtheorem}[see Theorem~\ref{thm:JCc=0vsJ=0} and Corollary~\ref{cor:JCc=0vsJ=0}]
\label{thm:Intro_J_CSTAR}
Let $G$ be a non-Hausdorff {\'e}tale groupoid. Assume that the source map (or equivalently the range map) restricts to a finite-to-one map $\overline{\Gu} \to \Gu$. Then $J \cap \Cc(G) = \gekl{0}$ if and only if $J = \gekl{0}$, and $J = \gekl{0}$ if and only if $G$ does not satisfy condition ($\cS_0$).   
\end{introtheorem}
\setlength{\parindent}{0cm} \setlength{\parskip}{0cm}

Theorem~\ref{thm:Intro_J_CSTAR} leads to a characterisation when $C^*_r(G)$ is simple in terms of a groupoid property, if $G$ is as in Theorem~\ref{thm:Intro_J_CSTAR} and has compact unit space or can be covered by countably many open bisections (see Corollary~\ref{cor:SimpCSTAR}). This provides a partial answer to \cite[Question~2]{CEPSS}. Theorem~\ref{thm:Intro_J_CSTAR} also implies that, if $G$ is as in Theorem~\ref{thm:Intro_J_CSTAR} and in addition ample, then simplicity of $\Cz G$ implies simplicity of $C^*_r(G)$ (see Corollary~\ref{cor:SimpSteinbCSTAR}). This provides a partial answer to \cite[Question~3]{CEPSS}. Theorems~\ref{thm:TheVeryFirstOne} and \ref{thm:Intro_J_CSTAR} also provide a conceptual explanation for the result in \cite{GNSV} that for groupoids attached to contracting self-similar groups, we have $J = \gekl{0}$ if and only if $J_{\Cz} = \gekl{0}$ (see Corollary~\ref{cor:SelfSimJvsJC}).
\setlength{\parindent}{0cm} \setlength{\parskip}{0.5cm}

Here is our main result describing the images of $\mfi$ and $\mfi_R$.
\setlength{\parindent}{0cm} \setlength{\parskip}{0cm}

\begin{introtheorem}[see Theorem~\ref{thm:CharCinC}]
\label{thm:the1stone}
Suppose that $G$ is a non-Hausdorff {\'e}tale groupoid with Hausdorff cover $\tiG$. 
\begin{enumerate}[label=(\Roman*)]
    \item A function $\tif \in C_c(\tiG)$ lies in the image of $\mfi$ if and only if $\tif$ is decomposable in the sense that $\tif(\bmg) = \sum_{g \in \bmg} \ti{f}(\iota(g))$ for all $\bmg \in \tiG$.
  
    \item Now further assume that $G$ is ample. Let $R$ be a ring. A function $\tif \in R \tiG$ lies in the image of $\mfi_R$ if and only if $\tif$ is decomposable in the sense that $\tif(\bmg) = \sum_{g \in \bmg} \ti{f}(\iota(g))$ for all $\bmg \in \tiG$.
\end{enumerate}
\end{introtheorem}
Similar arguments as for Theorem~\ref{thm:the1stone} lead to decomposition results for functions in $\Cc(G)$ or $RG$ (Corollary~\ref{cor:f=sum}), which are interesting in their own right and also crucial for Theorem~\ref{thm:TheVeryFirstOne}.
\setlength{\parindent}{0cm} \setlength{\parskip}{0.5cm}

Now assume that $G$ can be covered by countably many open bisections. The set of continuity points of $\iota$ forms a subgroupoid $G_C^C$ of $G$. Let $\tiG_{\ess}$ denote the closure of $\iota(G_C^C)$ in $\tiG$. $\tiG_{\ess}$ is the restriction of $\tiG$ to a closed invariant subset of the unit space of $\tiG$, and we denote by $\ti{J}$ the kernel of the canonical projection $C^*_r(\tiG) \onto C^*_r(\tiG_{\ess})$ given by restriction. If $G$ is in addition ample, let $\ti{J}_R$ denote the kernel of the canonical projection $R \tiG \onto R \tiG_{\ess}$. Using Hausdorff covers, we establish the following descriptions and sufficient vanishing criteria for singular ideals (see Proposition~\ref{prop:Jsing=CcapJ}, Remark~\ref{rem:JR=cap} as well as Corollaries~\ref{cor:tiXsing=0} and \ref{cor:SuffCondSingVan}):
\setlength{\parindent}{0cm} \setlength{\parskip}{0cm}

\begin{enumerate}[label=(\Roman*)]
    \item We have $C^*_r(\mfi)(J) = C^*_r(\mfi)(C^*_r(G)) \cap \ti{J}$.
\setlength{\parindent}{0.5cm} \setlength{\parskip}{0cm}

\begin{itemize}
    \item 
    If $\tiG_{\ess}^{(0)} = \tiG^{(0)}$, then $J = \gekl{0}$. 
    
    \item
    Suppose that $G$ has the property that $\# \, \bmg < \infty$ for all $\bmg \in \tiG$, and that for every decomposable function $\tif \colon \tiG \to \Cz$, $\tif(\bmg) = 0$ for all $\bmg \in \tiG_{\ess}$ implies that $\tif(\bmg) = 0$ for all $\bmg \in \tiG$. Then $J = \gekl{0}$.
\end{itemize}
\setlength{\parindent}{0cm} \setlength{\parskip}{0cm}

    \item Suppose that $G$ is ample. We have $\mfi_R(J_R) = \mfi(RG) \cap \ti{J}_R$.
\setlength{\parindent}{0.5cm} \setlength{\parskip}{0cm}

\begin{itemize}
    \item 
    If $\tiG_{\ess}^{(0)} = \tiG^{(0)}$, then $J_R = \gekl{0}$ for every ring $R$.    
\end{itemize}

\end{enumerate}
\setlength{\parindent}{0cm} \setlength{\parskip}{0cm}

As explained in \S~\ref{s:SingId}, these results give conceptual explanations for vanishing results for singular ideals in \cite{CEPSS,NS}.
\setlength{\parindent}{0cm} \setlength{\parskip}{0.5cm}

Next, recall that we say that $C^*_{\ess}(G)$ has the ideal intersection property if for every ideal $I \subseteq C^*_{\ess}(G)$, $I \cap C_0(\Gn) = \gekl{0}$ implies that $I = \gekl{0}$. Similarly, $C^*_r(\tiG_{\ess})$ is said to have the ideal intersection property if for every ideal $\tiI \subseteq C^*_r(\tiG_{\ess})$, $\tiI \cap C_0(\tiG_{\ess}^{(0)}) = \gekl{0}$ implies that $\tiI = \gekl{0}$. 
\begin{introtheorem}[See Theorem~\ref{thm:IIPNonHd}]
\label{thm:the2ndone}
Suppose that $G$ is a non-Hausdorff {\'e}tale groupoid with Hausdorff cover $\tiG$. Assume that $G$ can be covered by countably many open bisections. Then $C^*_{\ess}(G)$ has the ideal intersection property if and only if $C^*_r(\tiG_{\ess})$ has the ideal intersection property.
\end{introtheorem}
\setlength{\parindent}{0cm} \setlength{\parskip}{0cm}

We refer the reader to \S~\ref{s:IIP} for more details. As $\tiG_{\ess}$ is an {\'e}tale groupoid which is Hausdorff, Theorem~\ref{thm:the2ndone} reduces questions about the ideal intersection property for essential groupoid $C^*$-algebras to the Hausdorff case. For example, Theorem~\ref{thm:the2ndone} gives a conceptual explanation for results on simplicity of essential groupoid $C^*$-algebras in \cite{CEPSS,EP22,KM} by showing that they follow from well-known results for Hausdorff groupoids. Theorem~\ref{thm:the2ndone} also significantly simplifies arguments in \cite{KKLRU} required to characterise the ideal intersection property for essential $C^*$-algebras of non-Hausdorff {\'e}tale groupoids and even leads to a more general result (see Remarks~\ref{rem:Improve} and \ref{rem:KKLRU}).
\setlength{\parindent}{0cm} \setlength{\parskip}{0.5cm}

Let us turn to amenability and nuclearity. Recall that an {\'e}tale groupoid $G$ (possibly non-Hausdorff) is called topologically amenable if there exists a net $(\xi_i)$ in $\Cc(G)$ with $\xi_i \geq 0$ for all $i$ such that
\setlength{\parindent}{0cm} \setlength{\parskip}{0cm}

\begin{itemize}
\item $\sum_{g \in G^x} \xi_i(g) \leq 1$ for all $x \in \Gn$ and for all $i$,
\item $\sum_{g \in G^x} \xi_i(g) \to 1$ as $i \to \infty$ uniformly on compact subsets of $\Gn$ (with respect to $x \in \Gu$),
\item $\sum_{h \in G^{r(g)}} \vert \xi_i(g^{-1} h) - \xi_i(h) \vert \to 0$ as $i \to \infty$ uniformly on compact subsets of $G$ (with respect to $g \in G$).
\end{itemize}

Applied to an {\'e}tale groupoid which is Hausdorff, such as the Hausdorff cover $\tiG$ of a non-Hausdorff {\'e}tale groupoid $G$, this yields the classical notion of topological amenability as in \cite{AR} (note that $\Cc$ coincides with $C_c$ in that case).
\begin{introtheorem}[See Theorems~\ref{thm:amenGvstiG}, \ref{thm:nuc=>CtiG-amen} and \ref{thm:nuc<=CtiG-amen}]
\label{thm:the3rdone}
Suppose that $G$ is a non-Hausdorff {\'e}tale groupoid with Hausdorff cover $\tiG$. Then $C^*_r(G)$ is nuclear if and only if $\tiG$ is topologically amenable. If $G$ is topologically amenable, then so is $\tiG$, and the converse holds if $G$ is $\sigma$-compact. Moreover, if the canonical projection $C^*(\tiG) \onto C^*_r(\tiG)$ from the full to the reduced groupoid $C^*$-algebra of $\tiG$ is an isomorphism, then so is the canonical projection $C^*(G) \onto C^*_r(G)$.
\end{introtheorem}
We refer the reader to \S~\ref{s:AmenNuc} for more details. For non-Hausdorff {\'e}tale groupoids $G$, another approach to nuclearity of groupoid $C^*$-algebras attached to $G$ has been developed in \cite{BM23,BM25}. The first claim in Theorem~\ref{thm:the3rdone} can also be deduced from (an upcoming revised version of) \cite{BM25}, even though Hausdorff covers do not appear in \cite{BM25}. Our approach has the benefit of avoiding major technicalities in \cite{BM25}. Theorem~\ref{thm:the3rdone} effectively reduces questions about amenability and nuclearity for non-Hausdorff {\'e}tale groupoids to the Hausdorff case. Our arguments also show that, for a non-Hausdorff {\'e}tale groupoid $G$ which can be covered by countably many open bisections, nuclearity of $C^*_{\ess}(G)$ implies topological amenability of $\tiG_{\ess}$ (see Theorem~\ref{thm:nuc=>CtiG-amen}).
\setlength{\parindent}{0cm} \setlength{\parskip}{0.5cm}

An important step towards our main results is a careful analysis of the Hausdorff cover $\tiG$ of a non-Hausdorff {\'e}tale groupoid $G$. At a technical level, as we explain in Lemma~\ref{lem:DesctiG}, the relationship between $G$ and $\tiG$ is encoded in a canonical translation action of $G$ on $\tiG$ (see \cite[Definition~4.16]{MB} for the notion of translation actions, which are called algebraic morphisms in \cite{BEM}). In a sense, the work in \cite{Tim} was ahead of its time because the construction of essential groupoid $C^*$-algebras was not yet available at the time when \cite{Tim} appeared. An important insight for our work is that the Hausdorff cover offers a new perspective on the essential groupoid $C^*$-algebra, justifying this new construction by showing that it is very natural. In spirit, the idea is very similar to the one behind the construction of essential covers in \cite{BHL}. 
\setlength{\parindent}{0.5cm} \setlength{\parskip}{0cm}

The continuity properties of $\Cc$ functions provided by viewing them as functions on the Hausdorff cover lead to a decomposition result (see Theorem~\ref{thm:CharCinC} and Corollary~\ref{cor:f=sum}) which is crucial in the proof of Theorems~\ref{thm:TheVeryFirstOne}, \ref{thm:Intro_J_CSTAR} and \ref{thm:the1stone}. Theorem~\ref{thm:the2ndone} relies on the existence of a dense set of points where the canonical inclusion $\iota \colon G \into \tiG$ is continuous (and one-to-one). Theorem~\ref{thm:the3rdone} relies on three ingredients: The dual picture for the Hausdorff cover $\tiG$, allowing us to view elements in $C_c(\tiG)$ as bounded Borel functions on the original groupoid $G$, together with the idea of extending representations in \cite[\S~4]{BM25} (\cite{Ren87} and \cite[Appendix~B]{MW08} contain similar ideas, and our improved version appears in Lemma~\ref{lem:ExtRep}), and the insight in \cite{Ren15} that topological amenability is a Borel property.
\setlength{\parindent}{0cm} \setlength{\parskip}{0.5cm}

The Hausdorff cover is not only a natural construction which leads to a new approach to non-Hausdorff {\'e}tale groupoids, but importantly it is also computable. Indeed, by analysing the original non-Hausdorff groupoid, it is possible to describe the Hausdorff cover explicitly, as we demonstrate in a variety of examples arising from group bundles, groupoids of germs, as well as canonical actions of Thompson groups and self-similar groups (see \S~\ref{s:Ex}). In particular, we construct a concrete non-Hausdorff {\'e}tale groupoid $G$ for which there cannot be a $^*$-homomorphism $C^*_r(G) \to C^*_r(\tiG)$ which is injective in $K$-theory (see Example~\ref{ex:NotKInj}). This means that even though the Hausdorff cover is helpful for studying many interesting properties of the underlying non-Hausdorff groupoid, in $K$-theory, we might lose information when going from $G$ to $\tiG$. Note that Example~\ref{ex:NotKInj} shows that $G$ and $\tiG$ (and also $G$ and $\tiG_{\ess}$) are not Morita equivalent in general.

Our paper is structured as follows: \S~\ref{s:Pre} contains required background material and fixes notation. In \S~\ref{s:HdCover} we construct the Hausdorff cover and characterise the images of $\mfi$ and $\mfi_R$ (Theorem~\ref{thm:the1stone}). In \S~\ref{s:SingId}, we prove Theorems~\ref{thm:TheVeryFirstOne} and \ref{thm:Intro_J_CSTAR}, describe singular ideals and establish sufficient criteria for their vanishing. \S~\ref{s:IIP} discusses the ideal intersection property and contains the proof of Theorem~\ref{thm:the2ndone}. In \S~\ref{s:AmenNuc}, we discuss amenability and nuclearity and prove Theorem~\ref{thm:the3rdone}. Finally, we illustrate our main results with the help of examples in \S~\ref{s:Ex}. We discuss examples coming from group bundles and related constructions in \S~\ref{ss:ExGB}, Thompson groups in \S~\ref{ss:ExThompson} and self-similar groups in \S~\ref{ss:ExSelfSim}.

The fourth-named author would like to thank O. Mohsen for making us aware of the reference \cite{Tim}. Initially, we had constructed Hausdorff covers for non-Hausdorff {\'e}tale groupoids, inspired by \cite{BHL}, but then we realised that the construction has already been carried out in \cite{Tim}. We also thank D. Mart{\'i}nez for making us aware of \cite{BM25}. Before that, we had already established our results on amenability and nuclearity under the additional assumption of second countability, using Renault's disintegration theorem (see Remark~\ref{rem:Borel-amen=>nuc}). But then we saw the idea of extending representations in \cite[\S~4]{BM25} and realised that it allows to prove our results in full generality (see Theorem~\ref{thm:nuc<=CtiG-amen}).

Parts of this paper are based on contents of the second-named author's PhD thesis (currently in progress at the University of Glasgow).

\section{Preliminaries}
\label{s:Pre}

We will use the following notations: If $G$ is a topological space and $S \subseteq G$ a subset, we write $S^c \coloneq G \setminus S$, $\inte(S)$ stands for the interior of $S$ and $\overline{S}$ for the closure of $S$. For a function $f$ on $G$ with values in an abelian group with identity element $0$, we let $\osupp(f) \coloneq \menge{ g \in G }{ f(g) \neq 0 }$ be the strict support. We let $B_b(G)$ denote the set of bounded Borel measurable functions $f \colon G \to \Cz$. The set of bounded Borel measurable functions $G \to \Cz$ with strict support contained in a compact subset of $G$ will be denoted by $B_{c}(G)$. The closure of $B_c(G)$ in $B_b(G)$ with respect to the uniform norm $\Vert \cdot \Vert_{\infty}$ will be denoted by $B_0(G)$. Given a locally compact Hausdorff space $\tiG$, we denote by $C_c(\tiG)$ the set of continuous functions $\tiG \to \Cz$ with compact support, and we let $C_0(\tiG)$ denote the set of continuous functions $\tiG \to \Cz$ which vanish at infinity.

The focus of this paper will be groupoids, which are small categories whose morphisms are all invertible. As usual, we identify a groupoid with its set of morphisms $G$, and view its set of objects (also called units) $\Gn$ as a subset of $G$ by identifying objects with the corresponding identity morphisms. By definition, our groupoid $G$ comes with range and source maps $r \colon G \to \Gn$, $s \colon G \to \Gn$, a multiplication map
$G \leftindex_s{\times}_r G = \menge{(g_1,g_2)}{s(g_1) = r(g_2)} \to G, \, (g_1, g_2) \ma g_1 g_2$
and an inversion map $G \to G, \, g \ma g^{-1}$ such that $r(g^{-1}) = s(g)$, $s(g^{-1}) = r(g)$, $g g^{-1} = r(g)$ and $g^{-1} g = s(g)$. These structure maps satisfy a list of conditions so that $G$ becomes a small category (see for instance \cite[Chapter~I, Section~1]{Ren}). We will use the notations $G_x \coloneq s^{-1}(x)$, $G^x \coloneq r^{-1}(x)$ and $G_x^x \coloneq r^{-1}(x) \cap s^{-1}(x)$ for $x \in \Gu$. Moreover, we write $\Isot(G) \coloneq \menge{g \in G}{r(g) = s(g)}$. 
\setlength{\parindent}{0.5cm} \setlength{\parskip}{0cm}

A topological groupoid is a groupoid $G$ which is endowed with a topology such that range, source, multiplication and inversion maps are all continuous. Throughout this paper, a topological groupoid $G$ will always be assumed to be locally compact, and the unit space $\Gn \subseteq G$ will always be assumed to be Hausdorff in the subspace topology. A topological groupoid $G$ is called {\'e}tale if the range map (and hence also the source map) is a local homeomorphism. In this case, the range map (and also the source map) is open, and it follows that $\Gn$ is an open subset of $G$. An open subset $U \subseteq G$ is called an open bisection if the restricted range and source maps $r \vert_U \colon \: U \to r(U), \, g \ma r(g)$ and $s \vert_U \colon U \to s(U), \, g \ma s(g)$ are bijections. If $G$ is {\'e}tale, then $G$ has a basis for its topology consisting of open bisections. A topological groupoid $G$ is called ample if it is {\'e}tale and its unit space $\Gn$ is totally disconnected. An {\'e}tale groupoid is ample if and only if it has a basis for its topology consisting of compact open bisections. Let $G$ be an {\'e}tale groupoid, and let $\B$ denote the set of open bisections in $G$. Every $U \in \B$ is locally compact and Hausdorff with respect to the subspace topology coming from $G$, and $G = \bigcup_{U \in \B} U$. 
\setlength{\parindent}{0.5cm} \setlength{\parskip}{0cm}

Now $\Cc(G) \coloneq \lspan( \bigcup_{U \in \B} C_c(U) ) \subseteq B_c(G)$ becomes a $^*$-algebra with respect to point-wise addition, involution $f^*(g) = \overline{f(g^{-1})}$ and convolution $(f_1 * f_2) (g) = \sum_{(g_1, g_2) \, \in \, G \leftindex_{s}{\times}_r G, \, g_1 g_2 \, = \, g} f_1(g_1) f_2(g_2)$ for $g \in G$. Let us now explain the construction of reduced groupoid $C^*$-algebras. Given $x \in \Gu$, let $\pi_x \colon \Cc(G) \to \cL(\ell^2(G_x))$ denote the left regular representation corresponding to $x$, which is given by $\pi_x(f)(\xi) = f * \xi$ for $f \in \Cc(G)$ and $\xi \in \ell^2(G_x)$, where $(f * \xi)(g) =  \sum_{(g_1, g_2) \, \in \, G \leftindex_{s}{\times}_r G, \, g_1 g_2 \, = \, g} f(g_1) \xi(g_2)$ for $g \in G_x$. Then the reduced groupoid $C^*$-algebra $C^*_r(G)$ of $G$ is given by the completion of $\Cc(G)$ with respect to the norm $\Vert f \Vert_{C^*_r(G)} = \sup_{x \in \Gu} \Vert \pi_x(f) \Vert$. The full (or maximal, or universal) groupoid $C^*$-algebra $C^*(G)$ of $G$ is given by the completion of $\Cc(G)$ with respect to the norm $\Vert f \Vert_{C^*(G)} = \sup_{\pi} \Vert \pi(f) \Vert$, where the supremum is taken over all $^*$-algebra representations $\pi \colon \Cc(G) \to \cL(H)$, for some Hilbert space $H$. 
\setlength{\parindent}{0.5cm} \setlength{\parskip}{0cm}

If $G$ is an ample groupoid and $R$ is a ring, then $RG \coloneq \lspan( \menge{c_U}{c \in R, \, U \subseteq G \text{ compact open bisection}}$ is the associated Steinberg algebra, where $c_U$ denotes the function $G \to R$ which takes the value $c$ on $U$ and $0$ everywhere else, and we take the linear span in the $R$-module of functions $G \to R$. The product in $R G$ is given by convolution (the product of $c_U$ and $d_V$ is given by $(cd)_{UV}$). Note that the coefficient ring $R$ is often assumed to be commutative and unital, which is not necessary for our purposes.
\setlength{\parindent}{0cm} \setlength{\parskip}{0.5cm}

The following observation will be useful later on: Let $G$ be an {\'e}tale groupoid. Given $x \in \Gu$ and $f \in \Cc(G)$, define $\lambda_x(f) \coloneq \sum_{g \in G^x} f(g)$. 
\blemma
\label{lem:LambdaCont}
For every $f \in \Cc(G)$, the map $X \to \Cz, \, x \ma \lambda_x(f)$ is continuous.
\elemma
\setlength{\parindent}{0cm} \setlength{\parskip}{0cm}

\bproof
It suffices to prove the claim for $f \in C_c(U)$ for some open bisection $U \in \B$. For such a function $f$, we have that $\lambda_x(f) = f( (r \vert_U)^{-1}(x) )$ if $x \in r(U)$ and $\lambda_x(f) = 0$ if $x \notin r(U)$. As $r \vert_U$ is a homeomorphism $U \congto r(U)$, it follows that $r(U) \to \Cz, \, x \ma \lambda_x(f)$ lies in $C_c(r(U))$. This proves our claim.
\eproof
\setlength{\parindent}{0cm} \setlength{\parskip}{0.5cm}

\bremark
\label{rem:LambdaContFors}
Similarly, we have that, for every $f \in \Cc(G)$, the map $X \to \Cz, \, x \ma \sum_{g \in G_x} f(g)$ is continuous.
\eremark

\section{Construction of the Hausdorff cover groupoid}
\label{s:HdCover}

Let $G$ be a non-Hausdorff {\'e}tale groupoid. In the following, we describe the Hausdorff cover $\ti{G}$ of $G$, which has been constructed in \cite{Tim} (related constructions appeared in \cite{KS,Tu}). 

\bdefin
Let $\tiG$ be the closure of $G$ in $\fC(G) \setminus \gekl{\emptyset}$ in the Fell topology, where $\fC(G)$ is the space of closed subsets of $G$.
\edefin
By construction, $\tiG$ is Hausdorff. Elements of $\tiG$ are given by non-empty subsets $\bmg \subseteq G$ for which there exist nets $(g_i)$ in $G$ with the property that $\bmg$ is the set of all limit points of $(g_i)$ and that all cluster points of $(g_i)$ are limit points. Here, $g \in G$ is a limit point of $(g_i)$ if for every open set $U \subseteq G$ with $g \in U$, there exists $i_0$ such that $g_i \in U$ for all $i \geq i_0$, and cluster points of $(g_i)$ are precisely the limit points of subnets of $(g_i)$. A basis for the topology of $\tiG$ is given by sets of the form $\cU(U_i,K) \coloneq \menge{\bmg \in \tiG}{\bmg \cap U_i \neq \emptyset \ \forall \ 1 \leq i \leq n, \, \bmg \cap K = \emptyset}$, where the $U_i$ run through all open subsets of $G$ and $K$ runs through all compact subsets of $G$.

Let us present an alternative, dual picture for $\tiG$. Let $C^*(\Cc(G)) \subseteq B_0(G)$ be the smallest sub-$C^*$-algebra of $B_0(G)$ containing $\Cc(G)$ (involution and multiplication in $B_0(G)$ are defined point-wise here). 
\blemma
\label{lem:tiG=Spec}
The map $\tiG \to \Spec (C^*(\Cc(G)))$ sending $\bmg \in \tiG$ to the character $\chi_{\bmg}$ on $C^*(\Cc(G))$ determined by $\chi_{\bmg}(f) = f(g)$ if $f \in C_c(U)$ and $g \in \bmg \cap U$, for $U \in \B$, is a homeomorphism. Its inverse sends a character $\chi \in \Spec (C^*(\Cc(G)))$ to
$
 \bmg_{\chi} \coloneq \menge{g \in G}{\text{There is} \ U \in \B \text{ with } g \in U \text{ and } \chi \vert_{C_c(U)} = \ev_g}
$,
where $\ev_g$ is the character given by $C_c(U) \to \Cz, \, f \ma f(g)$.

The dual isomorphism $C_0(\tiG) \cong C^*(\Cc(G))$ yields the embedding $\mfi \colon \Cc(G) \into C^*(\Cc(G)) \congto C_0(\tiG)$, and we have 
\begin{equation}
\label{e:dec}
 \mfi(f)(\bmg) = \sum_{g \, \in \, \bmg} f(g)
\end{equation}
for all $\bmg \in \tiG$.
\elemma
\setlength{\parindent}{0cm} \setlength{\parskip}{0cm}

\bproof
Continuity of $\bmg \ma \chi_{\bmg}$ is straightforward. Given $\chi \in \Spec(C^*(\Cc(G)))$, the fact that $g_{\chi}$ is in $\tiG$ follows from density of $\menge{\ev_g}{g \in G}$ in $\Spec(C^*(\Cc(G)))$. Continuity of $\chi \ma \bmg_{\chi}$ is straightforward once one notices that the basic open sets $\cU(U_i,K)$, where the $U_i$ are open bisections and $K$ is a finite union of compact bisections each
contained in an open bisection, form a basis for the topology on $\tiG$. To see \eqref{e:dec}, observe that we have, for all $U \in \B$, $f \in C_c(U)$ and $\bmg \in \tiG$ that $\chi_{\bmg}(f) = \sum_{g \, \in \, \bmg} f(g)$.
\eproof
\setlength{\parindent}{0cm} \setlength{\parskip}{0.5cm}

In the following, we will make use of the above identification $\tiG \cong \Spec (C^*(\Cc(G)))$ without explicitly mentioning it. Let $\iota \colon G \into \tiG$ be the canonical inclusion sending an element $g$ to the singleton set $\gekl{g} \in \tiG$, so that the isomorphism $C_0(\tiG) \cong C^*(\Cc(G))$ dual to the homeomorphism $\tiG \cong \Spec (C^*(\Cc(G)))$ from Lemma~\ref{lem:tiG=Spec} is witnessed by the map $\tif \mapsto \tif \circ \iota$.

Because $G$ is a groupoid, elements of $\tiG$ (and also $\tiG$ itself) have more structure: Choose an element $\dot{g} \in \bmg \in \ti{G}$. Then it is easy to see that $r(g') = r(\dot{g})$ for all $g' \in \bmg$, and similarly $s(g') = s(\dot{g})$ for all $g' \in \bmg$. Moreover, given a net $(g_i)$ in $G$ with $\bmg = \lim g_i$, the net $(s(g_i))$ converges to $\ti{s}(\bmg) \coloneq \menge{g_1^{-1} g_2}{g_1, g_2 \in \bmg}$, and $\ti{s}(\bmg)$ is a subgroup of $G_{s(\dot{g})}^{s(\dot{g})}$. Similarly, $(r(g_i))$ converges to $\ti{r}(\bmg) \coloneq \menge{g_1 g_2^{-1}}{g_1, g_2 \in \bmg}$, and $\ti{r}(\bmg)$ is a subgroup of $G_{r(\dot{g})}^{r(\dot{g})}$. Our element $\bmg$ can then be written as $\bmg = \dot{g} \cdot \ti{s}(\bmg) = \ti{r}(\bmg) \cdot \dot{g}$. Here and in the sequel, we denote by $g \cdot \bmx$ the product of $\gekl{g}$ and $\bmx$ as subsets of $G$, i.e., $g \cdot \bmx = \menge{g g_x}{g_x \in \bmx}$. Moreover, let $\ti{G}^{(0)}$ denote the closure of $G^{(0)}$ in $\ti{G}$. Then $\ti{r}$ and $\ti{s}$ define range and source maps $\ti{G} \to \ti{G}^{(0)}$. Note that $\tis^{-1}(\iota(\Gu)) = \iota(G)$ and $\tir^{-1}(\iota(\Gu)) = \iota(G)$. Given two elements $\bmg, \bmh \in \ti{G}$, the product $\bmg \bmh$ is only defined if $\bmg = g \cdot \bmx$ and $\bmh = \bmx \cdot h$, for some composable $g, h \in G$ and $\bmx \in \tiGu$, in which case $\bmg \bmh = (gh) \cdot \bmy$ where $\bmy$ is such that $\bmh = \bmx \cdot h = h \cdot \bmy$, i.e., $\bmy = h^{-1} \cdot \bmx \cdot h$. Moreover, the inverse of $\bmg = g \cdot \bmx$ is given by $\bmg^{-1} = \bmx \cdot g^{-1} = g^{-1} \cdot (g \cdot \bmx \cdot g^{-1})$. In this way, $\ti{G}$ becomes a locally compact groupoid. 

Moreover, let us explain why $\ti{G}$ is {\'e}tale if $G$ is {\'e}tale. Given $U \in \B$, consider the ideal $C^*(\Cc(G))_U \coloneq \overline{ \lspan(C_c(U) \cdot C^*(\Cc(G))) } \subseteq C^*(\Cc(G))$ and the open subset $\tiU$ of $\tiG$ corresponding to $\Spec( C^*(\Cc(G))_U ) \subseteq \Spec( C^*(\Cc(G)))$ under the homeomorphism $\tiG \cong \Spec( C^*(\Cc(G)))$ from Lemma~\ref{lem:tiG=Spec}. By construction, $\tiU$ is an open subset of $\tiG$, and concretely we have $\tiU = \menge{\bmg \in \tiG}{\bmg \cap U \neq \emptyset}$. The inclusion $C_0(U) \subseteq C^*(\Cc(G))_U$ induces a proper continuous surjection $\pi_U \colon \tiU \onto U$. Moreover, the inclusion $\iota \colon G \into \tiG$ restricts to $U \into \tiU$. Now it follows that $\tir(\tiU)$ is open and that $\tir$ restricts to a homeomorphism $\tiU \cong \tir(\tiU)$. The reason is that $\tir(\tiU) = \widetilde{r(U)}$ is canonically homeomorphic to $\Spec (C^*(\Cc(G))_{r(U)})$, and that the restriction of $\tir$, $\tir_{\tiU} \colon \tiU \congto \tir(\tiU) = \widetilde{r(U)}$, is dual under Gelfand transform to the isomorphism $C^*(\Cc(G))_{r(U)} \cong C^*(\Cc(G))_{U}, \, f \ma f \circ r$. It follows that $\tiU$ is an open bisection. Moreover, we have $\tiG = \bigcup_{U \in \B} \ti{U}$. It also follows that $\tiG$ is {\'e}tale. In particular, continuity of the product map follows from $\tiU \cdot \ti{V}  = \widetilde{UV}$ and the properties of $\tiU$, $\ti{V}$ mentioned above. Note that, with the notation we just introduced, $\tiGu$ is canonically homeomorphic to $\Spec(C^*(\Cc(G))_{\Gu})$.

From now on, we write $X \coloneq G^{(0)}$ and $\ti{X} \coloneq \ti{G}^{(0)}$. Let $\prod (\C_c(G)) \coloneq \menge{f_1 \dotsm f_n}{n \in \Nz, \, f_m \in \C_c(G) \ \forall \ 1 \leq m \leq n}$. 
The following is an immediate consequence of our construction.
\blemma
\label{lem:SpanProd}
We have $C^*(\Cc(G)) = \clspan( \prod (\C_c(G)) )$, where the closure is taken in $B_0(G)$ with respect to $\norm{\cdot}_{\infty}$.
\elemma

\bcor
\label{cor:dim=0}
If $X$ is totally disconnected, then so is $\ti{X}$.
\ecor
\setlength{\parindent}{0cm} \setlength{\parskip}{0cm}

\bproof
Since $X$ is totally disconnected and $G$ is {\'e}tale, $\lspan(\menge{1_U}{U \subseteq G \text{ compact open bisection}})$ is dense in $\C_c(G)$. Hence Lemma~\ref{lem:SpanProd} implies that
\[
 C^*(\Cc(G)) = \clspan( \menge{1_{U_1} \dotsm 1_{U_n}}{n \in \Nz, \, U_m \subseteq G \text{ compact open bisection} \ \forall \ 1 \leq m \leq n} ).
\]
In particular, this shows that $C_0(\ti{X}) \cong C^*(\Cc(G))_X$ is given by
\[
 \clspan( \menge{1_{U_0} \cdot 1_{U_1} \dotsm 1_{U_n}}{n \in \Nz, \, U_m \subseteq G \text{ compact open bisection} \ \forall \ 0 \leq m \leq n, \, U_0 \subseteq X} )
\]
and thus contains a dense $^*$-algebra of projections. This implies that $\ti{X}$ is totally disconnected.
\eproof
\setlength{\parindent}{0cm} \setlength{\parskip}{0.5cm}

The following observation reduces the task of describing the Hausdorff cover $\ti{G}$ to the description of $\ti{X}$ and the projection $\pi \colon \ti{X} \onto X$, which sends $\bmx \in \ti{X}$ to $\pi(\bmx) \coloneq r(g) = s(g)$, where $g \in \bmx$ is arbitrary. Alternatively, we have $\gekl{\pi(\bmx)} = \bmx \cap X$. Note that $\pi$ coincides with $\pi_X$ as introduced above. Now set $G \leftindex_s{\times}_{\pi} \ti{X} \coloneq \menge{(g,\bmx) \in G \times \ti{X}}{s(g) = \pi(\bmx)}$, and define an equivalence relation $\sim$ on $G \leftindex_s{\times}_{\pi} \ti{X}$ by setting $(g,\bmx) \sim (g',\bmx')$ if and only if $\bmx = \bmx'$ and $g^{-1} g' \in \bmx$. Note that $(g,\bmx) \sim (g',\bmx)$ if and only if $g \cdot \bmx = g' \cdot \bmx$. We write $[g,\bmx]$ for the equivalence class of $(g,\bmx) \in G \leftindex_s{\times}_{\pi} \ti{X}$.
\blemma
\label{lem:DesctiG}
The map $q \colon G \leftindex_s{\times}_{\pi} \ti{X} \onto \tiG, \, (g,\bmx) \ma g \cdot \bmx$ is surjective, open and a local homeomorphism. Moreover, $G$ acts on $\tiG$ via $g.\bmg \coloneq q(g,\tir(\bmg)) \cdot \bmg$ if $s(g) = \pi(\tir(\bmg))$, and this is a translation action $G \acts \tiG$ in the sense of \cite[Definition~4.16]{MB} (i.e., $G \acts \tiG$ commutes with the canonical action $\tiG \stca \tiG$ by right translations). The map $q$ induces a homeomorphism $\ti{q} \colon (G \leftindex_s{\times}_{\pi} \ti{X}) / { }_{\sim} \congto \tiG$. 
\setlength{\parindent}{0.5cm} \setlength{\parskip}{0cm}

For every $U \in \B$, the open set $\ti{U} \subseteq \tiG$ constructed above satisfies $\ti{U} = \ti{q}(\menge{[g,\bmx]}{g \in U, \, \bmx \in \pi^{-1}(s(U))})$.
\elemma
\setlength{\parindent}{0cm} \setlength{\parskip}{0cm}

Note that part of this appeared in \cite[Proposition~7]{Tim}. Translation actions are called algebraic morphisms in \cite{BEM}.
\bproof
Clearly, $q$ is surjective. Let us show that $q$ is continuous. Given a net $(g_i)$ in $G$ converging to $g \in G$ and a net $(\bmx_i)$ in $\tiX$ converging to $\bmx \in \tiX$ with $s(g_i) = \pi(\bmx_i)$ and $s(g) = \pi(\bmx)$, there exists $U \in \B$ such that $g \in U$, and, without loss of generality, $g_i \in U$ for all $i$, and thus $\bmx \in \pi^{-1}(s(U))$ as well as $\bmx_i \in \pi^{-1}(s(U))$ for all $i$. So we have $g_i \cdot \bmx_i, g \cdot \bmx \in \tiU$. The homeomorphism $\tiU \cong \tis(\tiU)$ given by $\tis$ sends $g_i \cdot \bmx_i$ to $\bmx_i$ and $g \cdot \bmx$ to $\bmx$. As $(\bmx_i)$ converges to $\bmx$ in $\tis(\tiU)$, it follows that $(g_i \cdot \bmx_i)$ converges to $g \cdot \bmx$ in $\tiU$, as desired. Now it is straightforward to see that the action $G \acts \tiG$ given by $g.\bmg \coloneq q(g,\tir(\bmg)) \bmg = (g \cdot \tir(\bmg)) \bmg$ is well-defined and indeed a translation action. Thus \cite[Lemma~4.20]{MB} implies that $q$ is open and a local homeomorphism. Furthermore, it is clear that $\ti{q}$ is bijective, and it is a homeomorphism because $q$ is continuous and open. The remaining claims are straightforward.
\eproof
\setlength{\parindent}{0cm} \setlength{\parskip}{0.5cm}

\bremark
\label{rem:DesctiG}
Lemma~\ref{lem:DesctiG} tells us that to describe $\tiG$, it suffices to describe $\tiX$ as a topological space and the canonical projection $\pi \colon \tiX \onto X$. We will follow this approach when discussing examples in \S~\ref{s:Ex}.
\eremark

Let us study compact sets in $G$ and $\ti{G}$, as this will be useful later on.

\blemma
\label{lem:cpctGvstiG}
We have the following relationship between compact sets in $G$ and $\ti{G}$:
\setlength{\parindent}{0cm} \setlength{\parskip}{0cm}

\begin{enumerate}[label=(\roman*)]
\item For every compact subset $L \subseteq \ti{X}$ there exists a compact subset $K \subseteq X$ such that $L \subseteq \overline{\iota(K)}$.
\item For every compact subset $L \subseteq \ti{G}$ there exists a compact subset $K \subseteq G$ such that $L \subseteq \overline{\iota(K)}$, and therefore $L \cap \iota(G) \subseteq \iota(K)$.
\item For every compact subset $K \subseteq G$ there exists a compact subset $L \subseteq \ti{G}$ such that $\iota(K) \subseteq L$.
\end{enumerate}
\elemma
\setlength{\parindent}{0cm} \setlength{\parskip}{0cm}

\bproof
(i): Take a compact subset $L \subseteq \ti{X}$. Then $\pi(L)$ is a compact subset of $X$, and we can choose a compact subset $K \subseteq X$ such that $\pi(L) \subseteq \inte(K)$. Now we claim that $L \subseteq \overline{\iota(K)}$. Indeed, take $\bmx \in L$. As $\iota(X)$ is dense in $\ti{X}$, we can find a net $(x_j)$ in $X$ such that $\iota(x_j)$ converges to $\bmx$. Then $x_j = \pi(\iota(x_j))$ converges to $\pi(x)$. Hence $x_j$ lies eventually in $K$, and thus $\bmx = \lim_j \iota(x_j) \in \overline{\iota(K)}$.
\setlength{\parindent}{0cm} \setlength{\parskip}{0.5cm}

(ii): Take a compact subset $L \subseteq \ti{G}$. By compactness of $L$, we can find a finite cover $L \subseteq \bigcup_{k \in F} V_k$ for some compact subsets $V_k \subseteq \ti{U}_k$ and some $U_k \in \B$. Now take compact subsets $K_k \subseteq U_k$ with $\pi_{U_k}(V_k) \subseteq \inte(K_k)$. The same argument as in (i) shows that $V_k \subseteq \overline{\iota(K_k)}$. Hence $L \subseteq \bigcup_{k \in F} V_k \subseteq \bigcup_{k \in F} \overline{\iota(K_k)} \subseteq \overline{\iota(\bigcup_{k \in F} K_k)}$. Now set $K \coloneq \bigcup_{k \in F} K_k$.
\setlength{\parindent}{0.5cm} \setlength{\parskip}{0cm}

For the second part, let $g \in G$ be such that $\iota(g) \in L$. There exists a net $(g_i)$ in $K$ such that $\iota(g_i)$ converges to $\iota(g)$. By definition of the Fell topology, this means that the net $(g_i)$ converges to $g$ in the topology of $G$, and $g$ is the only cluster point. Compactness of $K$ implies that $(g_i)$ has a subnet converging in $K$. This limit must be $g$, so that $g \in K$.
\setlength{\parindent}{0cm} \setlength{\parskip}{0.5cm}

(iii): Take a compact subset $K \subseteq G$. Find a finite cover $K \subseteq \bigcup_{k \in F} K_k$, where $K_k$ are compact sets with $K_k \subseteq U_k$ for some $U_k \in \B$. Set $L_k \coloneq \pi_{U_k}^{-1}(K_k)$. Then $L_k$ is compact as $\pi_{U_k}$ is proper. Moreover, $K_k = \pi_{U_k}(\iota(K_k))$ implies that $\iota(K_k) \subseteq L_k$. It follows that $\iota(K) \subseteq \bigcup_{k \in F} \iota(K_k) \subseteq \bigcup_{k \in F} L_k$. Now set $L \coloneq \bigcup_{k \in F} L_k$.
\eproof
\setlength{\parindent}{0cm} \setlength{\parskip}{0.5cm}

\blemma
\label{lem:C*GintiG}
The canonical embedding $\mfi \colon \Cc(G) \into C_0(\tiG)$ from Lemma~\ref{lem:tiG=Spec} has image in $C_c(\tiG)$, and the resulting embedding $\Cc(G) \into C_c(\tiG)$, also denoted by $\mfi$, is a $^*$-algebra homomorphism (where multiplication in $\Cc(G)$ and $C_c(\tiG)$ is given by convolution). 
\setlength{\parindent}{0.5cm} \setlength{\parskip}{0cm}

Moreover, the map $\mfi \colon \Cc(G) \into C_c(\tiG)$ induces an injective $^*$-homomorphism $C^*_r(\mfi) \colon C^*_r(G) \into C^*_r(\tiG)$.
\elemma
\setlength{\parindent}{0cm} \setlength{\parskip}{0cm}

Note that $C^*_r(\mfi)$ is the restriction of the isomorphism $C^*(\Cc(G)) \congto C_0(\tiG)$ from Lemma~\ref{lem:C*GintiG} to $C^*_r(G)$.

\bproof
The first claim is an immediate consequence of Lemma~\ref{lem:cpctGvstiG}~(iii) and the definition of the relevant $^*$-algebra structures. Let us prove the second claim. Given $x \in \Gu$, let $\pi_x \colon \Cc(G) \to \cL(\ell^2(G_x))$ denote the left regular representation corresponding to $x$, and for $\bmx \in \tiGu$, let $\ti{\pi}_{\bmx} \colon C_c(\tiG) \to \cL(\ell^2(\tiG_{\bmx}))$ denote the left regular representation corresponding to $\bmx$. For every $x \in \Gu$, the representation $\pi_x$ is canonically unitarily equivalent to the representation given by $\Cc(G) \to C_c(\tiG) \to \cL(\ell^2(\tiG_{\iota(x)}))$, where the first arrow is given by $\mfi$ and the second arrow is given by $\ti{\pi}_{\iota(x)}$. This means that we have $\Vert \pi_x(f) \Vert = \Vert \ti{\pi}_{\iota(x)}(\mfi(f) )\Vert$ for every $x \in \Gu$ and $f \in \Cc(G)$. We conclude that
\[
 \Vert f \Vert_{C^*_r(G)} = \sup_{x \in \Gu} \Vert \pi_x(f) \Vert = \sup_{x \in \Gu} \Vert \ti{\pi}_{\iota(x)}(\mfi(f)) \Vert = \Vert \mfi(f) \Vert_{C^*_r(\tiG)}.
\]
For the last equality, we used that $\iota(\Gu)$ is dense in $\tiGu$ and that, since $\tiG$ is Hausdorff, taking the supremum over a dense subset of the unit space is enough to obtain the norm in $C^*_r(\tiG)$.
\eproof
\setlength{\parindent}{0cm} \setlength{\parskip}{0.5cm}

\bremark
\label{rem:SteinbIncl}
Lemma~\ref{lem:C*GintiG} has the following algebraic analogue: Let $G$ be an ample groupoid. Its Hausdorff cover $\tiG$ is again ample by Corollary~\ref{cor:dim=0}. Now let $R$ be a ring and form the Steinberg algebras $RG$ and $R \tiG$. Then we have a canonical inclusion of $R$-algebras $\mfi_R \colon RG \into R \tiG$ given by $c_U \ma c_{\ti{U}}$.
\eremark

\bremark
As we will see later (see Corollary~\ref{cor:FullC*Embed}), $\mfi$ also induces an injective $^*$-homomorphism at the level of full groupoid $C^*$-algebras.
\eremark

Let us introduce the following terminology, which is inspired by \eqref{e:dec}.
\bdefin
\label{def:dec}
Let $G$ be a non-Hausdorff {\'e}tale groupoid and $\tiG$ its Hausdorff cover. Let $R$ be a ring. A function $\ti{f} \colon \tiG \to R$ is called decomposable if for all $\bmg \in \tiG$, we have
\begin{equation}
\label{e:CcGinCctiG}
 \ti{f}(\bmg) = \sum_{g \in \bmg} \ti{f}(\iota(g)).
\end{equation}
\edefin
If a function $\tif$ is decomposable, then we mean that only finitely many summands on the right-hand side of \eqref{e:CcGinCctiG} are non-zero for $\tif$. Note that, if $\tif$ is decomposable, then we have $\ti{f}(\bmg) = \sum_i \ti{f}(\bmh_i)$ for all $\bmg \in \tiG$ and every family $\bmh_i \in \tiG$ with $\bmg = \coprod_i \bmh_i$ (as subsets of $G$).

The following is an immediate consequence of \eqref{e:dec} because $\norm{\cdot}_{C^*_r}$ dominates $\norm{\cdot}_{\infty}$.
\setlength{\parindent}{0cm} \setlength{\parskip}{0cm}

\blemma
\label{lem:CinC}
Let $G$ be a non-Hausdorff {\'e}tale groupoid and $\tiG$ its Hausdorff cover. Suppose that $\bmg \in \tiG$ satisfies $\# \, \bmg < \infty$. If $\tif \in C^*_r(\tiG)$ lies in $\im(C^*_r(\mfi))$, then $\tif(\bmg) = \sum_{g \in \bmg} \tif(\iota(g))$. In particular, if $\# \, \bmg < \infty$ for all $\bmg \in \tiG$, then every $\tif \in C^*_r(\tiG)$ which lies in $\im(C^*_r(\mfi))$ is decomposable (here we view $\tif$ as a function $\tiG \to \Cz$).
\elemma
\setlength{\parindent}{0cm} \setlength{\parskip}{0.5cm}

Let us now characterise the image $\im(\mfi)$ of $\mfi$ and, for ample groupoids $G$, the image $\im(\mfi_R)$ of $\mfi_R$.
\btheo
\label{thm:CharCinC}
Let $G$ be a non-Hausdorff {\'e}tale groupoid and $\tiG$ its Hausdorff cover.
\setlength{\parindent}{0cm} \setlength{\parskip}{0cm}

\begin{enumerate}[label=(\roman*)]
    \item Let $\ti{f} \in C_c(\tiG)$ be decomposable and set $f \coloneq \tif \circ \iota$. Assume that there exist $V_1, \dotsc, V_n \in \B$ and a compact set $K \subseteq V_1 \cup \dotso \cup V_n$ such that $\osupp(f) \subseteq K$. Then there are $f_1, \dotsc, f_n$ with $f_i \in C_c(V_i)$ for all $1 \leq i \leq n$ such that $f = f_1 + \dotso + f_n$.
    \item A function $\ti{f} \in C_c(\tiG)$ lies in $\im(\mfi)$ if and only if it is decomposable.
    \item Assume that $G$ is ample, and let $R$ be a ring. Let $\ti{f} \in R \tiG$ be decomposable and set $f \coloneq \tif \circ \iota$. Assume that there exist $V_1, \dotsc, V_n \in \B$ and a compact set $K \subseteq V_1 \cup \dotso \cup V_n$ such that $\osupp(f) \subseteq K$. Then there are $f_1, \dotsc, f_n \in RG$ with $f_i \in C_c(V_i,R)$ for all $1 \leq i \leq n$ such that $f = f_1 + \dotso + f_n$.    
    \item Assume that $G$ is ample, and let $R$ be a ring. A function $\ti{f} \in R \tiG$ lies in $\im(\mfi_R)$ if and only if it is decomposable.
\end{enumerate}
\etheo
\setlength{\parindent}{0cm} \setlength{\parskip}{0cm}

\bproof
Let us prove (i). Since $K$ is compact, there are open subsets $U_i \subseteq V_i$ such that the closure $\overline{U_i}^{V_i}$ of $U_i$ in $V_i$ is compact (for all $1 \leq i \leq n$) and we have $\osupp(f) \subseteq U_1 \cup \dotso \cup U_n$. 
\setlength{\parindent}{0.5cm} \setlength{\parskip}{0cm}

Let us now construct, inductively on $m \in \gekl{1, \dotsc, n-1}$, functions $f_i \in C_c(V_i)$, for $1 \leq i \leq m$, such that $\osupp(f - f_1 - \dotso - f_m) \subseteq U_{m+1} \cup \dotso \cup U_n$. Let $\check{V}_1 \coloneq V_1 \setminus (U_2 \cup \dotso \cup U_n)$. Our first claim is that $f$ is continuous on $\check{V}_1$. To see this, take a net $(g_k)$ in $\check{V}_1$ converging to $g \in \check{V}_1$. By properness of $\pi_{V_1}$, we may without loss of generality assume that $\iota(g_k)$ converges in $\tiG$, say to $\bmg \in \tiG$. Therefore, $f(g_k) = \ti{f}(\iota(g_k))$ converges to $\ti{f}(\bmg)$ by continuity of $\ti{f}$ on $\tiG$. As $\ti{f}$ is decomposable, we have $\ti{f}(\bmg) = \sum_{g' \, \in \, \bmg} \ti{f}(\iota(g')) = \sum_{g' \, \in \, \bmg} f(g')$. But for every $g' \in \bmg$, we have $g' \notin U_2 \cup \dotso \cup U_n$ since $g_k \notin U_2 \cup \dotso \cup U_n$ for all $k$. Since $\osupp(f) \subseteq U_1 \cup \dotso \cup U_n$, we conclude that $f(g') = 0$ for all $g' \in \bmg$ with $g' \neq g$. It follows that $\sum_{g' \, \in \, \bmg} f(g') = f(g)$ and hence $f(g_k)$ converges to $f(g)$, as desired.

Now take a compact subset $K_1 \subseteq V_1$ such that $\overline{U_1}^{V_1} \subseteq \inte(K_1)$. As $f$ is continuous on the closed subset $\check{K}_1 \coloneq K_1 \setminus (U_2 \cup \dotso \cup U_n)$, the Tietze extension theorem allows us to find a continuous function $f' \colon K_1 \to \Cz$ extending $f$ on $\check{K}_1$, i.e., $f' \equiv f$ on $\check{K_1}$. Let $f'' \in C_c(V_1)$ be such that $f'' \equiv 1$ on $\overline{U_1}^{V_1}$ and $f'' \equiv 0$ on $V_1 \setminus K_1$. Let $f_1 \colon V_1 \to \Cz$ be such that $f_1 \equiv f' \cdot f''$ on $K_1$ and $f_1 \equiv 0$ on $V_1 \setminus K_1$. Clearly, $f_1 \in C_c(V_1)$ and $f_1 \equiv f$ on $\check{V}_1$.

Our second claim is that $\osupp(f - f_1) \subseteq U_2 \cup \dotso \cup U_n$. Indeed, let $g \notin U_2 \cup \dotso \cup U_n$. If $g \in \check{V}_1$, then $f(g) = f_1(g)$ by construction, so that $(f-f_1)(g) = 0$. Now assume that $g \notin \check{V}_1$. Then $g \notin U_1 \cup \dotso \cup U_n$ and hence $f(g) = 0$. At the same time, $g \notin V_1$ and hence $f_1(g) = 0$. We conclude that $(f-f_1)(g) = 0$, as desired. Moreover, with $\ti{f}_1 \coloneq \mfi(f_1)$, we have that $\ti{f} - \ti{f}_1$ is decomposable and $f-f_1 = (\ti{f} - \ti{f}_1) \circ \iota$.
\setlength{\parindent}{0cm} \setlength{\parskip}{0.5cm}

Proceeding inductively, we obtain $f_1, \dotsc, f_{n-1} \in \Cc(G)$ with $f_i \in C_c(V_i)$ such that, with $f_n \coloneq f - f_1 - \dotso - f_{n-1}$, we have $\osupp(f_n) \subseteq U_n$ and $\ti{f}_n \coloneq \mfi(f_n)$ is decomposable. The above argument (our first claim) now implies that $f_n$ is continuous on $V_n$, and thus we conclude that $f_n \in C_c(V_n)$. Moreover, $f = f_1 + \dotso + f_n$, as desired.

Let us now prove (ii). The direction \an{$\Rarr$} follows from \eqref{e:dec}. To prove \an{$\Larr$}, let $\ti{f} \in C_c(\tiG)$ be decomposable, set $f \coloneq \ti{f} \circ \iota$, and apply (i). We obtain $f_i \in \Cc(G)$ such that $\ti{f} \circ \iota = f = f_1 + \dotso + f_n = \mfi(f_1 + \dotso + f_n) \circ \iota$. It follows that $\ti{f} = \mfi(f_1 + \dotsc + f_n)$, as desired. 

The proof of (iii) is analogous to the one of (i), so we just explain the main difference. Let $G$ be ample, $R$ a ring, and let $\ti{f} \in R G$ be decomposable. The argument for the first claim above shows that $f = \ti{f} \circ \iota$ is continuous on $\check{V}_1$. The main difference lies in the construction of $f_1$. As $V_1$ is totally disconnected, we can find a compact open set $K_1$ with $\overline{U_1}^{V_1} \subseteq K_1 \subseteq V_1$. As above, let $\check{K}_1 \coloneq K_1 \setminus (U_2 \cup \dotso \cup U_n)$. Then $\menge{\check{K}_1 \cap O}{O \subseteq K_1\text{ compact open}}$ is a basis of open sets in $\check{K}_1$. As $f$ is continuous on $\check{K}_1$ and only takes finitely many values, $f$ is locally constant on $\check{K}_1$. Hence there exist $c_1, \dotsc, c_m \in R$ and compact open subsets $O_1, \dotsc, O_m \subseteq K_1$ such that $f \vert_{\check{K}_1} \equiv \sum_{j=1}^m (c_j)_{\check{K}_1 \cap O_j}$. Now define $f_1 \coloneq \sum_{j=1}^m (c_j)_{O_j}$. Clearly, $f_1 \in C_c(V_1,R)$. By construction, $f \equiv f_1$ on $\check{K}_1$. Moreover, $f \equiv 0$ on $\check{V}_1 \setminus \check{K}_1$ and $f_1 \equiv 0$ on $\check{V}_1 \setminus \check{K}_1$. Thus $f \equiv f_1$ on $\check{V}_1$, and the same argument as for the second claim shows that $f - f_1 \equiv 0$ on $G \setminus (U_2 \cup \dotsc \cup U_n)$. Now proceed as above to complete the proof.

The proof of (iv) is analogous to the one of (ii).
\eproof 
\setlength{\parindent}{0cm} \setlength{\parskip}{0cm}

We record the following immediate consequence of Theorem~\ref{thm:CharCinC} for later use.
\bcor
\label{cor:f=sum}
Let $G$ be a non-Hausdorff {\'e}tale groupoid.
\setlength{\parindent}{0cm} \setlength{\parskip}{0cm}

\begin{enumerate}[label=(\roman*)]
    \item Let $f \in \Cc(G)$. Assume that there exist $V_1, \dotsc, V_n \in \B$ and a compact set $K \subseteq V_1 \cup \dotso \cup V_n$ such that $\osupp(f) \subseteq K$. Then there are $f_1, \dotsc, f_n$ with $f_i \in C_c(V_i)$ for all $1 \leq i \leq n$ such that $f = f_1 + \dotso + f_n$.
    \item Assume that $G$ is ample, and let $R$ be a ring. Let $f \in R G$. Assume that there exist $V_1, \dotsc, V_n \in \B$ and a compact set $K \subseteq V_1 \cup \dotso \cup V_n$ such that $\osupp(f) \subseteq K$. Then there are $f_1, \dotsc, f_n \in RG$ with $f_i \in C_c(V_i,R)$ for all $1 \leq i \leq n$ such that $f = f_1 + \dotso + f_n$.
\end{enumerate}
\ecor
Note that even though the statement of Corollary~\ref{cor:f=sum} does not mention the Hausdorff cover, its proof uses the Hausdorff cover because Theorem~\ref{thm:CharCinC} does.
\setlength{\parindent}{0cm} \setlength{\parskip}{0.5cm}

\section{Vanishing criteria for singular ideals}
\label{s:SingId}

Let $G$ be a non-Hausdorff {\'e}tale groupoid with unit space $X \coloneq \Gn$. Following \cite[\S~7]{KM} (see in particular \cite[Proposition~7.18]{KM}), recall that one way to define the singular ideal $J$ of $C^*_r(G)$ is to set
\[
 J \coloneq \menge{f \in C^*_r(G)}{s(\osupp(f)) \text{ has empty interior}}.
\]
If $G$ is ample and $R$ is a ring, then the singular ideal $J_R$ of $RG$ is defined as
\[
 J_R \coloneq \menge{f \in RG}{\osupp(f) \text{ has empty interior}}
\]
(see \cite{CEPSS,StSz21,StSz23}). Let us first reconcile these two seemingly different definitions of $C^*$-algebraic and algebraic singular ideals. In the following, given $f \in C^*_r(G)$, we write
$
S_\varepsilon(f) \coloneq \menge{g \in G}{\vert f(g) \vert > \varepsilon}
$.
\blemma
\label{lem:Sing=Sing}
Let $G$ be a non-Hausdorff {\'e}tale groupoid. 
\setlength{\parindent}{0cm} \setlength{\parskip}{0cm}

\begin{enumerate}[label=(\roman*)]
    \item Let $f \in \Cc(G)$ and $\varepsilon>0$ be such that $S_\varepsilon(f)$ has empty interior. Then $s(S_{\varepsilon}(f))$ is nowhere dense in $X$.
    \item If $f \in J$, then for all $\varepsilon > 0$, we have that $s(S_{\varepsilon}(f))$ is nowhere dense in $X$.
    \item An element $f \in C^*_r(G)$ lies in $J$ if and only if $\osupp(f)$ has empty interior.
    \item Let $G$ be ample and $R$ a ring. An element $f \in RG$ lies in $J_R$ if and only if $s(\osupp(f))$ has empty interior.
\end{enumerate}
\elemma
\setlength{\parindent}{0cm} \setlength{\parskip}{0cm}

\bproof
Let us prove (i). Assume for a contradiction that this is not the case, i.e., there is a non-empty open set $U \subseteq \overline{s(S_\varepsilon(f))}^X$. Choose the smallest positive integer $n$ such that there exist open bisections $U_1, \dotsc, U_n \subseteq G$, $f_i \in C_c(U_i)$, and a function $\chi \in C_c(U)$ satisfying $0 \le \chi \le 1$ and $\chi \equiv 1$ on some non-empty open subset $V \subseteq U$, such that $f * \chi = \sum_{i=1}^n f_i$. Note that $f * \chi \neq 0$. Set $S_i \coloneq \osupp(f_i)$. Our first claim is that, by minimality of $n$, we have $\vert (f * \chi) (g) \vert \le \varepsilon$ whenever $g \in S_i \cap S_j$ for $i \neq j$. Indeed, suppose that there are $i \neq j$ and $g \in S_i \cap S_j$ with $\vert (f * \chi) (g) \vert > \varepsilon$. Then $\vert f(g) \vert > \varepsilon$ and $\chi(s(g)) > 0$. Choose $\chi' \in C_c(s(S_i \cap S_j) \cap \osupp(\chi))$ such that $0 \leq \chi' \leq 1$ and $\chi' \equiv 1$ on some open set $W$ containing $s(g)$. Then $f * \chi' = (f * \chi) * (\chi^{-1} \chi') = \sum_{i=1}^n f_i * (\chi^{-1} \chi')$. But both $f_i * (\chi^{-1} \chi')$ and $f_j * (\chi^{-1} \chi')$ lie in $C_c(S_i \cap S_j)$, so that we reduced $n$, which is impossible. This proves our first claim. Our second claim is that $s(S_\varepsilon(f*\chi))$ is nowhere dense in $X$. For each $i$, set $W_i \coloneq S_\varepsilon(f_i) \cap (\bigcup_{j \neq i} S_j)$. Since $S_\varepsilon(f)$ has empty interior and $S_\varepsilon(f_i) \setminus W_i \subseteq S_\varepsilon(f)$, and since each $S_\varepsilon(f_i)$ is open, we conclude that $S_\varepsilon(f_i) \subseteq \overline{W_i}^{U_i}$. Hence $S_\varepsilon(f_i) \setminus W_i \subseteq \partial_{U_i}W_i \coloneq \overline{W_i}^{U_i} \setminus W_i$. Since each $W_i$ is open and contained in a compact subset of $U_i$, $\partial_{U_i} W_i$ is compact and has empty interior. Moreover, by our first claim, we have $s(S_\varepsilon(f * \chi)) = \bigcup_{i=1}^n s(S_\varepsilon(f_i) \setminus W_i) \subseteq \bigcup_{i=1}^n s(\partial_{U_i} W_i)$. Thus $s(S_\varepsilon(f * \chi))$ is nowhere dense, because it is contained in a finite union of closed sets with empty interior. This proves our second claim. But $V$ satisfies $V \subseteq U \subseteq \overline{s(S_\varepsilon(f))}^X$ and $V \cap s(S_\varepsilon(f)) = V \cap s(S_\varepsilon(f * \chi))$. Since $V$ is open, we have $V \subseteq \overline{s(S_\varepsilon(f * \chi))}^X$ and therefore $V = \emptyset$ by our second claim, which contradicts the assumption that $V$ was non-empty. Overall, we conclude that $s(S_\varepsilon(f))$ is nowhere dense in $X$. This proves (i). A similar argument (applied to $\osupp(f)$ instead of $S_{\varepsilon}(f)$ for $f \in RG$) yields (iv).
\setlength{\parindent}{0cm} \setlength{\parskip}{0.5cm}

Our third claim is that if $f \in C_r^*(G)$ is such that $\osupp(f)$ has empty interior, then $s(S_{\varepsilon}(f))$ is nowhere dense in $X$ for all $\varepsilon > 0$. Indeed, take $f_{\varepsilon/2} \in \Cc(G)$ with $\norm{f- f_{\varepsilon/2}}_\infty < \varepsilon/2$. Then $S_{\varepsilon/2}(f_{\varepsilon/2}) \subseteq \osupp(f)$ has empty interior. By (i), $s(S_{\varepsilon/2}({f_{\varepsilon/2}}))$ is nowhere dense in $X$. However, $S_{\varepsilon}(f) \subseteq S_{\varepsilon/2}(f_{\varepsilon/2})$, and hence $s(S_{\varepsilon}(f))$ is nowhere dense in $X$. Now (ii) immediately follows because for every $f \in J$, $\osupp(f)$ has empty interior. To show (iii), as \an{$\Rarr$} is clear, it suffices to prove \an{$\Larr$}. Let $f \in C_r^*(G)$ be such that $\osupp(f)$ has empty interior. As 
$
s(\osupp(f)) = \bigcup_{n=1}^\infty s( S_{1/n}(f) )
$, our third claim implies that $s(\osupp(f))$ is the countable union of nowhere dense sets, and thus has empty interior by the Baire category theorem.
\eproof
Note that Lemma~\ref{lem:Sing=Sing} strengthens \cite[Proposition~7.18]{KM} in our situation.
\setlength{\parindent}{0cm} \setlength{\parskip}{0.5cm}

\subsection{Complete characterisations of vanishing of singular ideals}

Let us now establish a characterisation for vanishing of $J \cap \Cc(G)$ as well as vanishing of $J_R$, for arbitrary rings $R$. We need the following notation: Given a ring $R$, we write $\ord(R) \coloneq \menge{\ord(x)}{0 \neq x \in R}$, where for $0 \neq x \in R$, $\ord(x)$ is the smallest positive integer $t$ such that $t \cdot x = 0$ in $R$, and $\ord(x) \coloneq 0$ if no such $t$ exists. Moreover, we set $R_0 \coloneq \Qz$ and $R_t \coloneq \Zz / t \Zz$ where $t$ is an integer with $t > 1$. Furthermore, given a non-negative integer $t \neq 1$, we say that an {\'e}tale groupoid $G$ satisfies condition ($\cS_t$) if there exist a positive integer $n$ and open bisections $U_i \subseteq G$ for $1 \leq i \leq n$ such that 
\setlength{\parindent}{0cm} \setlength{\parskip}{0cm}

\begin{itemize}
    \item $s(U_i) = s(U_j)$ for all $1 \leq i, j \leq n$,
    \item $U_i \setminus (\bigcup_{j \neq i} U_j)$ is not empty and has empty interior for all $1 \leq i \leq n$,
    \item $\lspan_{R_t}(\menge{b_I}{I \subseteq \gekl{1, \dotsc, n}, \, \# \, I > 1, \, (\bigcap_{i \in I} U_i) \setminus (\bigcup_{j \notin I} U_j) \neq \emptyset}) \neq R_t^n$, were $b_I \coloneq \sum_{i \in I} b_i$ and $b_i$ is the $i$-th standard basis vector in $R_t^n$ (all components of $b_i$ are zero except the $i$-th, which is $1$).
\end{itemize}

\btheo
\label{thm:CharSingIdVanish}
Let $G$ be a non-Hausdorff {\'e}tale groupoid. 
\setlength{\parindent}{0cm} \setlength{\parskip}{0cm}

\begin{enumerate}[label=(\roman*)]
    \item We have $J \cap \Cc(G) \neq \gekl{0}$ if and only if $G$ satisfies condition ($\cS_0$).
    \item Let $G$ be ample and $R$ a ring. We have $J_R \neq \gekl{0}$ if and only if there exists $t \in \ord(R)$ such that $G$ satisfies condition ($\cS_t$).
\end{enumerate}
\etheo
\setlength{\parindent}{0cm} \setlength{\parskip}{0cm}

\bproof
Let us prove (i). To see \an{$\Larr$}, let $U_i$ be as in condition ($\cS_0$). Write 
\[
 \Hz \coloneq \lspan_{\Cz}( \Big\{ b_I \colon \ I \subseteq \gekl{1, \dotsc, n}, \, \# \, I > 1, \, (\bigcap_{i \in I} U_i) \setminus (\bigcup_{j \notin I} U_j) \neq \emptyset \Big\} ).
\]
Condition ($\cS_0$) implies $\Hz \neq \Cz^n$ (compare dimensions), so there exists a non-zero linear map $\omega \colon \Cz^n \to \Cz$ such that $\Hz \subseteq \ker(\omega)$. For $1 \leq i \leq n$, write $a_i \coloneq \omega(b_i)$. Since $\omega$ is not zero, there exists $i_0 \in \gekl{1, \dotsc, n}$ with $a_{i_0} \neq 0$. Let $W \coloneq s(U_i)$ (note that $W$ does not depend on $i \in \gekl{1, \dotsc, n}$ by assumption). Choose a function $\psi \in C_c(W)$ such that $\psi(w_0) = 1$ for some $w_0 \in s(U_{i_0} \setminus (\bigcup_{j \neq i_0} U_j))$. For every $1 \leq i \leq n$, let $f_i \coloneq a_i (\psi \circ s \vert_{U_i})$. Then $f_i \in C_c(U_i)$. Set $f \coloneq \sum_{i = 1}^n f_i$. For the element $g_0 \in U_{i_0} \setminus (\bigcup_{j \neq i_0} U_j)$ with $s(g_0) = w_0$, we have $f(g_0) = f_{i_0}(g_0) = a_{i_0} \psi(w_0) = a_{i_0} \neq 0$. This shows that $f \neq 0$. If $g \in \bigcup_{i \neq j \, \in \, \gekl{1, \dotsc, n}} U_i \cap U_j$, then there exists $I \subseteq \gekl{1, \dotsc, n}$ with $\# \, I > 1$ such that $g \in (\bigcap_{i \in I} U_i) \setminus (\bigcup_{j \notin I} U_j)$. In that case, $b_I \in \Hz$ and thus $f(g) = \sum_{i \in I} f_i(g) = (\sum_{i \in I} a_i) \psi(s(g)) = \omega(b_I) \psi(s(g)) = 0$. This shows that $\osupp(f) \subseteq \bigcup_{i=1}^n (U_i \setminus (\bigcup_{j \neq i} U_i \cap U_j))$ and hence $\osupp(f)$ has empty interior because the $U_i$ are open. We conclude that $f \in J$ by Lemma~\ref{lem:Sing=Sing}~(iii).
\setlength{\parindent}{0cm} \setlength{\parskip}{0.5cm}

Let us prove \an{$\Rarr$}. Let $n$ be the smallest positive integer for which there exist open bisections $V_1, \dotsc, V_n$ and $f_i \in C_c(V_i)$ for $1 \leq i \leq n$ such that $f = \sum_{i=1}^n f_i$ lies in $J$ and $f \neq 0$. Note that $n > 1$ because for $n=1$, $\osupp(f)$ is open. Write $S_i \coloneq \osupp(f_i)$ and let $K_i$ be the closure of $S_i$ in $V_i$. Then 
\begin{equation}
\label{e:S-SNotEmpty}
 S_i \setminus (\bigcup_{j \neq i} S_j) \neq \emptyset \text{ for all } 1 \leq i \leq n,   
\end{equation}
because otherwise there exists $i$ such that $\osupp(f) \subseteq \bigcup_{j \neq i} K_j \subseteq \bigcup_{j \neq i} V_j$, and then Corollary~\ref{cor:f=sum}~(i) would lead to a contradiction to minimality of $n$. Moreover, minimality of $n$ implies that 
\begin{equation}
\label{e:f=0OnInt}
f(g) = 0 \text{ for all } g \in S_i \cap S_j \text{ for } i \neq j. 
\end{equation}
If not, find $g' \in S_i \cap S_j$ for $i \neq j$ such that $f(g') \neq 0$ and $\psi \in C_c(s(S_i \cap S_j))$ such that $\psi(s(g')) = 1$. Set $f' \coloneq f * \psi$. Then $f' \in J$, $f'(g') \neq 0$ and $f' = \sum_{i=1}^n f'_i$, where $f'_i \coloneq f_i * \psi$. For every $1 \leq k \leq n$, we have $S'_k \coloneq \osupp(f'_k) = S_k \cap s^{-1}(\osupp(\psi))$. Hence $S'_i = S'_j$, which contradicts \eqref{e:S-SNotEmpty}. Furthermore, proceeding inductively on $1 \leq k \leq n$, let us now construct $\psi_k \in C_c(\bigcap_{i=1}^k s(S_i))$ such that $f * \psi_k \neq 0$. For $k=1$, choose $g \in S_1 \setminus (\bigcup_{j \neq 1} S_j)$ and $\psi_1 \in C_c(s(S_1))$ such that $\psi_1(s(g)) = 1$. Now assume that we have found $\psi_{k-1}$. Then $f' = f * \psi_{k-1}$ is a non-zero function in $J$ with $f' = \sum_{i=1}^n f'_i$, where $f'_i = f_i * \psi_{k-1}$. Set $S'_i \coloneq \osupp(f'_i)$, for $1 \leq i \leq n$. We have $S'_i \neq \emptyset$ by minimality of $n$. Furthermore, \eqref{e:S-SNotEmpty} implies that $S'_k \setminus (\bigcup_{j \neq k} S'_j) \neq \emptyset$, so that we can choose $g_k \in S'_k \setminus (\bigcup_{j \neq k} S'_j)$ and $\varphi_k \in C_c(s(S'_k))$ with $\varphi_k(s(g_k)) = 1$. Then set $\psi_k \coloneq \psi_{k-1} \varphi_k$. This concludes the construction. Thus, by replacing $f$ by $f * \psi_n$, we arrange that the $S_i$ from above satisfy $s(S_i) = s(S_j)$ for all $i, j \in \gekl{1, \dotsc, n}$. Set $W \coloneq s(S_i)$. 

For $w \in W$ and $1 \leq i \leq n$, let $g_{i,w}$ be the unique element of $S_i$ with $s(g_{i,w}) = w$. Consider the subspaces $\Hz_{f,w} \coloneq \lspan_{\Cz}(\menge{b_I}{\sum_{i \in I} f_i(g_{i,w}) = 0})$ and $\Hz'_{f,w} \coloneq \menge{\sum_{i=1}^n a_i b_i \in \Cz^n}{\sum_{i=1}^n a_i f_i(g_{i,w}) = 0}$. By construction, we have $\Hz_{f,w} \subseteq \Hz'_{f,w}$. For every $g \in G$ with $f(g) \neq 0$, we have $\Hz'_{f,s(g)} \neq \Cz^n$ and thus $\Hz_{f,s(g)} \neq \Cz^n$. Since $X \to \Cz, \, w \ma \sum_{i \in I} f_i(g_{i,w}) = \sum_{g \in G_w} (\sum_{i \in I} f_i)(g)$ is continuous by Remark~\ref{rem:LambdaContFors}, the set $W_I \coloneq \menge{w \in X}{\sum_{i \in I} f_i(g_{i,w}) \neq 0}$ is an open subset of $W$. Set $\Iz_w \coloneq \menge{I \subseteq \gekl{1, \dotsc, n}}{b_I \notin \Hz'_{f,w}}$. Then $W_w \coloneq \bigcap_{I \in \Iz_w} W_I$ is an open neighbourhood of $w$ such that $\Hz_{f,w'} \subseteq \Hz_{f,w}$ for all $w' \in W_w$. Now take $g_0 \in G$ with $f(g_0) \neq 0$ and set $w_0 \coloneq s(g_0)$. Choose $\psi \in C_c(W_{w_0})$ such that $\psi(w_0) = 1$. Replace $f$ with $f' \coloneq f * \psi$. We have $f' = \sum_{i=1}^n f'_i$, where $f'_i = f_i * \psi$. Set $U_i \coloneq \osupp(f'_i) = S_i \cap s^{-1}(\menge{w \in W_{w_0}}{\psi(w) \neq 0})$. Since $\sum_{g \in G_w} (\sum_{i \in I} f_i * \psi) (g) = (\sum_{g \in G_w} (\sum_{i \in I} f_i)(g)) \psi(w)$, we deduce that $\Hz_{f',s(g)} = \Hz_{f,s(g)} \subseteq \Hz_{f,w_0}$ for all $g \in \bigcup_{i=1}^n U_i$.

We claim that these $U_i$ have all the desired properties. We have $s(U_i) = \osupp(\psi) = s(U_j)$ for all $i, j \in \gekl{1, \dotsc, n}$. Moreover, $U_i \setminus (\bigcup_{j \neq i} U_j)$ is not empty by \eqref{e:S-SNotEmpty}. By construction, $(\bigcup_i U_i) \setminus (\bigcup_{i \neq j} U_i \cap U_j)$ is contained in $\osupp(f')$. As $f' \in J$, this shows that $U_i \setminus (\bigcup_{i \neq j} U_i \cap U_j)$ has empty interior for each $i$. Moreover, for every $I \subseteq \gekl{1, \dotsc, n}$ with $\# \, I > 1$ and $g \in (\bigcap_{i \in I} U_i) \setminus (\bigcup_{j \notin I} U_j)$, we have $f'(g) = 0$ by \eqref{e:f=0OnInt}. Hence $0 = f'(g) = \sum_{i \in I} f'_i(g) = \sum_{i \in I} f'_i(g_{i,s(g)})$ and thus $b_I \in \Hz_{f',s(g)} \subseteq \Hz_{f,w_0}$. This shows that $\lspan_{\Cz}(\menge{b_I}{I \subseteq \gekl{1, \dotsc, n}, \, \# \, I > 1, \, (\bigcap_{i \in I} U_i) \setminus (\bigcup_{j \notin I} U_j) \neq \emptyset}) \subseteq \Hz_{f,w_0} \neq \Cz^n$, as desired.

Let us now prove (ii). Given a ring $R$, let $A$ be the underlying abelian group of $R$. Let $\Cc(G,A) \coloneq \lspan(\menge{c_U}{c \in A, \, U \subseteq G \text{ compact open bisection}})$ and let $\J_A$ be the subgroup of $\Cc(G,A)$ of those functions $f \in \Cc(G,A)$ with the property that $\osupp(f)$ has empty interior. Clearly, we have $J_R \cong \J_A$ as abelian groups. In other words, whether or not we have $J_R \neq \gekl{0}$ only depends on the underlying abelian group of $R$. As $A = \ilim A_{\rm fg}$, where $A_{\rm fg}$ runs through all finitely generated subgroups of $A$ (ordered by inclusion), we have $\J_A = \ilim \J_{A_{\rm fg}}$. Thus $\J_A \neq \gekl{0}$ if and only if $\J_{A_{\rm fg}} \neq \gekl{0}$ for some finitely generated subgroup $A_{\rm fg}$ of $A$. Now every finitely generated subgroup of $A$ is isomorphic to $\bigoplus_{l=1}^L \Zz / t_l \Zz$ for some non-negative integer $L$ and $t_l \in \ord(R)$. Moreover, $\J_{\bigoplus_{l=1}^L \Zz / t_l \Zz}$ embeds into $ \bigoplus_{l=1}^L \J_{\Zz / t_l \Zz}$. We conclude that $J_R \neq \gekl{0}$ if and only if there exists $t \in \ord(R)$ such that $J_{\Zz / t \Zz} \neq \gekl{0}$.

To prove (ii), it thus suffices to show that for a non-negative integer $t \neq 1$ and the ring $R = \Zz / t \Zz$, we have $J_R \neq \gekl{0}$ if and only if $G$ satisfies condition ($\cS_t$). For $t = 0$, it is easy to see that $J_{\Zz} \neq \gekl{0}$ if and only if $J_{\Qz} \neq \gekl{0}$, and the same argument as for (i) shows that $J_{\Qz} \neq \gekl{0}$ if and only if $G$ satisfies condition ($\cS_0$). Now consider the case $t > 0$. That $J_R \neq \gekl{0}$ implies that $G$ satisfies condition ($\cS_t$) follows by the same argument as for the corresponding implication in (i), but now using Corollary~\ref{cor:f=sum}~(ii). It remains to show that if $G$ satisfies condition ($\cS_t$), then $J_R \neq \gekl{0}$, where $R = \Zz / t \Zz$. Define $\Hz$ as in (i), but with $\lspan_R$ instead of $\lspan_\Cz$. Then we have $\Hz \neq R^n$, and once we construct a non-zero $R$-linear map $\omega \colon R^n \to R$ with $\Hz \subseteq \ker(\omega)$, the same argument as for \an{$\Larr$} in (i) applies. We have $R^n / \Hz \cong \bigoplus_{m=1}^M \Zz / \tau_m \Zz$ for some positive integer $M$ and integers $\tau_m > 1$ dividing $t$. Set $\tau \coloneq \tau_1$ and let $\omega$ be the composition $R^n \onto R^n / \Hz \congto \bigoplus_{m=1}^M \Zz / \tau_m \Zz \onto \Zz / \tau \Zz \into \Zz / t \Zz = R$, where the first map is the canonical projection, the second map is an arbitrary isomorphism of abelian groups, the third map is given by the canonical projection and the fourth map is given by multiplication with $t/\tau$. Then $\omega$ is clearly $R$-linear, and we have $\Hz \subseteq \ker(\omega)$ as well as $\omega \neq 0$ by construction.
\eproof
\setlength{\parindent}{0cm} \setlength{\parskip}{0cm}

\bremark
If $G$ is ample, then the $U_i$ in condition ($\cS_t$) can always be arranged to be compact open bisections. Moreover, the last part of condition ($\cS_0$), that the span of the $b_I$ is a proper subspace of $\Qz^n$, can be decided by computing determinants of (finitely many) matrices whose columns are given by certain $b_I$.
\eremark

\bremark
\label{rem:OnlyChar}
Let $\Kz$ be a field. Then clearly $\ord(\Kz)$ only consists of the characteristic of $\Kz$. Hence, Theorem~\ref{thm:CharSingIdVanish} tells us that whether or not $J_{\Kz} = \gekl{0}$ only depends on the characteristic of $\Kz$. This explains \cite[Theorem~5.9]{StSz21}.
\eremark
\setlength{\parindent}{0cm} \setlength{\parskip}{0.5cm}

The following is an immediate consequence of Theorem~\ref{thm:CharSingIdVanish} and \cite[Theorem~A]{StSz21}.
\bcor
\label{cor:SimpSt}
Let $G$ be a non-Hausdorff ample groupoid and $\Kz$ a field of characteristic $p$. Then $\Kz G$ is simple if and only if $G$ is minimal, effective and does not satisfy condition ($\cS_p$).
\ecor

The following is another immediate consequence of Theorem~\ref{thm:CharSingIdVanish}.
\bcor
\label{cor:JC}
Let $G$ be a non-Hausdorff ample groupoid. Then $J \cap \Cc(G) = \gekl{0}$ if and only if $J_{\Cz} = \gekl{0}$.
\ecor

Let us now characterise vanishing of $C^*$-algebraic singular ideals under a finiteness assumption. Let $\overline{X}$ be the closure of $X$ in $G$. Given $x \in X$, we set $\oX(x) \coloneq \oX \cap G_x^x$. 
\btheo
\label{thm:JCc=0vsJ=0}
Let $G$ be a non-Hausdorff {\'e}tale groupoid. Assume that $\oX(x)$ is finite for all $x \in X$. Then $J \cap \Cc(G) = \gekl{0}$ implies $J = \gekl{0}$.
\etheo
\setlength{\parindent}{0cm} \setlength{\parskip}{0cm}

\begin{proof}
Take $x \in X$, let $\langle \oX(x) \rangle$ be the subgroup of $G_x^x$ generated by $\oX(x)$, and set $m \coloneq \# \, \oX(x)$. Let us first show
\begin{equation}
\label{e:NoXx02}
 \# \, \langle \oX(x) \rangle \leq m!.
\end{equation}
To see this, enumerate $\oX(x) = \gekl{\gamma_1, \dotsc, \gamma_m}$. Write $\langle \gamma_i \rangle \coloneq \menge{\gamma_i^{\epsilon}}{\epsilon \in \Zz}$ for the cyclic group generated by $\gamma_i$. Clearly, we have $\langle \gamma_i \rangle \subseteq \oX(x)$. Now we claim that every $g \in \langle \oX(x) \rangle$ is of the form $g = h_1 \dotso h_{\ell}$ for some $h_j \in \langle \gamma_{i_j} \rangle$ with $h_l \notin \langle \gamma_{i_k} \rangle$ for every $k \neq l$. Since the $h_j$ have to be pairwise distinct elements of $\oX(x)$, this will then imply $\# \, \langle \oX(x) \rangle \leq m!$ (and $\ell \leq m$). So it suffices to prove the claim. Given $g \in \langle \oX(x) \rangle$, let $\ell$ be the smallest positive integer such that $g$ can be written as $g = h_1 \dotso h_{\ell}$ for some $h_j \in \langle \gamma_{i_j} \rangle$. Assume that there are $1 \leq k < l \leq \ell$ with $h_l \in \langle \gamma_{i_k} \rangle$. Set $h'_j \coloneq h_j$ for $1 \leq j < k$, $h'_k \coloneq h_k h_l$, $h'_j \coloneq h_l^{-1} h_j h_l$ for $k < j < l$ and $h'_j \coloneq h_{j+1}$ for $l \leq j \leq \ell - 1$. We have $h'_k \in \langle \gamma_{i_k} \rangle$ as well as $h'_j \in h_l^{-1} \langle \gamma_{i_j} \rangle h_l = \langle h_l^{-1} \gamma_{i_j} h_l \rangle$ for $k < j < l$, and $h_l^{-1} \gamma_{i_j} h_l \in \oX(x)$ because $\oX(x)$ is closed under conjugation by elements of $G_x^x$. Moreover, we have $g = h'_1 \dotsm h'_{\ell - 1}$ by construction. However, this contradicts minimality of $\ell$. This concludes the proof of \eqref{e:NoXx02}.
\setlength{\parindent}{0cm} \setlength{\parskip}{0.5cm}

For $f \in C^*_r(G)$, write
$
S_\varepsilon(f) \coloneq \menge{g \in G}{\vert f(g) \vert > \varepsilon}.
$
Assume that $J \neq \gekl{0}$. Then we can find $f \in J$ (i.e., $s(\osupp(f))$ has empty interior) with $f(g_0) \neq 0$ for some $g_0 \in G$, and by modifying $f$ if necessary, we can always arrange that $f(x_0) \neq 0$ for some $x_0 \in X$. Let us fix $x_0$ throughout the proof. Since $\langle \oX(x_0) \rangle$ is finite by \eqref{e:NoXx02}, we can enumerate $\langle \oX(x_0) \rangle = \gekl{g_1, \dotsc , g_n}$ (where $n \le m!$). For each $1 \leq i \leq n$, let $U_i$ be an open bisection containing $g_i$, such that $s(U_i) = s(U_j)$ for all $1 \leq i, j \leq n$. Write $V \coloneq s(U_i)$ (this does not depend on $i$). Choose $\psi \in C_c(V)$ such that $\psi(x_0)=1$, and denote $\omega_i \coloneq f(g_i)$. For every $1 \leq i \leq n$, let $f^\psi_i \coloneq \omega_i (\psi \circ s \vert_{U_i})$ so that $f^\psi_i \in C_c(U_i)$. Set $f^\psi \coloneq \sum_{i = 1}^n f^\psi_i$. Then it is clear that $f^\psi \in \Cc(G)\setminus \gekl{0}$. Let us show that $f^{\psi}$ lies in $J$ for some $\psi$.

Set $S_i \coloneq \osupp(f_i^\psi)$. Given $I \subseteq \gekl{1, \dotsc, n}$, set $\check{S}_I \coloneq (\bigcap_{i \in I} S_i) \setminus (\bigcup_{j \notin I} S_j)$ as well as $W_I \coloneq \inte \Big( \overline{ s(\check{S}_I) }^X \Big)$. We can modify $\psi$ such that 
\begin{equation}
\label{e:x0WI}
 W_I \neq \emptyset \text{ implies that } x_0 \in \overline{W_I}^X \text{ for all } I \subseteq \gekl{1, \dotsc, n}.
\end{equation}
Indeed, choose a regular open subset $S$ of $X$ (i.e., $S = \inte \big( \overline{S}^X \big)$) such that $x_0 \in S$ and $S \cap \overline{W_I}^X = \emptyset$ for all $I \subseteq \gekl{1, \dotsc, n}$ with $W_I \neq \emptyset$ and $x_0 \notin \overline{W_I}^X$. Choose $\varphi \in C_c(S)$ with $\varphi(x_0) = 1$, and consider $\dot{f} \coloneq f^\psi * \varphi = \sum_{i=1}^n f_i^\psi * \varphi = \sum_{i=1}^n \omega_i \big( (\psi \varphi) \circ s \vert_{U_i}\big) = f^{\psi \varphi}$. Then $\dot{S}_i \coloneq \osupp(f_i^\psi * \varphi)$ is given by $\dot{S}_i = S_i \cap s^{-1}\big(\osupp(\varphi)\big)$, so that $(\dot{S}_I)\check{\vrule height1.3ex width0pt} \coloneq (\bigcap_{i \in I} \dot{S}_i) \setminus (\bigcup_{j \notin I} \dot{S}_j)$ satisfies $(\dot{S}_I)\check{\vrule height1.3ex width0pt} = \check{S}_I \cap s^{-1}\big(\osupp(\varphi)\big)$. Therefore, with $\dot{W}_I \coloneq \inte \Big( \overline{s((\dot{S}_I)\check{\vrule height1.3ex width0pt} \, )}^X \Big)$, we have $\dot{W}_I = \inte \Big( \overline{s(\check{S}_I) \cap \osupp(\varphi)}^X \Big) \subseteq \inte \Big( \overline{s(\check{S}_I)}^X \cap \overline{S}^X \Big) \subseteq W_I \cap \inte \big( \overline{S}^X \big) = W_I \cap S = \emptyset$ whenever $x_0 \notin \overline{W_I}^X$, and $x_0 \notin \overline{W_I}^X$ if and only if $x_0 \notin \overline{\dot{W}_I}^X$.

For each $x \in V$, let $g_x^i$ be the unique element in $U_i$ satisfying $s(g_x^i) = x$. We claim that the map $\chi \colon V \to \Cz$ given by
\begin{equation}
\label{e:CtsAtx_0}
\chi(x) \coloneq \sum_{i=1}^n |f(g_x^i) - f^\psi(g_x^i)|
\end{equation}
is continuous at $x_0$. Let $(x_\alpha)$ be a net in $V$ converging to $x_0$. In order to prove that $\lim_{\alpha} \chi(x_\alpha) = \chi(x_0)$, it suffices to find a subnet $(x_\beta)$ satisfying $\lim_{\beta} \chi(x_\beta) = \chi(x_0)$ (since $(x_\alpha)$ is arbitrary). As the canonical projection $\pi \colon \tiX \onto X$ is proper, there exists a subnet $(x_\beta)$ for which the net $(\iota(x_\beta))$ converges in $\tiG$ (i.e., in the Fell topology). For each $1 \leq i \leq n$ let $g_{\beta}^i$ be the unique element of $U_i$ such that $s(g_{\beta}^i) = x_{\beta}$. Then the net $(\iota(g_\beta^i))$ also converges in $\tiG$, say to the limit $\bmg_i$. Note that $\bmg_i \subseteq \langle \oX(x_0) \rangle$. We have
\[
\lim_{\beta} f(g_\beta^i) = \big(C_r^*(\mfi)(f)\big)(\bmg_i) 
= \sum_{g \in \bmg_i} f(g)
= \sum_{g \in \bmg_i} f^\psi(g)
= \mfi(f^{\psi})(\bmg_i)
= \lim_{\beta} f^\psi(g_\beta^i)
\ \ \ \ \ \ \text{for each } 1 \leq i \leq n. 
\]
Here we used Lemma~\ref{lem:CinC} for the second and Theorem~\ref{thm:CharCinC}~(ii) for the fourth equality. Therefore
\[
\lim_{\beta} \chi(x_\beta) = \lim_{\beta} \sum_{i=1}^n |f(g_\beta^i) - f^\psi(g_\beta^i)| = 0 = \chi(x_0).
\]
This concludes the proof that the map $\chi$ from \eqref{e:CtsAtx_0} is continuous at $x_0$. Set $\cI \coloneq \menge{I \subseteq \gekl{1, \dotsc, n}}{W_I \neq \emptyset}$. We claim that
\begin{equation}
\label{e:sum=0}
 \sum_{i \in I} f(g_i)= 0 \text{ for all } I \in \cI.
\end{equation}
To see this, take an arbitrary $\varepsilon > 0$. By continuity of $\chi$ at $x_0$, there exists an open neighbourhood $W_\varepsilon$ of $x_0$ such that $|f(g) - f^\psi(g)| < \varepsilon/2$ whenever $g \in s^{-1}(W_\varepsilon) \cap (\bigcup_{i=1}^n S_i)$. In particular, $W_\varepsilon \cap \overline{s(S_{\varepsilon}(f^\psi))}^X \subseteq \overline{W_\varepsilon \cap s(S_{\varepsilon}(f^\psi))}^X \subseteq \overline{s(S_{\varepsilon/2}(f))}^X$ has empty interior by Lemma~\ref{lem:Sing=Sing}~(ii). Now set $\check{W}_I \coloneq (W_\varepsilon \cap W_I) \setminus \overline{s(S_{\varepsilon}(f^\psi))}^X$ for each $I \in \cI$. Then $\check{W}_I$ is open and dense in $W_\varepsilon \cap W_I$ because $W_{\varepsilon} \cap \overline{s(S_{\varepsilon}(f^\psi))}^X$ has empty interior. Therefore $x_0$ lies in the closure of $\check{W}_I$ in $X$ by \eqref{e:x0WI}. As $s(\check{S}_I)$ is dense in $W_I$ and because $\check{W}_I$ is open, we conclude that there is a net $(x_{\lambda})$ in $\check{W}_I \cap s(\check{S}_I)$ converging to $x_0$ in $X$. Let $g_{\lambda}$ be the unique element of $\check{S}_I$ such that $s(g_{\lambda}) = x_{\lambda}$. 
By construction, we have that $( \iota(g_{\lambda}))$ converges to $\bmg = \menge{g_i}{i \in I}$ in $\tiG$ (i.e., in the Fell topology). Therefore, 
\[
 \lim_{\lambda} f^\psi (g_{\lambda}) =  \lim_{\lambda} \mfi(f^\psi)(\iota(g_{\lambda})) = \mfi(f^\psi)(\bmg) = \sum_{i \in I} \mfi(f^\psi)(\iota(g_i)) = \sum_{i \in I} f^\psi(g_i) = \sum_{i \in I} f(g_i).
\]
For the third equality, we used that $\mfi(f^\psi)$ is decomposable by Theorem~\ref{thm:CharCinC}~(ii). Now $s(g_{\lambda}) \in \check{W}_I$ implies that $g_{\lambda} \notin S_{\varepsilon}(f^\psi)$, so that $\abs{f^\psi(g_{\lambda})} \leq \varepsilon$. We conclude that
$
\big\vert \sum_{i \in I} f(g_i) \big\vert \leq \varepsilon
$. However, $\varepsilon > 0$ was arbitrary, and therefore \eqref{e:sum=0} holds.

Let us now show that $f^\psi \in J$. By the definition of $f^\psi$, for $I \in \cI$ and $g \in \check{S}_I$ we have $f^\psi(g) = \sum_{i \in I} f_i^\psi(g) = \sum_{i \in I} \omega_i \psi(s(g)) = \sum_{i \in I} f(g_i) \psi(s(g)) = 0$ using \eqref{e:sum=0}. It follows that $\osupp(f^\psi) \subseteq \bigcup_{I \, \in \, \cP(\gekl{1, \dotsc, n}) \setminus \cI} \check{S}_I$ and therefore $s(\osupp(f^\psi)) \subseteq \bigcup_{I \, \in \, \cP(\gekl{1, \dotsc, n}) \setminus \cI} s(\check{S}_I)$ (here $\cP$ stands for power set). Hence $s(\osupp(f^\psi))$ has empty interior because it is contained in a finite union of nowhere dense sets ($s(\check{S}_I)$ is nowhere dense because $W_I = \emptyset$ for $I \in \cP(\gekl{1, \dotsc, n}) \setminus \cI$). This shows $f^\psi \in J$ and thus $J \cap \Cc(G) \neq \gekl{0}$, as desired.
\end{proof}
\setlength{\parindent}{0cm} \setlength{\parskip}{0cm}

Note that the finiteness condition ($\# \, \oX(x) < \infty$ all $x \in X$) in Theorem~\ref{thm:JCc=0vsJ=0} is satisfied by all groupoids of contracting self-similar groups (see Corollary~\ref{cor:SelfSimFinite}) and all {\'e}tale groupoids with finite isotropy groups.
\setlength{\parindent}{0cm} \setlength{\parskip}{0.5cm}

The following is an immediate consequence of Theorem~\ref{thm:JCc=0vsJ=0} and Theorem~\ref{thm:CharSingIdVanish}~(i).
\bcor
\label{cor:JCc=0vsJ=0}
Let $G$ be a non-Hausdorff {\'e}tale groupoid. Assume that $\oX(x)$ is finite for all $x \in X$. Then $J \cap \Cc(G) = \gekl{0}$ if and only if $J = \gekl{0}$. Moreover, $J = \gekl{0}$ if and only if $G$ does not satisfy condition ($\cS_0$).
\ecor

We also obtain the following immediate consequence of Corollary~\ref{cor:JCc=0vsJ=0} and \cite[Theorem~A]{KKLRU} (see also Remark~\ref{rem:Improve}).
\bcor
\label{cor:SimpCSTAR}
Let $G$ be a non-Hausdorff {\'e}tale groupoid. Assume that $\oX(x)$ is finite for all $x \in X$. Further suppose that $G$ has compact unit space, or that $G$ can be covered by countably many open bisections. Then the reduced groupoid $C^*$-algebra $C^*_r(G)$ is simple if and only if $G$ does not satisfy condition ($\cS_0$), is minimal and has no essentially confined amenable sections of isotropy groups in the sense of \cite[Definition~7.1]{KKLRU}.
\ecor
\setlength{\parindent}{0cm} \setlength{\parskip}{0cm}

This provides a partial answer to \cite[Question~2]{CEPSS}.
\setlength{\parindent}{0cm} \setlength{\parskip}{0.5cm}

Moreover, Corollaries~\ref{cor:JCc=0vsJ=0}, \ref{cor:SimpSt} and \cite[Theorem~7.26]{KM} immediately imply the following result. 
\bcor
\label{cor:SimpSteinbCSTAR}
Let $G$ be a non-Hausdorff ample groupoid. Assume that $\oX(x)$ is finite for all $x \in X$. If the complex Steinberg algebra $\Cz G$ is simple, then the reduced groupoid $C^*$-algebra $C^*_r(G)$ is simple.
\ecor
\setlength{\parindent}{0cm} \setlength{\parskip}{0cm}

This provides a partial answer to \cite[Question~3]{CEPSS}.

\setlength{\parindent}{0cm} \setlength{\parskip}{0.5cm}

\bquestions
In general, i.e., without the hypothesis that $\oX(x)$ is finite for all $x \in X$, the following two questions remain open.
\setlength{\parindent}{0cm} \setlength{\parskip}{0cm}
\begin{enumerate}[label=(\Roman*)]
    \item Which groupoid property characterises vanishing of the $C^*$-algebraic singular ideal?
    \item When does $J \cap \Cc(G) = \gekl{0}$ imply $J = \gekl{0}$?
\end{enumerate}
The following is another natural, related question.
\begin{enumerate}
    \item[(III)] When is $J \cap \Cc(G)$ dense in $J$?
\end{enumerate}
\equestions
\setlength{\parindent}{0cm} \setlength{\parskip}{0.5cm}

\subsection{Description of the singular ideal and sufficient vanishing criteria}
\label{ss:DescSingId}

In this subsection, we assume that $G$ can be covered by countably many open bisections. Let us now describe the singular ideal of $C^*_r(G)$ using the Hausdorff cover $\ti{G}$. As above, we write $X \coloneq G^{(0)}$ and $\ti{X} \coloneq \ti{G}^{(0)}$. Recall that we constructed a proper continuous surjection $\pi \colon \ti{X} \onto X$ characterised by $\gekl{\pi(\bmx)} = \bmx \cap X$ as well as a canonical (non-continuous) inclusion $\iota \colon X \into \tiX$. Now we define the set of dangerous points $D$ as $D \coloneq \pi(\tiX \setminus \iota(X))$. Note that $D = s(\overline{X} \setminus X)$, where $\overline{X}$ is the closure of $X$ in $G$. We also set $C \coloneq X \setminus D$. Note that $C$ is precisely the set of points where $\iota \colon X \into \ti{X}$ is continuous, and $\iota(C)$ is precisely the set of points in $\tiX$ where $\pi$ is one-to-one. Moreover, \cite[Lemma~7.15]{KM} implies that $C$ is dense in $X$. As explained in \cite[\S~7]{KM}, the singular ideal $J$ of $C^*_r(G)$ can be equivalently described as the set of those $f \in C^*_r(G)$ which satisfy $\osupp(f) \subseteq r^{-1}(D) = s^{-1}(D)$. The essential $C^*$-algebra of $G$ is then defined as $C^*_\ess(G) = C^*_r(G) / J$.

\bdefin
Define $\ti{X}_\sing \coloneq \inte(\pi^{-1}(D)) \subseteq \ti{X}$ and $\ti{X}_\ess \coloneq \ti{X} \setminus \ti{X}_\sing$.
\edefin
\setlength{\parindent}{0cm} \setlength{\parskip}{0cm}

By construction, $\ti{X}_\ess$ coincides with the closure of $\iota(C)$ in $\ti{X}$, i.e., $\ti{X}_\ess = \overline{\iota(C)} \subseteq \ti{X}$. It is straightforward to see that $\pi^{-1}(D)$ is $\ti{G}$-invariant. Hence $\ti{X}_\sing$ is an open invariant subset of $\ti{X}$, and thus $\ti{X}_\ess$ is a closed invariant subset of $\ti{X}$. Also note that the canonical projection $\pi \colon \: \ti{X} \onto X$ restricts to a surjective map $\ti{X}_\ess \onto X$. This follows from the observation that $C \subseteq \pi(\tiX_\ess)$ and that $C$ is dense in $X$.
\setlength{\parindent}{0cm} \setlength{\parskip}{0.5cm}

\bdefin
Define $\ti{G}_\sing \coloneq \ti{G}_{\ti{X}_\sing}^{\ti{X}_\sing} = \menge{\bmg \in \tiG}{\ti{r}(\bmg) \in \tiX_\sing, \, \ti{s}(\bmg) \in \tiX_\sing}$ and $\ti{G}_\ess \coloneq \ti{G}_{\ti{X}_\ess}^{\ti{X}_\ess} =  \ti{G} \setminus \ti{G}_\sing$.
\edefin
\setlength{\parindent}{0cm} \setlength{\parskip}{0cm}

Note that $\ti{G}_\sing$ is an open subgroupoid of $\ti{G}$ and $\ti{G}_\ess$ is a closed subgroupoid of $\ti{G}$. Moreover, $\tiG_\ess$ is the closure of $\iota(G_C^C)$ in $\tiG$, and $G_C^C$ are precisely the set of points where the canonical inclusion $\iota \colon G \into \tiG$ is continuous.
\setlength{\parindent}{0cm} \setlength{\parskip}{0.5cm}

\bdefin
Define $\ti{J} \coloneq \menge{\tif \in C^*_r(\ti{G})}{\osupp(\tif) \subseteq \ti{G}_\sing}$.
\edefin
\setlength{\parindent}{0cm} \setlength{\parskip}{0cm}

In other words, $\ti{J}$ is the ideal of $ C^*_r(\ti{G})$ consisting of those $\tif \in C^*_r(\ti{G})$ with the property that $\tif \vert_{\ti{G}_\ess} = 0$. Hence, by construction, we have a canonical isomorphism $C^*_r(\ti{G}) / \ti{J} \cong C^*_r(\ti{G}_\ess)$ given by restriction to $\ti{G}_\ess$ (and we do not need exactness here).
\setlength{\parindent}{0cm} \setlength{\parskip}{0.5cm}

The following is our description of the singular ideal $J$ in terms of $\ti{G}$.
\bprop
\label{prop:Jsing=CcapJ}
If $G$ can be covered by countably many open bisections, then $C^*_r(\mfi)(J) = C^*_r(\mfi)(C^*_r(G)) \cap \ti{J}$.
\eprop
\setlength{\parindent}{0cm} \setlength{\parskip}{0cm}

\bproof
As mentioned above, the ideal $J$ consists of precisely those $f \in C^*_r(G)$ which are non-zero only in points of $s^{-1}(D)$. In other words, a function $f \in C^*_r(G)$ lies in $J$ if and only if $\iota^{-1}(\osupp(\tif)) \subseteq s^{-1}(D)$, where $\tif \coloneq C^*_r(\mfi)(f)$. We claim that this is equivalent to the statement that $\osupp(\tif) \subseteq \tis^{-1}(\pi^{-1}(D))$. Indeed, if the latter holds, and we are given $g \in G$ with $\tif(\iota(g)) \neq 0$, then $\iota(g) \in \tis^{-1}(\pi^{-1}(D))$, and thus $s(g) = \pi(\iota(s(g))) = \pi(\tis(\iota(g))) \in D$. Conversely, suppose that $\iota^{-1}(\osupp(\tif)) \subseteq s^{-1}(D)$, and take $\bmg \in \osupp(\tif)$. If $\tis(\bmg) \notin \pi^{-1}(D)$, then $\tis(\bmg) \in \pi^{-1}(C) = \iota(C)$ (because $\pi$ is one-to-one on $\iota(C)$). But then $\bmg$ lies in $\iota(s^{-1}(C))$, which contradicts $\iota^{-1}(\osupp(\tif)) \subseteq s^{-1}(D)$.
\setlength{\parindent}{0.5cm} \setlength{\parskip}{0cm}

Now $\tif$ is a actually a continuous function on $\ti{G}$ because $\ti{G}$ is Hausdorff. Thus our statement that $\osupp(\tif) \subseteq \tis^{-1}(\pi^{-1}(D))$ is equivalent to $\osupp(\tif) \subseteq \inte(\tis^{-1}(\pi^{-1}(D))) = \tis^{-1}(\inte(\pi^{-1}(D))) = \ti{G}_\sing$, as desired.
\eproof
\setlength{\parindent}{0cm} \setlength{\parskip}{0.5cm}

\bremark
\label{rem:JR=cap}
Let us now explain the analogue of Proposition~\ref{prop:Jsing=CcapJ} in the algebraic setting. Let $G$ be an ample groupoid which can be covered by countably many open bisections. In that case, the singular ideal $J_R$ of $RG$ (see for instance \cite{CEPSS,StSz21,StSz23}) can be equivalently described as the set of $f \in RG$ with the property that $\osupp(f) \subseteq r^{-1}(D) = s^{-1}(D)$. Now define $\ti{J}_R \coloneq \menge{\tif \in R \tiG}{\osupp(\tif) \subseteq \tiG_\sing}$. Similar arguments as for Proposition~\ref{prop:Jsing=CcapJ} show that $\mfi_R(J_R) = \mfi_R(RG) \cap \ti{J}_R$.
\eremark

Proposition~\ref{prop:Jsing=CcapJ} immediately gives the following sufficient condition for vanishing of $J$ in terms of $\ti{G}$.

\bcor
\label{cor:tiXsing=0}
Assume that $G$ can be covered by countably many open bisections, and that $\ti{X}_\sing = \emptyset$, i.e., $\pi^{-1}(D)$ has empty interior in $\tiX$. Then the singular ideal $J$ vanishes. If $G$ is ample, then $J_R$ vanishes for any ring $R$.
\ecor

In \cite{CEPSS}, another sufficient condition for vanishing of the singular ideal was given in the case of non-Hausdorff ample groupoids. Namely, it is shown in \cite{CEPSS} that for a non-Hausdorff ample groupoid $G$, the singular ideal vanishes if all compact open subsets of $G$ are regular open. Recall that an open subset $U$ is called regular open if $U = {\rm int} \big( \overline{U} \big)$. It turns out that our sufficient conditions are equivalent for non-Hausdorff ample groupoids.
\blemma
\label{lem:CO=RO<->tiXsing=0}
Let $G$ be ample. All compact open subsets of $G$ are regular open if and only if $\inte(\pi^{-1}(D)) = \emptyset$.
\elemma
\setlength{\parindent}{0cm} \setlength{\parskip}{0cm}

\bproof
Let us prove \an{$\Rarr$}. Suppose that all compact open subsets of $G$ are regular open, but that $\inte(\pi^{-1}(D)) \neq \emptyset$. As $\iota(X)$ is dense in $\tiX$, we can find a compact open bisection $U$ and some compact open set $K$ such that $\bmU \coloneq \menge{\bmg \in \ti{G}}{\bmg \cap U \neq \emptyset, \, \bmg \cap K = \emptyset}$ is contained in $\pi^{-1}(D)$. We claim that $U \cap \big( \overline{K} \big)^c \neq \emptyset$. Indeed, $\bmU$ is non-empty, so contains $\iota(x)$ for some $x \in X$. If $x \notin \overline{K}$, then we are done. Suppose that $x \in \overline{K}$. Since $K$ is regular open by assumption and $x \notin K$, we conclude that $x \notin {\rm int} \big(\overline{K}\big)$. In other words, $x \in \big( {\rm int} \big(\overline{K}\big) \big)^c = \overline{\big( \overline{K} \big) ^c}$. Hence there exists a net $(x_i)$ in $\big( \overline{K} \big) ^c$ with $\lim x_i = x$. As $x \in U$ it follows that $x_i \in U$ eventually. We conclude that $\emptyset \neq U \cap \big( \overline{K} \big) ^c \subseteq \iota^{-1}(\bmU) \subseteq \iota^{-1}(\pi^{-1}(D))$. As $C$ is dense in $X$, we must have $C \cap (U \cap \big( \overline{K} \big) ^c) \neq \emptyset$. But $\iota(C) \cap \pi^{-1}(D) = \emptyset$, a contradiction.
\setlength{\parindent}{0.5cm} \setlength{\parskip}{0cm}

Now let us show \an{$\Larr$}. Assume that $K$ is a compact open subset of $G$ such that $K$ is not regular open, i.e., $U \coloneq {\rm int} \big( \overline{K} \big)$ properly contains $K$. It follows that $U \setminus K$ is non-empty, but has empty interior. Consider $\bmU \coloneq \menge{\bmg \in \ti{G}}{\bmg \cap U \neq \emptyset, \, \bmg \cap K = \emptyset}$. Then $\bmU$ is non-empty, and since $U \setminus K$ has empty interior, $\iota^{-1}(\bmU) \subseteq U \setminus K$ must have empty interior. Hence $\bmU \cap \iota(G_C^C) = \emptyset$ and thus $\tis(\bmU) \subseteq \pi^{-1}(D)$. Thus $\inte(\pi^{-1}(D)) \neq \emptyset$, as desired.
\eproof
\setlength{\parindent}{0cm} \setlength{\parskip}{0.5cm}

\bremark
Lemma~\ref{lem:CO=RO<->tiXsing=0}, together with the example of the groupoid attached to the canonical self-similar action of the Grigorchuk group in \cite[\S~5.6]{CEPSS}, immediately implies that our sufficient condition in Corollary~\ref{cor:tiXsing=0} is not necessary for vanishing of the singular ideal. In Example~\ref{ex:Grig}, we concretely work out $\inte(\pi^{-1}(D))$ in this case.
\eremark

Let us also relate our sufficient condition to a condition in \cite{NS}. In \cite[Proposition~1.12]{NS}, Neshveyev and Schwartz introduce the set $D_0$ of extremely dangerous points and show that if $D_0 = \emptyset$, then the singular ideal vanishes. Recall that $D_0$ is the set of points $x \in X$ with the following property: There exist elements $g_1, \dotsc, g_n \in G_x^x \setminus \gekl{x}$, open bisections $U_1, \dotsc, U_n$ and an open neighbourhood $U \subseteq X$ of $x$ such that $g_k \in U_k$ for all $1 \leq k \leq n$ and $U \setminus (U_1 \cup \dotso \cup U_n)$ has empty interior. Note that $D_0 \subseteq D$. The condition $D_0 = \emptyset$ turns out to be equivalent to our sufficient condition that $\inte(\pi^{-1}(D)) = \emptyset$. This means that Corollary~\ref{cor:tiXsing=0} explains \cite[Proposition~1.12]{NS}.
\blemma
Assume that $G$ can be covered by countably many open bisections. We have $\iota(X) \cap (\inte(\pi^{-1}(D))) = \iota(D_0)$. Moreover, $\inte(\pi^{-1}(D)) = \emptyset$ if and only if $D_0 = \emptyset$.
\elemma
\setlength{\parindent}{0cm} \setlength{\parskip}{0cm}

\bproof
To prove $\iota(X) \cap (\inte(\pi^{-1}(D))) \subseteq \iota(D_0)$, take $x \in X$ with $\iota(x) \in \inte(\pi^{-1}(D))$. Then there are open sets $V_i \subseteq G$ and a compact set $K \subseteq G$ such that there exists $i$ with $V_i \subseteq X$, and
\[
\iota(x) \in \cU(V_i,K) \coloneq \menge{\bmx \in \tiX}{\bmx \cap V_i \neq \emptyset \ \forall \ i, \, \bmx \cap K = \emptyset} \subseteq \pi^{-1}(D).
\]
Then $x \in (\bigcap_i V_i) \setminus K$. Let $V \coloneq \bigcap_i V_i$. Then $V \subseteq X$, $x \in V$ and $V \setminus K \subseteq D$ has empty interior. Moreover, for every $y \in K$, we can find an open bisection $V_y \subseteq G$ such that $y \in V_y$ and, if $r(y) \neq x$ or $s(y) \neq x$, then there is an open set $U_y \subseteq X$ with $x \in U_y$ and $V_y \cap U_y = \emptyset$, and if $r(y) = x = s(y)$, then $x \notin V_y$. As $K$ is compact, we can cover $K$ by finitely many such $V_y$, say $K \subseteq \bigcup_{j \in F} V_{y_j}$ (where $F$ is a finite index set). Set $\check{F} \coloneq \menge{j \in F}{r(y_j) \neq x \text{ or } s(y_j) \neq x}$ and $F_x \coloneq F \setminus \check{F}$. Set $U \coloneq V \cap \bigcap_{j \in \check{F}} U_{y_j}$. Then $U \subseteq X$ is an open neighbourhood of $x$. Moreover, $U \cap K \subseteq \bigcup_{j \in F_x} V_{y_j}$ and thus $U \setminus (\bigcup_{j \in F_x} V_{y_j}) \subseteq U \setminus K$ has empty interior. Comparing this with the definition of $D_0$ as recalled above, we see that $x \in D_0$, as desired.
\setlength{\parindent}{0.5cm} \setlength{\parskip}{0cm}

Let us prove $\inte(\pi^{-1}(D)) \supseteq \iota(D_0)$. Take $x \in D_0$. First of all, note that we can always arrange that the bisections $U_k$ in the definition of $D_0$ to satisfy $r(U_k) = U$ for all $1 \leq k \leq n$ (just replace $U$ by $\bigcap_k r(U_k)$ and then $U_k$ by $r^{-1}(U) \cap U_k$). Now choose an open set $V$ and a compact set $K$ such that $x \in V \subseteq K \subseteq U$, and set
$K_k \coloneq r^{-1}(K) \cap U_k$. Clearly, $K_k$ are compact sets. By construction, $\iota(x) \in \cU(V,\bigcup_k K_k) \coloneq \menge{\bmg \in \tiG}{\bmg \cap V \neq \emptyset, \, \bmg \cap (\bigcup_k K_k) = \emptyset}$. We claim that $\cU(V,\bigcup_k K_k) \subseteq \pi^{-1}(D)$. Indeed, if not, then there would exist $c \in C$ with $\iota(c) \in \cU(V,\bigcup_k K_k)$. As $\iota \colon X \into \tiX$ is continuous at $c$, it follows that there exists an open set $W \subseteq X$ with $c \in W$, $W \subseteq V$ and $W \cap (\bigcup_k K_k) = \emptyset$. But this is impossible because $V \setminus (\bigcup_k K_k) = V \setminus (\bigcup_k U_k)$ has empty interior. Therefore, we conclude that $\iota(x) \in \inte(\pi^{-1}(D))$, as desired.
\setlength{\parindent}{0cm} \setlength{\parskip}{0.5cm}

The second claim follows from the first because $\iota(X)$ is dense in $\tiX$.
\eproof
\setlength{\parindent}{0cm} \setlength{\parskip}{0.5cm}

Let us formulate another sufficient condition for vanishing of singular ideals, which applies to cases where $\inte(\pi^{-1}(D)) \neq \emptyset$. The following is now an immediate consequence of Proposition~\ref{prop:Jsing=CcapJ} and Lemma~\ref{lem:CinC}.
\bcor
\label{cor:SuffCondSingVan}
Let $G$ be a non-Hausdorff {\'e}tale groupoid and $\tiG$ its Hausdorff cover. Assume that $G$ can be covered by countably many open bisections. 
Suppose that $G$ has the property that $\# \, \bmg < \infty$ for all $\bmg \in \tiG$, and that for every decomposable function $\tif \colon \tiG \to \Cz$, $\tif(\bmg) = 0$ for all $\bmg \in \tiG_{\ess}$ implies that $\tif(\bmg) = 0$ for all $\bmg \in \tiG$. Then $J = \gekl{0}$. 
\ecor
\setlength{\parindent}{0cm} \setlength{\parskip}{0cm}

Corollary~\ref{cor:SuffCondSingVan} can be used to explain why for the groupoid attached to the Grigorchuk group (see Example~\ref{ex:Grig}), we have $J = \gekl{0}$ even though $\inte(\pi^{-1}(D)) \neq \emptyset$.
\setlength{\parindent}{0cm} \setlength{\parskip}{0.5cm}

\section{The ideal intersection property}
\label{s:IIP}

In the following, let $G$ be a non-Hausdorff {\'e}tale groupoid and $\tiG$ its Hausdorff cover. As above, we set $X \coloneq \Gn$ and $\tiX \coloneq \tiG^{(0)}$. We assume throughout this section that $G$ can be covered by countably many open bisections. Let us start with a few simple observations.
\blemma
\label{lem:Min}
$G$ is minimal if and only if $\ti{G}_\ess$ is minimal. Moreover, $G$ is topologically transitive if and only if $\tiG_\ess$ is topologically transitive.
\elemma
\setlength{\parindent}{0cm} \setlength{\parskip}{0cm}

Here $G$ is called topologically transitive if there exists $x \in X$ such that the orbit $G.x$ is dense in $X$.

\bproof
Let $x \in X$ be arbitrary, and choose $\bmx \in \tiX_\ess$ with $\pi(\bmx) = x$. Suppose that $\tiG_\ess.\bmx$ is dense in $\tiX$. Then we claim that $G.x$ is dense in $X$. Indeed, we have $X = \pi(\tiX_\ess) = \pi \Big( \overline{\tiG_\ess.\bmx} \Big) \subseteq \overline{G.x}$. Conversely, take $\bmx \in \tiX_\ess$ and suppose that $x \coloneq \pi(\bmx)$ has dense orbit in $G$, i.e., $X = \overline{G.x}$. Then for every $y \in C$, there exists a net $g_i.x$ converging to $y$. As $\pi$ is proper, we may after passing to a subnet assume that $g_i.\bmx$ converges to some $\bmy \in \tiX_\ess$. It follows that $\pi(\bmy) = y \in C$ and thus $\bmy = \iota(y)$ (because $\pi$ is one-to-one on $\iota(C)$). This shows that $\iota(C) \subseteq \overline{\tiG_\ess.\bmx}$, and thus, $\overline{\tiG_\ess.\bmx} = \tiX_\ess$.
\eproof
\setlength{\parindent}{0cm} \setlength{\parskip}{0.5cm}

\blemma
\label{lem:TPTFEFF}
The following are equivalent:
\setlength{\parindent}{0cm} \setlength{\parskip}{0cm}

\begin{enumerate}[label=(\roman*)]
    \item $G$ is topologically principal,
    \item $G$ is topologically free,
    \item $\ti{G}_\ess$ is topologically principal,
    \item $\ti{G}_\ess$ is topologically free,
    \item $\tiG_{\ess}$ is effective.
\end{enumerate}
\elemma
\setlength{\parindent}{0cm} \setlength{\parskip}{0cm}

Recall that $G$ is topologically principal if $\menge{x \in X}{G_x^x = \gekl{x}}$ is dense in $X$, $G$ is topologically free if for every open bisection $U \subseteq G \setminus X$, $\inte(s(U \cap \Isot(G))) = \emptyset$, and $G$ is effective if $\inte(\Isot(G)) = X$.

\bproof
Since $G$ can be covered by countably many open bisections, the same is true for $\tiG_{\ess}$, and thus \an{(i) $\LRarr$ (ii)} and \an{(iii) $\LRarr$ (iv)} follow from \cite[Proposition~2.24]{KM}. \an{(iv) $\LRarr$ (v)} follows from \cite[Lemma~2.23]{KM} as $\tiG$ is Hausdorff. Let us show \an{(i) $\Rarr$ (iv)}. If $G$ is topologically principal, then for every $y \in C$ there exists a net $(x_i)$ converging to $y$ in $X$ such that $G_{x_i}^{x_i} = \gekl{x_i}$. Take $\bmx_i \in \tiX_\ess$ with $\pi(\bmx_i) = x_i$. Then, since $\pi$ is proper, we may assume after passing to a subnet that $(\bmx_i)$ converges in $\tiX_\ess$, and the limit has to be $\iota(y)$. Moreover, $G_{x_i}^{x_i} = \gekl{x_i}$ implies that $\ti{G}_{\bmx_i}^{\bmx_i} = \gekl{\bmx_i}$. As $\iota(C)$ is dense in $\tiX_\ess$, we conclude that $\tiG_\ess$ is topologically free. Finally, let us prove \an{(v) $\Rarr$ (ii)}. If $G$ is not topologically free, then there exists an open bisection $U \subseteq G$ with $U \cap X = \emptyset$ and there exists a non-empty open subset $V \subseteq X$ such that $V \subseteq s(U \cap \Isot(G))$. Define $\bmU \coloneq \menge{\bmg \in \tiG_\ess}{\bmg \in \ti{U}, \, \pi(\ti{s}(\bmg)) \in V}$. We claim that $\bmU \not\subseteq \tiX_\ess$ and $\bmU \subseteq \Isot(\tiG_\ess)$. For the first claim, take $g \in U \cap \Isot(G)$ with $s(g) \in V$. Then $g \notin X$ as $U \cap X = \emptyset$. It follows that $\iota(g) \in \bmU$, and $\iota(g) \notin \tiX$ as $g \notin X$. To show the second claim, take $\bmg \in \bmU$. As $\iota(U)$ is dense in $\tiU$, we can find a net $(g_i)$ in $U$ such that $(\iota(g_i))$ converges to $\bmg$, and without loss of generality we may assume that $\pi(\ti{s}(\iota(g_i))) \in V$ for all $i$. As $\pi(\ti{s}(\iota(g_i))) = \pi(\iota(s(g_i))) = s(g_i)$, we conclude that $s(g_i) \in V$. Since $U$ is a bisection, it follows that $g_i \in \Isot(G)$ and hence $\iota(g_i) \in \Isot(\tiG_\ess)$. We conclude that $\bmg \in \Isot(\tiG_\ess)$, as desired. This shows that $\tiG_\ess$ is not effective.
\eproof
\setlength{\parindent}{0cm} \setlength{\parskip}{0.5cm}

Our goal now is to prove that $C_0(X) \subseteq C^*_\ess(G)$ has the ideal intersection property if and only if $C_0(\ti{X}_\ess) \subseteq C^*_r(\ti{G}_\ess)$ has the ideal intersection property. Recall that $C_0(X) \subseteq C^*_\ess(G)$ has the ideal intersection property if for every ideal $I$ of $C^*_\ess(G)$, $I \cap C_0(X) = \gekl{0}$ implies that $I = \gekl{0}$ (and similarly for $C_0(\ti{X}_\ess) \subseteq C^*_r(\ti{G}_\ess)$). We start with some preparations. In the following, we will view $C_0(X)$ as a subalgebra of $C_0(\tiX_\ess)$ and $C^*_\ess(G)$ as a subalgebra of $C^*_r(\tiG_\ess)$ using the canonical embedding $C^*_r(\mfi)$. Moreover, recall that $C \subseteq X$ are those points such that the canonical projection $\pi \colon \tiX \onto X$ is one-to-one on $\iota(C)$.

\blemma
\label{lem:UniExt:X-tiXess}
For all $x \in C$, $\ev_x \colon C_0(X) \to \Cz, \, f \ma f(x)$ has a unique extension (as a state) from $C_0(X)$ to $C_0(\ti{X}_\ess)$, and that extension is given by $\ev_{\iota(x)} \colon C_0(\tiX_{\ess}) \to \Cz, \, \tif \ma \tif(\iota(x))$.
\elemma
\setlength{\parindent}{0cm} \setlength{\parskip}{0cm}

\bproof
Suppose that $\ti{\varphi}$ is a state on $C_0(\ti{X}_\ess)$ extending $\ev_x$. Then $\ti{\varphi}$ must vanish on the maximal ideal $C_0(\ti{X}_\ess \setminus \gekl{\iota(x)})$. It follows that $\ti{\varphi}$ is given by evaluation at $x \in \ti{X}_\ess$.
\eproof
\setlength{\parindent}{0cm} \setlength{\parskip}{0.5cm}

\blemma
\label{lem:tiG=GtiX=tiXG}
We have $C^*_r(\ti{G}) = \clspan(C^*_r(G) C_0(\ti{X})) = \clspan(C_0(\ti{X}) C^*_r(G))$. Moreover, we also have $C^*_r(\ti{G}_\ess) = \clspan(C^*_\ess(G) C_0(\ti{X}_\ess)) = \clspan(C_0(\ti{X}_\ess) C^*_\ess(G))$.
\elemma
\setlength{\parindent}{0cm} \setlength{\parskip}{0cm}

\bproof
Certainly, the first claim implies the second. So it suffices to prove the first claim. Let $U \subseteq G$ be an open bisection. Consider the open bisection $\ti{U} = \menge{\bmg \in \ti{G}}{\bmg \cap U \neq \emptyset}$. Since $\ti{G} = \bigcup_{U \in \B} \ti{U}$, it suffices to prove that $C_c(\ti{U}) \subseteq C_0(\ti{X}) C_c(U)$ and $C_c(\ti{U}) \subseteq C_c(U) C_0(\ti{X})$. 
\setlength{\parindent}{0.5cm} \setlength{\parskip}{0cm}

Let $\pi_U \colon \ti{U} \onto U$ be the proper continuous surjection from \S~\ref{s:HdCover}. Take $f \in C_c(\ti{U})$ with $\overline{\osupp(f)} = \ti{K}$. Set $K \coloneq \pi_U(\ti{K})$. Then $K$ is a compact subset of $U$. Choose $e \in C_c(U)$ such that $e \equiv 1$ on $K$. Define $f_s$ as the composition $\ti{s}(\ti{U}) \congto \ti{U} \to \Cz$, where the first map is the inverse of $\ti{s} \vert_{\tiU}$, and the second map is $f$. View $f_s$ as an element of $C_0(\ti{X})$. Then $f = e * f_s$. Similarly, define $f_r$ as the composition $\ti{r}(\ti{U}) \cong \ti{U} \to \Cz$, where the first map is the inverse of $\ti{r} \vert_{\tiU}$, and the second map is $f$. View $f_r$ as an element of $C_0(\ti{X})$. Then $f = f_r * e$.
\eproof
\setlength{\parindent}{0cm} \setlength{\parskip}{0.5cm}

\btheo
\label{thm:IIPNonHd}
Let $G$ be a non-Hausdorff {\'e}tale groupoid which can be covered by countably many open bisections. $C_0(X) \subseteq C^*_\ess(G)$ has the ideal intersection property if and only if $C_0(\ti{X}_\ess) \subseteq C^*_r(\ti{G}_\ess)$ has the ideal intersection property.
\etheo
\setlength{\parindent}{0cm} \setlength{\parskip}{0cm}

\bproof
Write $A \coloneq C^*_\ess(G)$ and $\tiA \coloneq C^*_r(\ti{G}_\ess)$. Let us first assume that $C_0(X) \subseteq A$ has the ideal intersection property. Take an ideal $\tiI$ of $\tiA$, and assume that $\tiI \cap C_0(\ti{X}_\ess) = \gekl{0}$. Then $\tiI \cap C_0(X) = \gekl{0}$. Thus the ideal intersection property for $C_0(X) \subseteq A$ implies $\tiI \cap A = \gekl{0}$. This means that the embedding $A \into \tiA$ induces an embedding $A \into \tiA/\tiI$. Let $x \in C$ be arbitrary and $\varphi_x = \spkl{\pi_x(\cdot) \delta_x,\delta_x}$ be the state on $A$ induced from $\ev_x$ on $C_0(X)$, where $\pi_x$ is the left regular representation of $A$ attached to $x$. Extend $\varphi_x$ to a state $\ti{\varphi}_x$ on $\tiA/\tiI$. Then $\ti{\varphi}_x \vert_{C_0(\ti{X}_\ess)}$ extends $\varphi_x \vert_{C_0(X)} = \ev_x$. Lemma~\ref{lem:UniExt:X-tiXess} implies that $\ti{\varphi}_x \vert_{C_0(\ti{X}_\ess)}$ is given by evaluation at $\iota(x)$. In particular, $C_0(\ti{X}_\ess)$ lies in the multiplicative domain of $\ti{\varphi}_x$. We claim that this implies that $\ti{\varphi}_x = \spkl{\ti{\pi}_{\iota(x)}(\cdot) \ti{\delta}_{\iota(x)},\ti{\delta}_{\iota(x)}}$, where $\ti{\pi}_{\iota(x)}$ is the left regular representation of $\tiA$ attached to $\iota(x)$. By Lemma~\ref{lem:tiG=GtiX=tiXG}, it suffices to check this equality for products of the form $a * \tif$ where $a \in A$ and $\tif \in C_0(\ti{X}_\ess)$. And indeed, we have
\[
 \ti{\varphi}_x (a * \tif) = \ti{\varphi}_x(a) \ti{\varphi}_x(\tif) = \varphi_x(a) \tif(\iota(x)) = \spkl{\ti{\pi}_{\iota(x)}(a * \tif)\ti{\delta}_{\iota(x)},\ti{\delta}_{\iota(x)}}.
\]
This shows that $\tiI \subseteq \ker(\ti{\pi}_{\iota(x)})$ for all $x \in C$. But since $\iota(C)$ is dense in $\ti{X}_\ess$, we have that $\bigcap_{x \in C} \ker(\ti{\pi}_{\iota(x)}) = \gekl{0}$. Hence $\tiI = \gekl{0}$, as desired.
\setlength{\parindent}{0cm} \setlength{\parskip}{0.5cm}

Now suppose that $C_0(\ti{X}_\ess) \subseteq \tiA$ has the ideal intersection property. Take an ideal $I$ of $A$ with $I \cap C_0(X) = \gekl{0}$. Let $x \in C$ be arbitrary. Extend $\ev_x$ to a state $\varphi_x$ from $C_0(X)$ to $A / I$. Let $\dot{\varphi}_x$ be the composition $A \onto A / I \to \Cz$, where the first map is the canonical projection and the second map is given by $\varphi_x$. Now extend $\dot{\varphi}_x$ further to a state $\ti{\varphi}_x$ on $\tiA$. Observe that $\ti{\varphi}_x \vert_{C_0(\ti{X}_\ess)}$ extends $\dot{\varphi}_x \vert_{C_0(X)} = \ev_x$. Hence Lemma~\ref{lem:UniExt:X-tiXess} implies that $\ti{\varphi}_x \vert_{C_0(\ti{X}_\ess)}$ is given by evaluation at $\iota(x)$. In particular, $C_0(\ti{X}_\ess)$ lies in the multiplicative domain of $\ti{\varphi}_x$. Define $\tiI \coloneq \menge{a \in \tiA}{\ti{\varphi}_x(\tiA a \tiA) = 0 \ \forall \ x \in C}$. Clearly, $\tiI$ is an ideal of $\tiA$. We claim that the following hold: (i) $I \subseteq \tiI$ and (ii) $\tiI \cap C_0(\ti{X}_\ess) = \gekl{0}$.

Let us prove (i). Take an element $a \in I$. In order to show that $a \in \tiI$, it suffices to show that $\ti{\varphi}_x(\tif_1 * a' * a * a'' * \tif_2) = 0$ for all $a', a'' \in A$ and $\tif_1, \tif_2 \in C_0(\ti{X}_\ess)$ because of Lemma~\ref{lem:tiG=GtiX=tiXG}. And indeed, we have
\[
 \ti{\varphi}_x(\tif_1 * a' * a * a'' * \tif_2) = \ti{\varphi}_x(\tif_1) \ti{\varphi}_x(a' * a * a'') \ti{\varphi}_x(\tif_2) = \ti{\varphi}_x(\tif_1) \dot{\varphi}_x(a' * a * a'') \ti{\varphi}_x(\tif_2) = 0.
\]
Here we used that $C_0(\ti{X}_\ess)$ lies in the multiplicative domain of $\ti{\varphi}_x$ for the first equality. The last equality follows because $a' a a'' \in I$.
\setlength{\parindent}{0.5cm} \setlength{\parskip}{0cm}

Now we prove (ii). Suppose that $\tif \in C_0(\ti{X}_\ess)$ lies in $\tiI$. Then $\ti{\varphi}_x(\tif) = 0$ for all $x \in C$, i.e., $\tif(\iota(x)) = 0$ for all $x \in C$. But $C$ is dense in $\ti{X}_\ess$, and $\tif$ is continuous. Hence $\tif = 0$, as desired.

Now (ii) implies that $\tiI = \gekl{0}$ since $C_0(\ti{X}_\ess) \subseteq \tiA$ has the ideal intersection property. And (i) implies that $I = \gekl{0}$, as desired.
\eproof
\setlength{\parindent}{0cm} \setlength{\parskip}{0.5cm}

\bremark
Theorem~\ref{thm:IIPNonHd} together with Lemmas~\ref{lem:Min} and \ref{lem:TPTFEFF} explain results on simplicity of essential $C^*$-algebras of non-Hausdorff {\'e}tale groupoids in \cite{CEPSS,EP22,KM}.
\eremark

Let us now explain the connection between our results and the ones in \cite{KKLRU}, where characterisations of the ideal intersection property for groupoid $C^*$-algebras were obtained using Hamana-Furstenberg boundary groupoids. In the following, let $G$ be a non-Hausdorff {\'e}tale groupoid which can be covered by countably many open bisections, let $X = \Gu$, and let $\tiG$ be the Hausdorff cover of $G$ with $\tiX = \tiGu$. Let $\pi \colon \tiX \onto X$ be the canonical projection and $\iota \colon X \into \tiX$ the canonical inclusion. Furthermore, let $\partial X$ be the Hamana-Furstenberg boundary of $G$ (see \cite[\S~4]{KKLRU} for details), which comes with an $G$-action $G \acts \partial X$ and anchor map $\mfq \colon \partial X \to X$, and set $\partial G \coloneq G \ltimes \partial X$. The induced map $\partial G \onto G$ will be denoted by $\mfq$ as well. Let $\widetilde{\partial G}$ be the Hausdorff cover of $\partial G$ and set $\widetilde{\partial X} \coloneq \widetilde{\partial G}^{(0)}$. Let $\pi_{\partial} \colon \widetilde{\partial X} \onto \partial X$ be the canonical projection and $\iota_{\partial} \colon \partial X \into \widetilde{\partial X}$ the canonical inclusion. Finally, let $\partial \tiX$ be the Hamana-Furstenberg boundary of $\tiG$ and $\ti{\mfq} \colon \partial \tiX \to \tiX$ the anchor map for the $\tiG$-action $\tiG \acts \partial \tiX$, and set $\partial \tiG \coloneq \tiG \ltimes \partial \tiX$.
\bprop
\label{prop:mfq}
Projection onto the $G$-coordinate gives rise to a continuous surjective map $\mfp \colon \widetilde{\partial X} \onto \tiX$ sending $\menge{(g,y)}{g \in \bmx} \subseteq \partial G$ to $\bmx \subseteq G$. For every $\bmy \in \widetilde{\partial X}$, we have a group isomorphism $\bmy \congto \mfp(\bmy), \, (g,y) \ma g$. Moreover, $\mfp$ induces a bijection $\pi_{\partial}^{-1}(y) \congto \pi^{-1}(\mfq(y))$ for every $y \in \partial X$. Finally, $\tiG$ acts on $\widetilde{\partial X}$ with anchor map $\mfp$ and action given by $\bmg.\bmy \coloneq \menge{(g g_x g^{-1},y)}{g_x \in \bmx}$, where $\bmg = g \cdot \bmx$ and $\bmy = \menge{(g_x,y)}{g_x \in \bmx}$.
\eprop
\setlength{\parindent}{0cm} \setlength{\parskip}{0cm}

\bproof
It is straightforward to see that the map $\mfp$ is well-defined. To show continuity, take open sets $U_i \subseteq G$ and a compact set $K \subseteq G$ and consider $\cU(U_i,K) \coloneq \menge{\bmx \in \tiX}{\bmx \cap U_i \neq \emptyset, \, \bmx \cap K = \emptyset}$. Such sets form a basis for the topology of $\tiX$. Now set $\hat{U}_i \coloneq U_i \times \mfq^{-1}(s(U_i))$ and $\hat{K} \coloneq K \times \mfq^{-1}(s(K))$, and define $\cU(\hat{U}_i,\hat{K}) \coloneq \menge{\bmy \in \widetilde{\partial X}}{\bmy \cap \hat{U}_i \neq \emptyset, \, \bmy \cap \hat{K} = \emptyset}$. Then $\cU(\hat{U}_i,\hat{K})$ is open, and we have $\mfp^{-1}(\cU(U_i,K)) = \cU(\hat{U}_i,\hat{K})$. The remaining claims are easy to check.
\eproof
\setlength{\parindent}{0cm} \setlength{\parskip}{0.5cm}

Now let $C, D \subseteq X$ be as above, and define $D_{\partial} \coloneq s(\overline{\partial X} \setminus \partial X) \subseteq \partial X$ and $C_{\partial} \coloneq \partial X \setminus D_{\partial}$. Let $J$ be the singular ideal in $C^*_r(G)$ and $J_{\partial}$ the singular ideal in $C^*_r(\partial G)$.
\bcor
\label{cor:C*essEmbed}
We have $\mfq^{-1}(C) = C_{\partial}$ and $\mfq^{-1}(D) = D_{\partial}$. The map $\mfq^* \colon C^*_r(G) \to C^*_r(\partial G)$ dual to $\mfq$, given by $f \ma f \circ \mfq$, satisfies $(\mfq^*)^{-1}(J_{\partial}) = J$. In particular, $\mfq^*$ induces an injective $^*$-homomorphism $C^*_{\ess}(G) \to C^*_{\ess}(\partial G)$.
\ecor
\setlength{\parindent}{0cm} \setlength{\parskip}{0cm}

\bproof
The first claim follows from Proposition~\ref{prop:mfq} because $C$ consists precisely of those points $x \in X$ for which $\# \, \pi^{-1}(x) = 1$, and analogously for $C_{\partial}$. Now $\mfq^{-1}(D) = D_{\partial}$ implies $(\mfq^*)^{-1}(J_{\partial}) = J$ because of the descriptions of singular ideals as recalled at the beginning of \S~\ref{ss:DescSingId}.
\eproof
\setlength{\parindent}{0cm} \setlength{\parskip}{0.5cm}

\bremark
\label{rem:Improve}
This improves the results on embeddings of essential groupoid $C^*$-algebras in \cite[\S~5]{KKLRU} (see also \cite[Remark~6.2]{KKLRU}). As a consequence, we can replace $\sigma$-compactness by the condition that our non-Hausdorff groupoids can be covered by countably many open bisections in \cite[Theorems~A, B and D]{KKLRU}.
\eremark

\bprop
\label{prop:HFbdCover=CoverHFbd}
Assume that $X$ is compact. Then we have an isomorphism of topological groupoids $(\widetilde{\partial G})_{\ess} \cong \partial (\tiG_{\ess})$ induced by an $\tiG_{\ess}$-equivariant homeomorphism $(\widetilde{\partial X})_{\ess} \cong \partial (\tiX_{\ess})$.
\eprop
\setlength{\parindent}{0cm} \setlength{\parskip}{0cm}

\bproof
First of all, we have an $G$-action $G \acts \partial (\tiX_{\ess})$ via the anchor map $\pi \circ \ti{\mfq} \colon \partial (\tiX_{\ess}) \to \tiX_{\ess} \to X$ and action given by $g.\ti{y} \coloneq (g \cdot \ti{\mfq}(\ti{y})). \ti{y}$. Injectivity of $C(\partial X)$ allows us to extend the map $C(X) \to C(\partial (\tiX_{\ess}))$ dual to $\pi \circ \ti{\mfq}$ to a $G$-ucp map $C(\partial (\tiX_{\ess})) \to C(\partial X)$. Composing with the canonical map $C(\partial X) \to C((\widetilde{\partial X})_{\ess})$, we obtain a $G$-ucp map $\phi \colon C(\partial (\tiX_{\ess})) \to C((\widetilde{\partial X})_{\ess})$.
\setlength{\parindent}{0.5cm} \setlength{\parskip}{0cm}

We also have an $\tiG_{\ess}$-action $\tiG_{\ess} \acts \partial (\widetilde{\partial X})_{\ess}$ whose anchor map $(\widetilde{\partial X})_{\ess} \to \tiX_{\ess}$ is given by the restriction of $\mfp$ from Proposition~\ref{prop:mfq} to $(\widetilde{\partial X})_{\ess} = \overline{\iota_{\partial}(C_{\partial})}$, and the action is given by the restriction of the $\tiG$-action $\tiG \acts \widetilde{\partial X}$ from Proposition~\ref{prop:mfq}. Injectivity of $C(\partial (\tiX_{\ess}))$ allows us to extend the map $C(\tiX_{\ess}) \to C((\widetilde{\partial X})_{\ess})$ to a $G$-ucp map $\psi \colon C((\widetilde{\partial X})_{\ess}) \to C(\partial (\tiX_{\ess}))$.

Now using rigidity of the extensions $C(\tiX_{\ess}) \into C(\partial (\tiX_{\ess}))$ and $C(X) \into C(\partial X)$, and making use of the fact that the canonical projections $\tiX_{\ess} \onto X$ and $(\widetilde{\partial X})_{\ess} \onto \partial X$ are one-to-one on dense subsets, it is straightforward to see that $\phi$ and $\psi$ are inverse to each other.
\eproof
\setlength{\parindent}{0cm} \setlength{\parskip}{0.5cm}

\bremark
\label{rem:KKLRU}
Using Proposition~\ref{prop:HFbdCover=CoverHFbd}, it is straightforward to see that the characterisations of the ideal intersection property for $C^*_\ess(G)$ and $C^*_r(\tiG_{\ess})$ in \cite{KKLRU} are equivalent (under the assumption that $G$ is topologically transitive, or more generally, that $\menge{y \in \partial X}{(\partial G)_y^y \text{ is amenable}}$ is dense in $\partial X$). This applies in particular to the characterisations of the ideal intersection property in terms of essentially confined amenable sections of isotropy subgroups. 
\setlength{\parindent}{0.5cm} \setlength{\parskip}{0cm}

All in all, this explains the connection between Theorem~\ref{thm:IIPNonHd} and results characterising the ideal intersection property for groupoid $C^*$-algebras in \cite{KKLRU}.
\eremark
\setlength{\parindent}{0cm} \setlength{\parskip}{0.5cm}

\section{Amenability and nuclearity}
\label{s:AmenNuc}

Let us first recall some definitions. Throughout this section, let $G$ be a non-Hausdorff {\'e}tale groupoid with Hausdorff cover $\tiG$, and let $\iota \colon G \into \ti{G}$ be the canonical inclusion. Let $X \coloneq \Gu$ and $\tiX \coloneq \tiGu$. Moreover, for $x \in X$ and a function $\xi$ on $G$, recall that $\lambda_x(\xi) = \sum_{g \in G^x} \xi(g)$ (and $\ti{\lambda}_{\bmx}$ is defined similarly, for $\bmx \in \ti{X}$). Here we only consider $\xi$ for which the sum on the right hand side converges absolutely.

\bdefin
\label{def:CG-amen}
$G$ is called $\Cc(G)$-amenable if there exists a net $(\xi_i)$ in $\Cc(G)$ with $\xi_i \geq 0$ for all $i$ such that
\setlength{\parindent}{0cm} \setlength{\parskip}{0cm}

\begin{enumerate}
\item[(i)] $\lambda_x(\xi_i) \leq 1$ for all $x \in X$ and for all $i$,
\item[(ii)] $\lambda_x(\xi_i) \to 1$ as $i \to \infty$ uniformly on compact subsets of $X$ (with respect to $x \in X$),
\item[(iii)] $\sum_{h \in G^{r(g)}} \vert \xi_i(g^{-1} h) - \xi_i(h) \vert \to 0$ as $i \to \infty$ uniformly on compact subsets of $G$ (with respect to $g \in G$).
\end{enumerate}
\edefin
\setlength{\parindent}{0cm} \setlength{\parskip}{0.5cm}

\bdefin
\label{def:CtiG-amen}
$G$ is called $C_c(\ti{G})$-amenable if there exists a net $(\eta_i)$ in $C_c(\ti{G})$ with $\eta_i \geq 0$ for all $i$ such that the net $(\xi_i)$ given by $\xi_i \coloneq \eta_i \circ \iota$ satisfies (i), (ii) and (iii) from Definition~\ref{def:CG-amen}.
\edefin

\bdefin
\label{def:tiG-amen}
$\ti{G}$ is called topologically amenable if there exists a net $(\eta_i)$ in $C_c(\ti{G})$ with $\eta_i \geq 0$ for all $i$ such that
\setlength{\parindent}{0cm} \setlength{\parskip}{0cm}

\begin{enumerate}
\item[(i)] $\ti{\lambda}_{\bmx}(\eta_i) \leq 1$ for all $\bmx \in \ti{X}$ and for all $i$,
\item[(ii)] $\ti{\lambda}_{\bmx}(\eta_i) \to 1$ as $i \to \infty$ uniformly on compact subsets of $\ti{X}$ (with respect to $\bmx \in \tiX$),
\item[(iii)] $\sum_{\bmh \in \ti{G}^{r(\bmg)}} \vert \eta_i(\bmg^{-1} \bmh) - \eta_i(\bmh) \vert \to 0$ as $i \to \infty$ uniformly on compact subsets of $\ti{G}$ (with respect to $\bmg \in \ti{G}$).
\end{enumerate}
\edefin
\setlength{\parindent}{0cm} \setlength{\parskip}{0.5cm}

If $G$ is $\sigma$-compact, then so is $\tiG$ and we can replace nets by sequences in Definitions~\ref{def:CG-amen}, \ref{def:CtiG-amen} and \ref{def:tiG-amen} (and condition~(i) in Definitions~\ref{def:CG-amen}, \ref{def:CtiG-amen} and \ref{def:tiG-amen} is not necessary in this case). To define Borel amenability, for technical reasons we only consider sequences, following \cite{Ren15}.

\bdefin
$G$ is called Borel amenable if there exists a sequence $(\xi_i)$ of Borel functions on $G$ with $\xi_i \ge 0$ for all $i$ such that
\setlength{\parindent}{0cm} \setlength{\parskip}{0cm}

\begin{enumerate}
\item[(i)] $\lambda_x(\xi_i) \leq 1$ for all $x \in X$ and for all $i$,
\item[(ii)] $\lambda_x(\xi_i) \to 1$ as $i \to \infty$ for all $x \in X$,
\item[(iii)] $\sum_{h \in G^{r(g)}} \vert \xi_i(g^{-1} h) - \xi_i(h) \vert \to 0$ as $i \to \infty$ for all $g \in G$.
\end{enumerate}
\edefin
\setlength{\parindent}{0cm} \setlength{\parskip}{0.5cm}

\btheo
\label{thm:amenGvstiG}
Let $G$ be a non-Hausdorff {\'e}tale groupoid. Consider the following statements:
\setlength{\parindent}{0cm} \setlength{\parskip}{0cm}

\begin{enumerate}
\item[(i)] $G$ is $\Cc(G)$-amenable.
\item[(ii)] $G$ is $C_c(\ti{G})$-amenable.
\item[(iii)] $\ti{G}$ is topologically amenable.
\item[(iv)] $G$ is Borel amenable.
\end{enumerate}
We always have (i) $\Rarr$ (ii) $\LRarr$ (iii). If $G$ is $\sigma$-compact, then we also have (iii) $\Rarr$ (iv) $\Rarr$ (i).
\etheo
\setlength{\parindent}{0cm} \setlength{\parskip}{0cm}

\bproof
To see \an{(i) $\Rarr$ (ii)}, let $\eta_i \coloneq \mfi(\xi_i)$. Then $\xi_i = \eta_i \circ \iota$ by construction. Hence $(\eta_i)$ is a net satisfying the conditions in Definition~\ref{def:CtiG-amen}.
\setlength{\parindent}{0.5cm} \setlength{\parskip}{0cm}

\an{(ii) $\LRarr$ (iii)} follows by continuity of $\eta_i$ (viewed as functions on $\tiG$) and Lemma~\ref{lem:cpctGvstiG}.

Now assume that $G$ is $\sigma$-compact.
\setlength{\parindent}{0.5cm} \setlength{\parskip}{0cm}

Let us prove \an{(iii) $\Rarr$ (iv)}. As mentioned above, we now may replace nets by sequences in (iii). Now (iv) follows from (iii) because for every function $\eta \in C_c(\ti{G})$, $\eta \circ \iota$ is a bounded Borel function on $G$.

\an{(iv) $\Rarr$ (i)} follows from \cite[Theorem~2.14]{Ren15}. 
\eproof
\setlength{\parindent}{0cm} \setlength{\parskip}{0.5cm}

Now we turn to the relationship between amenability and nuclearity. For {\'e}tale groupoids which are Hausdorff, it is well known that amenability of such a groupoid is equivalent to nuclearity of its reduced groupoid $C^*$-algebra (see \cite[Chapter~6]{AR} and also \cite[Theorem~5.6.18]{BO}). For non-Hausdorff groupoids, we have the following result, whose proof follows the one of \cite[Theorem~5.6.18]{BO}.

\btheo
\label{thm:nuc=>CtiG-amen}
Let $G$ be a non-Hausdorff {\'e}tale groupoid. If $C^*_r(G)$ is nuclear, then $G$ is $C_c(\ti{G})$-amenable.

If $G$ can be covered by countably many open bisections and $C^*_{\ess}(G)$ is nuclear, then $\tiG_{\ess}$ is topologically amenable.
\etheo
\setlength{\parindent}{0cm} \setlength{\parskip}{0cm}

\bproof
We construct a net $(\zeta_i)$ of non-negative functions in the subalgebra of $B_0(G)$ generated by $\Cc(G)$ such that $\norm{\zeta_i}_{\ell^2(G_x)}^2 \le 1 $ for all $x \in G^{(0)}$, and $(\zeta_i^* *  \zeta_i)(g) \rightarrow 1$ uniformly on compact subsets of $G$ (with respect to $g \in G$). Defining $\xi_i(g) \coloneq \zeta_i(g^{-1})^2$ (for all $g \in G)$, the net $(\xi_i)$ will satisfy the conditions of Definition~\ref{def:CG-amen}. The image of the net $(\xi_i)$ in $C_0(\tiG)$ under the canonical isomorphism $C_0(\tiG) \cong C^*(\Cc(G))$ will lie in $C_c(\tiG)$ by Lemma~\ref{lem:C*GintiG}.
\setlength{\parindent}{0.5cm} \setlength{\parskip}{0cm}

Let $K \subseteq G$ be compact, $\varepsilon>0$. Let $U_1, \dotsc , U_m$ be an open covering of $K$ such that, for every $l$, there exists an open bisection $V_l$ containing a compact superset of $U_l$. Let $\underline{\pi} \coloneq \bigoplus_{x \in X} \pi_x \colon \Cc(G) \rightarrow C_r^*(G)$ denote the left regular representation. Let $f_l \in  C_c(V_l) \subseteq \Cc(G)$ be such that $0 \le f_l \le 1$ and $f_l = 1$ on $U_l$. Then $\norm{\underline{\pi}(f_l)}_{C_r^*(G)} \le 1$ and $(f_l^* * f_l)(x) = 1$ for $x \in s(U_l)$. By nuclearity, there exist positive integer $n$ and c.c.p. maps $\psi \colon C_r^*(G) \rightarrow M_n(\mathbb{C})$ and $\phi \colon M_n(\Cz) \rightarrow C_r^*(G)$ such that
$
\norm{(\phi \circ \psi) (\underline{\pi}(f_l)) - \underline{\pi}(f_l)}_{C_r^*(G)} <\varepsilon
$
for every $l$.

In the following, let $\ell_n^2$ be the Hilbert space $\Cz^n$ with the canonical inner product, and let $\{u_i\}_{i=1}^n$ denote the standard basis. For fixed $i, j \in \lge 1, \dotsc , n\rge$ let $e_{i,j} \in M_n(\Cz)$ be the matrix which has value 1 in the $(i,j)$-th entry and is 0 elsewhere, and if $a_{i,j} \in C_r^*(G)$ for all $i, j \in \lge 1, \dotsc , n\rge$, let $[a_{i,j}] \in M_n(C_r^*(G))$ denote the matrix with $(i,j)$-th entry $a_{i,j}$. Equip the algebraic tensor product $\ell_n^2 \otimes \ell_n^2 \otimes C_r^*(G)$ with the structure of a right pre-Hilbert $C_r^*(G)$-module (whose canonical inner product $\spkl{\cdot,\cdot}$ is linear in the second variable). 
Define $[b_{i,j}] \coloneq [\phi(e_{i,j})]^\frac{1}{2} \in M_n(C_r^*(G))$ and set
$\eta_\phi \coloneq \sum_{j,k} u_j \otimes u_k \otimes b_{k,j} \in \ell_n^2 \otimes \ell_n
^2 \otimes C_r^*(G)$.
For any $a \in M_n(\mathbb{C})$ we have $\phi(a) = \langle \eta_\phi, (a \otimes 1 \otimes 1) \eta_\phi \rangle$. Approximate each $b_{k, j}$ by $\zeta_{k,j} \in \Cc(G)$ so that $\norm{b_{k, j} - \underline{\pi}(\zeta_{k, j})}_{C_r^*(G)} \le \frac{\varepsilon}{n}$ and set
$\eta_\phi ' \coloneq \sum_{j,k} u_j \otimes u_k \otimes \zeta_{k,j} \in \ell_n^2 \otimes \ell_n^2 \otimes \Cc(G)$.
Note that $\eta_\phi'$ is within $\varepsilon$ of $\eta_\phi$ with respect to the $C_r^*(G)$-valued inner product. Moreover, since $\eta_\phi$ has norm at most $1$,
we may without loss of generality assume that $\eta_\phi'$ also has norm at most $1$.
    
Define the non-negative function $\zeta \colon G \to [0,\infty)$ by $\zeta(g)^2 \coloneq \sum_{j,k} \vert \zeta_{k,j}(g) \vert^2$ for $g \in G$, and observe that $\zeta$ lies in the subalgebra of $B_0(G)$ generated by $\Cc(G)$. Then
\[
\norm{\zeta}_{\ell^2(G_x)}^2 = \sum_{g \in G_x} \vert \zeta(g) \vert^2 = \sum_{j,k} (\zeta_{k,j}^* * \zeta_{k,j})(x) = \langle \eta_\phi', \eta_\phi' \rangle(x) \le 1
\]
for every $x \in G^{(0)}$, since the supremum norm is dominated by the $C_r^*(G)$-norm on $\Cc(G)$.
\setlength{\parindent}{0cm} \setlength{\parskip}{0.5cm}
    
For each $l$, write $a_l = \psi(\underline{\pi}(f_l))$, and note that $a_l$ is a contraction. Then, using the right pre-Hilbert $C_r^*(G)$-module structure of $\ell_n^2 \otimes \ell_n^2 \otimes C_r^*(G)$ we have
\begin{align*}      
    &\langle (a_l \otimes 1 \otimes 1) \eta_\phi', \ \eta_\phi' \cdot \underline{\pi}(f_l) \rangle
    =\langle (a_l \otimes 1 \otimes 1) \eta_\phi', \ \eta_\phi' \rangle * \underline{\pi}(f_l)
    =\langle \eta_\phi', \ (a_l \otimes 1 \otimes 1) \eta_\phi' \rangle^* * \underline{\pi}(f_l)\\
    \approx_{2\varepsilon} \ &\langle \eta_\phi, \ (a_l \otimes 1 \otimes 1) \eta_\phi \rangle^* * \underline{\pi}(f_l) 
    = \phi(a_l)^* * \underline{\pi}(f_l)
    \approx_\varepsilon \underline{\pi}(f_l)^* * \underline{\pi}(f_l) 
\end{align*}
where $\approx_{\varepsilon}$ means \an{within $\varepsilon$ with respect to the $C_r^*(G)$-norm}. Both sides of the above expression are elements of $\Cc(G)$, so that, for $x \in s(U_l)$, we have
\[
\vert 1 - \langle (a_l \otimes 1 \otimes 1) \eta_\phi', \eta_\phi' \cdot \underline{\pi}(f_l) \rangle(x) \vert  = \vert (f_l^* * f_l)(x) - \langle (a_l \otimes 1 \otimes 1) \eta_\phi', \eta_\phi' \cdot \underline{\pi}(f_l) \rangle(x) \vert \leq 3\varepsilon.
\]
Fix an arbitrary $g \in K$. Assume $g \in U_l$, and let $x = s(g)$. Any element of $\ell_n^2 \otimes \ell_n^2 \otimes \Cc(G)$ can be viewed as a function from $G$ to $\ell_n^2 \otimes \ell_n^2$ in the obvious way. 
For $h \in G_x$ observe that
\[
(\eta_\phi' \cdot \underline{\pi}(f_l)) (h) = \sum_{\ti{g} \in G_x} \eta_\phi'(h \ti{g}^{-1}) f_l(\ti{g}) = \eta_\phi'(h g^{-1}) \in \ell_n^2 \otimes \ell_n^2.
\]
This implies
\begin{align*}
    1 - 3\varepsilon &\le \vert \langle (a_l \otimes 1 \otimes 1) \eta_\phi', \eta_\phi' \cdot \underline{\pi}(f_l) \rangle(x) \vert
    = \bigg\vert \sum_{h \in G_x} \langle (a_l \otimes 1) \eta_\phi'(h), (\eta_\phi' \cdot \underline{\pi}(f_l)) (h) \rangle_{\ell_n^2 \otimes \ell_n^2} \bigg\vert\\
    &\le \sum_{h \in G_x} \norm{\eta_\phi'(h)}_{\ell_n^2 \otimes \ell_n^2} \norm{\eta_\phi'(h g^{-1})}_{\ell_n^2 \otimes \ell_n^2} 
    = \sum_{h \in G_x} \zeta(h) \zeta(h g^{-1}) = \sum_{\ti{h} \in G_{r(g)}} \zeta(\ti{h} g) \zeta(\ti{h})
    = (\zeta^* * \zeta) (g).
\end{align*}
But $\norm{\zeta}_{\ell^2(G_y)}^2 \le 1$ for all $y \in G^{(0)}$, so
$
1-3\varepsilon \le (\zeta^* * \zeta) (g) \le 1
$,
as desired. This proves our first claim.

The proof for the second claim is similar: Let $K$ and $f_l$ be as above, and let $\dot{\underline{\pi}} \colon \Cc(G) \rightarrow C^*_{\ess}(G)$ be the composite of $\underline{\pi}$ and the canonical projection $C^*_r(G) \onto C^*_{\ess}(G)$. Given arbitrary $\varepsilon > 0$, nuclearity of $C^*_{\ess}(G)$ implies the existence of c.c.p. maps $\dot{\psi} \colon C^*_{\ess}(G) \rightarrow M_n(\mathbb{C})$ and $\dot{\phi} \colon M_n(\mathbb{C}) \rightarrow C^*_{\ess}(G)$ such that
\[
\norm{(\dot{\phi} \circ \dot{\psi}) (\dot{\underline{\pi}}(f_l)) - \dot{\underline{\pi}}(f_l)}_{C^*_{\ess}(G)} <\varepsilon
\]
for every $l$. Now the same argument as above shows that $1-3\varepsilon \le (\zeta^* * \zeta) (g) \le 1$ for all $g \in (G_C^C) \cap K$. This shows that $\tiG_{\ess}$ is topologically amenable.
\eproof
\setlength{\parindent}{0cm} \setlength{\parskip}{0.5cm}

The following is a generalisation of the corresponding result about extending representations in \cite[\S~4]{BM25}. (Note that similar ideas already appear in \cite{Ren87} and \cite[Appendix~B]{MW08}.) Our proof follows \cite[\S~4]{BM25}, but we give details because we do not impose the constraint that $G$ can be covered by countably many open bisections. 

\blemma
\label{lem:ExtRep}
Let $G$ be a non-Hausdorff {\'e}tale groupoid. Every $^*$-algebra representation $\pi \colon \Cc(G) \ri \cL(H)$ extends to a $^*$-algebra representation $\ti{\pi} \colon B_c(G) \ri \cL(H)$. Moreover, the commutators of the ranges of $\pi$ and $\ti{\pi}$ are equal.
\elemma
\setlength{\parindent}{0cm} \setlength{\parskip}{0cm}

Before proving the lemma we introduce some terminology and notation. 
Given a Hausdorff locally compact space $U$, a \emph{positive Radon measure} on $U$ is a positive Borel measure which is inner regular on open sets, outer regular on all Borel sets, and finite on compacts. A \emph{complex Radon measure} on $U$ is a complex Borel measure of the form $\mu = \rho \ {\rm d}\nu$, where $\rho \colon U \ri \Tz = \menge{z \in \Cz}{\vert z \vert = 1}$ is a Borel-measurable function, and $\nu$ is a finite positive Radon measure on $U$. Concretely $\mu(B) = \int_B \rho \ {\rm d}\nu$ for Borel subsets $B \subseteq U$. Here $\rho$ and $\nu$ are essentially unique, in the sense that if $\mu = \ti{\rho} \ {\rm d}\ti{\nu}$ is another such representation, then $\nu = \ti{\nu}$, and $\rho= \ti{\rho}$ for $\nu$-almost all points in $U$.
Any complex Radon measure $\mu = \rho \ {\rm d}\nu$ on $U$ has an associated linear functional $\overline{\mu}$ on $C_c(U)$ given by $\overline{\mu}(h) = \int_U h \rho \ {\rm d}\nu$. By the Riesz-Markov-Kakutani representation theorem, the mapping $\mu \mapsto \overline{\mu}$ gives a one-to-one correspondence between complex Radon measures on $U$ and linear functionals $C_c(U) \ri \mathbb{C}$ which are $\norm{\cdot}_{\infty}$-bounded. Further details can be found in \cite[Chapter~7]{Fol}.

\bproof[Proof of Lemma~\ref{lem:ExtRep}]
Let $\pi \colon \Cc(G) \ri \cL(H)$ be a $^*$-algebra representation. Fix $\xi, \eta \in H$ and define the linear functional $\pi_{\xi, \eta} \colon \Cc(G) \ri \mathbb{C}$ by $h \mapsto \big\langle \xi, \pi(h) \eta \big\rangle$. (By convention the inner product $\spkl{\cdot,\cdot}$ on $H$ is linear in the second variable.) For any open bisection $U \in \B$, the restriction of $\pi$ to $C_c(U)$ is bounded with respect to the supremum norm, and hence the restriction of $\pi_{\xi, \eta}$ to $C_c(U)$ is as well. By the Riesz-Markov-Kakutani representation theorem, there exist a Borel function $\rho_{\xi, \eta}^U \colon U \ri \mathbb{T}$ and finite positive Radon measure $\nu_{\xi, \eta}^U$ on $U$ such that
$
\pi_{\xi, \eta}(h) = \int_U h \rho_{\xi, \eta}^U \ {\rm d}\nu_{\xi, \eta}^U
$
for all $h \in C_c(U)$. By uniqueness of $\rho_{\xi, \eta}^U$ and $\nu_{\xi, \eta}^U$, for any open bisections $U_1, U_2 \in \B$, the measures $\nu_{\xi, \eta}^{U_1}, \nu_{\xi, \eta}^{U_2}$ and the Borel functions $\rho_{\xi, \eta}^{U_1}, \rho_{\xi, \eta}^{U_2}$ agree on the intersection $U_1 \cap U_2$; that is, $\nu_{\xi, \eta}^{U_1}(B) = \nu_{\xi, \eta}^{U_2}(B)$ for any Borel set $B \subseteq U_1 \cap U_2$, and $\rho_{\xi, \eta}^{U_1}(g) = \rho_{\xi, \eta}^{U_2}(g)$ for $\nu_{\xi, \eta}^{U_1}$-almost every $g \in U_1 \cap U_2$. Then, for any finite set of open bisections $\cF = \{U_1, \dotsc , U_n\} \subset \B$ there exists (see for example \cite[Lemma~A.1]{MW08}) a unique finite positive Borel measure $\nu_{\xi, \eta}^\cF$ on $\bigcup_{i=1}^n U_i$ such that for any $U_i \in \cF$, $\nu_{\xi, \eta}^\cF$ agrees with $\nu_{\xi, \eta}^{U_i}$ on Borel subsets of $U_i$. 
Similarly, there exists a Borel function $\rho_{\xi, \eta}^\cF \colon \bigcup_{i=1}^n U_i \ri \mathbb{T}$ such that for each $i$, $\rho_{\xi, \eta}^\cF(g) = \rho_{\xi, \eta}^{U_i}(g)$ for $\nu_{\xi, \eta}^\cF$-almost all $g \in U_i$.
\setlength{\parindent}{0cm} \setlength{\parskip}{0.5cm}

Define the extension $\ti{\pi}_{\xi, \eta} \colon B_c(G) \ri \mathbb{C}$ as follows: given $f \in B_c(G)$, let $\cF = \lge U_1, \dotsc , U_n\rge \subset \B$ be a finite collection of open bisections covering $\osupp(f)$, and define
$\ti{\pi}_{\xi, \eta}(f) \coloneq \int_{\bigcup_{i=1}^n U_i} f \rho_{\xi, \eta}^\cF \ {\rm d}\nu_{\xi, \eta}^\cF$. The complex Borel measures $\rho_{\xi, \eta}^\cF \ {\rm d}\nu_{\xi, \eta}^\cF$ agree on intersections for different choices of finite subset $\cF \subset \B$, hence it is clear that $\ti{\pi}_{\xi, \eta}$ is well-defined and linear. For fixed $\eta \in H$, the map $\xi \mapsto \pi_{\xi, \eta}$ is anti-linear, therefore, for fixed finite $\cF \subset \B$ the map $\xi \mapsto \rho_{\xi, \eta}^\cF \ {\rm d}\nu_{\xi, \eta}^\cF$ is also anti-linear (because these complex Borel measures restrict to complex Radon measures on open bisections, and complex Radon measures are determined by their action on compactly supported continuous functions). It follows that the map $\xi \mapsto \ti{\pi}_{\xi, \eta}$ is anti-linear. Similarly, for fixed $\xi \in H$, the map $\eta \mapsto \ti{\pi}_{\xi, \eta}$ is linear. We can now define the extension $\ti{\pi} \colon B_c(G) \ri \cL(H)$ of $\pi$: given $f \in B_c(G)$, let $\ti{\pi}(f) \in \cL(H)$ be the unique operator such that
$
\ti{\pi}_{\xi, \eta}(f) = \big\langle \xi, \ti{\pi}(f) \eta \big\rangle
$
for all $\xi, \eta \in H$.  It is clear that $\ti{\pi}_{\xi, \eta}$ extends $\pi_{\xi, \eta}$ for each $\xi, \eta \in H$, and hence $\ti{\pi}(f) = \pi(f)$ for all $f \in \Cc(G)$. We show that the map $\ti{\pi} \colon B_c(H) \ri \cL(H)$ is a $^*$-algebra representation. Linearity of $\ti{\pi}$ is clear, therefore, it suffices to prove that $\ti{\pi}$ is multiplicative and $^*$-preserving for functions in $B_c(G)$ whose strict support is contained in some open bisection. First we show that $\ti{\pi}(f^*) = \ti{\pi}(f)^*$ for such $f \in B_c(G)$. Fix $U \in \B$ and let $\xi, \eta \in H$. Observe that $\ti{\pi}(\dot{f}^*) = \pi(\dot{f}^*) = \pi(\dot{f})^* = \ti{\pi}(\dot{f})^*$ for any $\dot{f} \in C_c(U)$ so that 
\setlength{\parindent}{0cm} \setlength{\parskip}{0cm}

\begin{align*}
\int_U \dot{f}(g) \rho_{\eta, \xi}^U(g) \ {\rm d}\nu_{\eta, \xi}^U(g)
&= \big\langle \eta, \pi(\dot{f}) \xi \big\rangle
= \overline{ \big\langle \xi, \pi(\dot{f}^*) \eta \big\rangle }
= \overline{ \int_{U^{-1}} \dot{f}^*(g) \rho_{\xi, \eta}^{U^{-1}}(g) \ {\rm d}\nu_{\xi, \eta}^{U^{-1}}(g) } \\
&= \int_{U^{-1}} \dot{f}(g^{-1}) \overline{\rho_{\xi, \eta}^{U^{-1}}(g)}\ {\rm d}\nu_{\xi, \eta}^{U^{-1}}(g)
= \int_U \dot{f}(g) \Big(\rho_{\xi, \eta}^{U^{-1}}\Big)^*(g)\ {\rm d}\big(\nu_{\xi, \eta}^{U^{-1}} \big)^{-1}(g)
\end{align*}
where the positive Radon measure $\big(\nu_{\xi, \eta}^{U^{-1}}\big)^{-1}$ on $U$ is defined by $\big(\nu_{\xi, \eta}^{U^{-1}}\big)^{-1}(B) \coloneq \nu_{\xi, \eta}^{U^{-1}}(B^{-1})$ for Borel $B \subseteq U$. Since $\dot{f} \in C_c(U)$ was arbitrary, the complex Radon measures $\rho_{\eta, \xi}^U \ {\rm d}\nu_{\eta, \xi}^U$ and $\Big(\rho_{\xi, \eta}^{U^{-1}}\Big)^*\ {\rm d}\big(\nu_{\xi, \eta}^{U^{-1}}\big)^{-1}$ agree on Borel subsets of $U$. Therefore, replacing $\dot{f}$ in the calculations above with some $f \in B_c(G)$ satisfying $\osupp(f) \subseteq U$ shows that $\big\langle \eta, \ti{\pi}(f) \xi \big\rangle = \overline{ \big\langle \xi, \ti{\pi}(f^*) \eta \big\rangle }$. The Hilbert space elements $\xi, \eta  \in H$ were arbitrary, hence $\ti{\pi}(f^*) = \ti{\pi}(f)^*$, as required.
\setlength{\parindent}{0cm} \setlength{\parskip}{0.5cm}

Next we show that $\ti{\pi}$ is multiplicative on $B_c(G)$. Let $U_1, U_2 \in \B$ be arbitrary open bisections. Whenever $f_1, f_2 \in B_c(G)$ have strict supports contained in $U_1$ and $U_2$ respectively, the convolution $f_1 * f_2$ has strict support contained in the open bisection $V \coloneq U_1 U_2$ and satisfies $(f_1 * f_2)(g) = f_1(gU_2^{-1})f_2(U_1^{-1}g)$ for $g \in V$. Here we use $gU_2^{-1}$ and $U_1^{-1}g$ to denote the unique elements of the sets $\menge{gg_2^{-1}}{g_2 \in U_2, s(g_2) = s(g)}$ and $\menge{g_1^{-1}g}{g_1 \in U_1, r(g_1) = r(g)}$ respectively. 
Fix $\dot{f}_1 \in C_c(U_1)$ and observe that multiplicativity of $\pi$ ensures the following for any $\dot{f}_2 \in C_c(U_2)$ and $\xi, \eta \in H$:
\begin{align*}
&\int_{U_2} \dot{f}_2(g) \rho_{\pi(\dot{f}_1)^*\xi, \eta}^{U_2}(g) \ {\rm d}\nu_{\pi(\dot{f}_1)^*\xi, \eta}^{U_2}(g)
= \big\langle \pi(\dot{f}_1)^*\xi, \pi(\dot{f}_2) \eta \big\rangle
= \big\langle \xi, \pi(\dot{f}_1)\pi(\dot{f}_2) \eta \big\rangle
= \big\langle \xi, \pi(\dot{f}_1 * \dot{f}_2) \eta \big\rangle \\
= \ & \int_V (\dot{f}_1 * \dot{f}_2)(g) \rho_{\xi, \eta}^V(g) \ {\rm d}\nu_{\xi, \eta}^V(g)
= \int_V \dot{f}_1 (gU_2^{-1}) \dot{f}_2(U_1^{-1}g) \rho_{\xi, \eta}^V(g) \ {\rm d}\nu_{\xi, \eta}^V(g)
= \int_{U_2} \dot{f}_2(g) \ {\rm d}\mu(g)
\end{align*}
where the complex measure $\mu$ on $U_1^{-1}V \subseteq U_2$ is the pushforward of the complex measure $\dot{f}_1 \Big( (\cdot) U_2^{-1} \Big) \rho_{\xi, \eta}^V \ {\rm d}\nu_{\xi, \eta}^V$ on $V$ under the homeomorphism $V \ri U_1^{-1}V; g \mapsto U_1^{-1}g$. The set of complex Radon measures is invariant under multiplication by bounded Borel functions, so $\mu$ is a complex Radon measure. Extend $\mu$ to $U_2$ by 0 and observe that the extension is still a complex Radon measure (because $U_1^{-1}V$ is open). Since $\dot{f}_2 \in C_c(U_2)$ was arbitrary, the measures $\mu$ and $\rho_{\pi(\dot{f}_1)^*\xi, \eta}^{U_2} \ {\rm d}\nu_{\pi(\dot{f}_1)^*\xi, \eta}^{U_2}$ agree on Borel subsets of $U_2$. Therefore, replacing $\dot{f}_2$ in the calculations above with some $f_2 \in B_c(G)$ satisfying $\osupp(f_2) \subseteq U_2$ shows that $\big\langle \xi, \ti{\pi}(\dot{f}_1)\ti{\pi}(f_2) \eta \big\rangle = \big\langle \xi, \ti{\pi}(\dot{f}_1 * f_2) \eta \big\rangle$. The Hilbert space elements $\xi, \eta  \in H$ were arbitrary, hence $\ti{\pi}(\dot{f}_1 * f_2) = \ti{\pi}(\dot{f}_1) \ti{\pi}(f_2)$. Taking adjoints and reversing the roles of $U_1$ and $U_2$ implies that $\ti{\pi}(f_1 * \dot{f}_2) = \ti{\pi}(f_1) \ti{\pi}(\dot{f}_2)$ whenever $\dot{f}_2 \in C_c(U_2)$ and $f_1 \in B_c(G)$ satisfies $\osupp(f_1) \subseteq U_1$. Repeating an identical argument to the one above with $\dot{f}_1$ replaced by $f_1$ gives $\ti{\pi}(f_1 * f_2) = \ti{\pi}(f_1) \ti{\pi}(f_2)$ for any $f_2 \in B_c(G)$ satisfying $\osupp(f_2) \subseteq U_2$. This completes the proof that the extension $\ti{\pi}$ is a $^*$-algebra representation.

We now address the ``moreover" part of the lemma. Let $a \in \cL(H)$ commute with the range of $\pi$, and let $\xi, \eta \in H$, $U \in \B$. Observe that
\begin{align*}
&\int_U \dot{f}(g) \rho_{a^*(\xi), \eta}^U(g) \ {\rm d}\nu_{a^*(\xi), \eta}^U(g)
= \big\langle a^*(\xi), \pi(\dot{f}) (\eta) \big\rangle
= \big\langle \xi, \big( a\pi(\dot{f}) \big) (\eta) \big\rangle
= \big\langle \xi, \big( \pi(\dot{f}) a \big) (\eta) \big\rangle \\
= \ &\overline{\big\langle a(\eta), \pi(\dot{f}^*) (\xi) \big\rangle }
=  \int_{U^{-1}} \dot{f}(g^{-1}) \overline{\rho_{a(\eta), \xi}^{U^{-1}}(g)} \ {\rm d}\nu_{a(\eta), \xi}^{U^{-1}}(g) 
= \int_U \dot{f}(g) \Big( \rho_{a(\eta), \xi}^{U^{-1}} \Big)^*(g) \ {\rm d}\big(\nu_{a(\eta), \xi}^{U^{-1}}\big)^{-1}(g)
\end{align*}
for all $\dot{f} \in C_c(U)$. Similar arguments to earlier parts of the proof show that $\big\langle \xi, \big( a\ti{\pi}(f) \big) (\eta) \big\rangle = \big\langle \xi, \big( \ti{\pi}(f) a \big) (\eta) \big\rangle$ for any $\xi, \eta \in H$ and $f \in B_c(G)$ satisfying $\osupp(f) \subseteq U$. The Hilbert space elements $\xi, \eta \in H$ were arbitrary, hence $a$ commutes with $\ti{\pi}(f)$. Any element of $B_c(G)$ can be written as a finite linear combination of ones supported inside some open bisection, hence $a$ commutes with the range of $\ti{\pi}$.
\eproof
\setlength{\parindent}{0cm} \setlength{\parskip}{0.5cm}

\bcor
\label{cor:FullC*Embed}
Let $A$ be a $C^*$-algebra. Then $\mfi \odot \id_A$ induces an isometric $^*$-homomorphism $C^*(\mfi) \otimes \id_A \colon C^*(G) \otimes_{\max} A \into C^*(\tiG) \otimes_{\max} A$. In particular, $\mfi$ induces an isometric $^*$-homomorphism $C^*(\mfi) \colon C^*(G) \into C^*(\tiG)$.
\ecor
\setlength{\parindent}{0cm} \setlength{\parskip}{0cm}

\bproof
By universality of the full groupoid $C^*$- algebra there exists a $^*$-homomorphism $C^*(\mfi) \colon C^*(G) \ri C^*(\tiG)$ extending $\mfi$. By universality of the maximal norm on $C^*(G) \odot A$, there is an induced $^*$-homomorphism $C^*(\mfi) \otimes \id_A \colon C^*(G) \otimes_{\max} A \ri C^*(\tiG) \otimes_{\max} A$. In order to show that this map is an isometry, it suffices to prove that $\norm{x}_{C^*(G) \otimes_{\max} A} \le \norm{C^*(\mfi) \odot \id_A(x)}_{C^*(\tiG) \otimes_{\max} A}$ for all $x \in \Cc(G) \odot A$. For the remainder of the proof, we identify $\Cc(G)$ with its image under $\mfi$, so that we view $\Cc(G) \odot A$ as a subset of both $C^*(G) \odot A$ and $C^*(\tiG) \odot A$.
\setlength{\parindent}{0.5cm} \setlength{\parskip}{0cm}

We claim that for any non-degenerate $^*$-algebra representation $\phi \colon C^*(G) \odot A \ri \cL(H)$ there exists a $^*$-algebra representation $\ti{\phi} \colon C^*(\tiG) \odot A \ri \cL(H)$ agreeing with $\phi$ on $\Cc(G) \odot A$. Any non-degenerate $^*$-algebra representation $C^*(G) \odot A \ri \cL(H)$ is of the form $\pi \times  \sigma$ where $\pi \colon C^*(G) \ri \cL(H)$ and $\sigma \colon A \ri \cL(H)$ are non-degenerate $^*$-algebra representations with commuting ranges (see \cite[Theorem~3.2.6]{BO}). Here $\pi \times \sigma \colon C^*(G) \odot A \ri \cL(H)$ is given by $(\pi \times \sigma)(b \otimes a) \coloneq \pi(b) \sigma(a)$ for $b \in C^*(G)$, $a \in A$. By Lemma~\ref{lem:ExtRep} there exists a $^*$-algebra representation $\ti{\pi} \colon B_c(G) \ri \cL(H)$ extending $\pi$, and whose range commutes with the range of $\sigma$. By Lemma~\ref{lem:cpctGvstiG}~(ii), the image of $C_c(\tiG)$ under the canonical isomorphism $C_0(\tiG) \cong C^*(\Cc(G))$ will be contained in $B_c(G)$. This canonical isomorphism is $^*$-preserving and multiplicative with respect to convolution on $C_c(\tiG)$. Hence, $\ti{\pi}$ restricts to a $^*$-algebra representation of $C_c(\tiG)$, which in turn induces a $^*$-algebra representation $C^*(\tiG) \ri \cL(H)$ (also denoted $\ti{\pi}$) whose range commutes with the range of $\sigma$. Therefore (by \cite[Proposition~3.1.17]{BO}) the map $\ti{\pi} \times \sigma \colon C^*(\tiG) \odot A \ri \cL(H)$ is a well-defined $^*$-algebra representation agreeing with $\pi \times \sigma$ on $\Cc(G) \odot A$. This proves the claim.

Take $x \in \Cc(G) \odot A$. Let $\phi \colon C^*(G) \odot A \ri \cL(H)$ be an arbitrary non-degenerate $^*$-algebra representation, and let $\ti{\phi} \colon C^*(\tiG) \odot A \ri \cL(H)$ be a $^*$-algebra representation agreeing with $\phi$ on $\Cc(G) \odot A$. Then $\norm{\phi(x)}_{\cL(H)} = \norm{\ti{\phi}(x)}_{\cL(H)} \le \norm{x}_{C^*(\tiG) \otimes_{\max} A}$. Taking the supremum over all non-degenerate $^*$-algebra representations $\phi$, we get $\norm{x}_{C^*(G) \otimes_{\max} A} \le \norm{x}_{C^*(\tiG) \otimes_{\max} A}$ as required.
\eproof
\setlength{\parindent}{0cm} \setlength{\parskip}{0.5cm}

\btheo
\label{thm:nuc<=CtiG-amen}
If the canonical projection $C^*(\tiG) \onto C^*_r(\tiG)$ is an isomorphism, then so is the canonical projection $C^*(G) \onto C^*_r(G)$.
\setlength{\parindent}{0.5cm} \setlength{\parskip}{0cm}

If $\tiG$ is topologically amenable, then the canonical projection $C^*(G) \onto C^*_r(G)$ is an isomorphism, and $C^*(G)$ as well as $C^*_r(G)$ are nuclear.
\etheo
\setlength{\parindent}{0cm} \setlength{\parskip}{0cm}

\bproof
For any $f \in \Cc(G)$ we have $\norm{f}_{C_r^*(G)} = \norm{\mfi(f)}_{C_r^*(\tiG)}$ and $\norm{f}_{C^*(G)} = \norm{\mfi(f)}_{C^*(\tiG)}$ by Lemma~\ref{lem:C*GintiG} and Corollary~\ref{cor:FullC*Embed} respectively. The first claim follows.
\setlength{\parindent}{0cm} \setlength{\parskip}{0.5cm}

Let $A$ be an arbitrary $C^*$-algebra. We have the following commutative diagram
\[
\xymatrix{
C^*(G) \otimes_{\max} A \ar@{->>}[d] \ar[rr]^{C^*(\mfi) \otimes \id_A} & & C^*(\tiG) \otimes_{\max} A \ar@{->>}[d] \\
C_r^*(G) \otimes_{\min} A \ar[rr]_{C^*_r(\mfi) \otimes \id_A} & & C_r^*(\tiG) \otimes_{\min} A
}
\]
where the vertical maps are the canonical projections. Now, suppose that the Hausdorff groupoid $\tiG$ is topologically amenable. By \cite[Corollary~5.6.17]{BO} the canonical projection $C^*(\tiG) \onto C^*_r(\tiG)$ is an isomorphism. Moreover, by \cite[Theorem~5.6.18]{BO} the $C^*$-algebra $C_r^*(\tiG)$ is nuclear, and therefore the vertical map on the right is an isomorphism. As $C^*(\mfi) \otimes \id_A$ is injective by Corollary~\ref{cor:FullC*Embed}, it follows that the vertical map on the left is injective, as required.
\eproof
\setlength{\parindent}{0cm} \setlength{\parskip}{0.5cm}

Let us now summarise our findings. The following is an immediate consequence of Theorems~\ref{thm:amenGvstiG}, \ref{thm:nuc=>CtiG-amen}, \ref{thm:nuc<=CtiG-amen}, \cite[Theorem~5.6.18]{BO} and \cite[Corollary~5.6.17]{BO}.
\bcor
\label{cor:AmenNuc}
Let $G$ be a non-Hausdorff {\'e}tale groupoid. Consider the following statements:
\setlength{\parindent}{0cm} \setlength{\parskip}{0cm}

\begin{enumerate}
\item[(i)] $G$ is $\Cc(G)$-amenable.
\item[(ii)] $G$ is $C_c(\ti{G})$-amenable.
\item[(iii)] $\ti{G}$ is topologically amenable.
\item[(iv)] $G$ is Borel amenable.
\item[(v)] $C^*_r(G)$ is nuclear.
\item[(vi)] $C^*_r(\tiG)$ is nuclear.
\item[(vii)] $C^*(G)$ is nuclear.
\item[(viii)] $C^*(\tiG)$ is nuclear.
\item[(ix)] The canonical projection $C^*(G) \onto C^*_r(G)$ is an isomorphism.
\item[(x)] The canonical projection $C^*(\tiG) \onto C^*_r(\tiG)$ is an isomorphism.
\end{enumerate}

We always have (i) $\Rarr$ (ii) $\LRarr$ (iii) $\LRarr$ (v) $\LRarr$ (vi) $\LRarr$ (vii) $\LRarr$ (viii) $\Rarr$ (x) $\Rarr$ (ix). If $G$ is $\sigma$-compact, then the statements (i) - (viii) are equivalent.
\ecor
\setlength{\parindent}{0cm} \setlength{\parskip}{0.5cm}

\bremark
\label{rem:Borel-amen=>nuc}
Let $G$ be a non-Hausdorff {\'e}tale groupoid. Assume that $G$ is second countable. By results of \cite{AR} and an application of Renault's disintegration theorem \cite{Ren87} (see also \cite[Appendix~B]{MW08}) in the non-Hausdorff setting, it can be shown directly that Borel amenability of $G$ implies that both $C^*(G)$ and $C^*_r(G)$ are nuclear. Similarly, by an analogue of \cite[Proposition 6.1.8]{AR} it can be shown directly that Borel amenability of $G$ implies that the canonical projection $C^*(G) \onto C^*_r(G)$ is an isomorphism.
\eremark

\bremark
\label{rem:BM}
As mentioned in the introduction, (ii) $\LRarr$ (v) $\LRarr$ (vii) $\Rarr$ (ix) in Corollary~\ref{cor:AmenNuc} can also be deduced from (an upcoming revised version of) \cite{BM25}. Our approach using Hausdorff covers, which do not feature in \cite{BM25}, avoids some of the main technicalities in \cite{BM25} by reducing to the Hausdorff case.
\eremark

\section{Examples}
\label{s:Ex}

\subsection{Group bundles and related examples}
\label{ss:ExGB}

\bex
\label{ex:NcupGamma}
Let $\Gamma$ be a non-trivial group and let $G$ be the disjoint union of $\Nz = \gekl{1, 2, 3, \dotsc}$ and $\Gamma$, $G = \Nz \cup \Gamma$. The unit space of $G$ is given by $X = \Nz \cup \gekl{1_\Gamma}$, where $1_\Gamma$ is the identity element of $\Gamma$. The topology of $X$ is such that $X$ is homeomorphic to the one-point compactification $\Nz^+ = \Nz \cup \gekl{\infty}$ of $\Nz$ via the map which is identity on $\Nz$ and sends $1_\Gamma$ to $\infty$. To describe the topology on $G$, equip $\Nz$ with the discrete topology, and a subset of $G$ is a neighbourhood of $\gamma \in \Gamma$ if and only if it contains a set of the form $(U \setminus \gekl{\infty}) \cup \gekl{\gamma}$, where $U$ is an open subset of $\Nz^+$ containing $\infty$.

In this example, the Hausdorff cover is given by the topological disjoint union $\tiG = \Nz^+ \amalg \Gamma$, where $\infty \in \Nz^+$ corresponds to the subset $\Gamma \subseteq G$. We have $\tiG_\sing = \Gamma$ and $\tiG_\ess = \Nz^+$. If $\Gamma$ is amenable, the image of $C^*_r(G)$ in $C^*_r(\tiG)$ is given by the pullback (in other words, the limit) of the following diagram
\[
 \xymatrix{
 ? \ar@{..>}[d] \ar@{..>}[r] & C^*_r(\Gamma) \ar[d]^{\chi_{\rm tr}} \\
 C(\Nz^+) \ar[r]^{\ev_{\infty}} & \Cz
 }
\]
where $\ev_\infty$ denotes the evaluation map at $\infty$ and $\chi_{\rm tr}$ is the trivial character. If $\Gamma$ is not amenable, then the canonical embedding $C^*_r(\mfi) \colon C^*_r(G) \into C^*_r(\tiG)$ is surjective. It follows that the image of the singular ideal $J$ in $C^*_r(\tiG)$ is given by $\ker(\chi_{\rm tr})$ if $\Gamma$ is amenable and by $C^*_r(\Gamma)$ if $\Gamma$ is not amenable. In both cases, we conclude that the canonical embedding $C^*_\ess(G) \into C^*_r(\tiG_\ess)$ induced by $C^*_r(\mfi)$ is surjective, i.e., an isomorphism.
\eex

\bex
Let us present an example where $C^*_r(\mfi)$ is no longer surjective. Let $G = (\Nz \times \Gamma) \cup (\Gamma \times \Gamma') \cup (\Nz \times \Gamma')$, where $\Gamma$ and $\Gamma'$ are non-trivial groups. The unit space of $G$ is given by $X = \Nz \cup \gekl{(1_\Gamma,1_{\Gamma'})} \cup \Nz \cong (\Nz^+ \amalg \Nz^+) / {}_{\infty \sim \infty'}$, where $\infty'$ denotes the point at infinity of the second copy of $\Nz^+$. To describe the topology on $G$, equip $\Nz \times \Gamma$ and $\Nz \times \Gamma'$ with the discrete topologies, and a subset of $G$ is a neighbourhood of $(\gamma,\gamma') \in \Gamma \times \Gamma'$ if and only if it contains a set of the form $((U \setminus \gekl{\infty}) \times \gekl{\gamma}) \cup \gekl{(\gamma,\gamma')} \cup ((V \setminus \gekl{\infty'}) \times \gekl{\gamma'})$, where $U$ and $V$ are open subsets of $\Nz^+$ containing $\infty$. 

In this example, the Hausdorff cover is given by $\tiG = (\Nz^+ \times \Gamma) \amalg (\Gamma \times \Gamma') \amalg (\Nz^+ \times \Gamma')$. Here $(\infty, \gamma) \in \Nz^+ \times \Gamma$ corresponds to the subset $\gekl{\gamma} \times \Gamma' \subseteq G$, and $(\infty, \gamma') \in \Nz^+ \times \Gamma'$ corresponds to the subset $\Gamma \times \gekl{\gamma'} \subseteq \Gamma \times \Gamma'$. We have $\tiG_\sing = \Gamma \times \Gamma'$ and $\tiG_\ess = (\Nz^+ \times \Gamma) \amalg (\Nz^+ \times \Gamma')$. If $\Gamma$ and $\Gamma'$ are amenable, then the image of $C^*_r(\mfi)$ is given by the limit of the following diagram
\[
 \xymatrix{
 & & ? \ar@{..>}[dll] \ar@{..>}[d] \ar@{..>}[drr] & & \\
 C(\Nz^+) \otimes C^*_r(\Gamma) \ar[dr]^{\ev_{\infty}} & & C^*_r(\Gamma \times \Gamma') \ar[dl]_{\pi_{\Gamma}} \ar[dr]^{\pi_{\Gamma'}} & & C(\Nz^+) \otimes C^*_r(\Gamma') \ar[dl]_{\ev_{\infty}} \\
 & C^*_r(\Gamma) & & C^*_r(\Gamma') &
 }
\]
where $\ev_\infty$ denotes the evaluation map at $\infty$ and $\pi_{\Gamma}$, $\pi_{\Gamma'}$ are induced by the canonical projection maps $\Gamma \times \Gamma' \onto \Gamma$, $\Gamma \times \Gamma' \onto \Gamma'$. We conclude that $C^*_r(\mfi)$ is no longer surjective.
\eex

\bex
Let us discuss the first example from \cite{Exel}. Consider the group $\Zz/2 \times \Zz/2$ and denote by $\alpha$ the generator of the first copy of $\Zz/2$ and by $\beta$ the generator of the second copy of $\Zz/2$. Consider the action of $\Zz/2 \times \Zz/2$ on $X \coloneq [-1,1] \times \gekl{0} \cup \gekl{0} \times [-1,1]$, where $\alpha$ acts via $[-1,1] \times \gekl{0} \to [-1,1] \times \gekl{0}, \, (x,0) \ma (-x,0)$ and identity on $\gekl{0} \times [-1,1]$, and $\beta$ acts via identity on $[-1,1] \times \gekl{0}$ and via $\gekl{0} \times [-1,1] \to \gekl{0} \times [-1,1], \, (0,y) \ma (0,-y)$. The groupoid $G$ is then given by the groupoid of germs for this action, $G = G((\Zz/2 \times \Zz/2) \acts X)$. Then $G = (\Zz/2 \ltimes ([-1,1] \times \gekl{0} \setminus \gekl{(0,0)})) \cup (\Zz/2 \ltimes (\gekl{0} \times [-1,1] \setminus \gekl{(0,0)})) \cup ((\Zz/2 \times \Zz/2) \ltimes \gekl{(0,0)})$. 

In this example, the Hausdorff cover is given by 
\[
 \tiG = (\Zz/2 \ltimes ([-1,1] \times \gekl{0})) \amalg (\Zz/2 \ltimes (\gekl{0} \times [-1,1])) \amalg ((\Zz/2 \times \Zz/2) \ltimes \gekl{(0,0)}),
\]
with unit space $\tiX = ([-1,1] \times \gekl{0}) \amalg (\gekl{0} \times [-1,1]) \amalg \gekl{(0,0)}$. Here $(0,0) \in [-1,1] \times \gekl{0} \subseteq \tiX$ corresponds to the subset $\gekl{[(1_{\Zz/2},1_{\Zz/2}),(0,0)], [(1_{\Zz/2},\beta),(0,0)]} \subseteq G$ and $(0,0) \in \gekl{0} \times [-1,1] \subseteq \tiX$ corresponds to the subset $\gekl{[(1_{\Zz/2},1_{\Zz/2}),(0,0)], [(\alpha, 1_{\Zz/2}),(0,0)]} \subseteq G$. The set of dangerous points is given by $D = \gekl{(0,0)} \subseteq X$, and $\pi^{-1}(D)$ consists of the three points $(0,0) \in [-1,1] \times \gekl{0}$, $(0,0) \in \gekl{0} \times [-1,1]$ and $\iota(0,0)$. It is easy to see that $\iota(D)$ becomes an isolated point in $\tiX$ and we have $\iota(D) = \inte(\pi^{-1}(D))$. The image of $C^*_r(\mfi)$ is given by the limit of the following diagram
\[
 \xymatrix{
 & ? \ar@{..>}[dl] \ar@{..>}[d] \ar@{..>}[dr] & \\
 C^*(\Zz/2 \ltimes ([-1,1] \times \gekl{0})) \ar[d]_{\ev_{(0,0)}} & C^*(\Zz/2 \times \Zz/2) \ar[dl]^{\pi_x} \ar[dr]_{\pi_y} & C^*(\Zz/2 \ltimes (\gekl{0} \times [-1,1])) \ar[d]^{\ev_{(0,0)}}\\
 C^*(\Zz/2) & & C^*(\Zz/2)
 }
\]
where $\ev_{(0,0)}$ denotes the evaluation map at $(0,0)$ and $\pi_x$, $\pi_y$ are induced by the canonical projection maps $\Zz/2 \times \Zz/2 \onto \Zz/2$ onto the first and second factor, respectively.
\eex

\bex
Let us discuss the second example from \cite{Exel}. Consider $I = \gekl{0, \dotsc, n}$ for $n \geq 4$ and the space $X = ([0,1] \times I) / { }_\sim$, where we make the identifications $(0,i) \sim (0,j)$ for all $i, j \in I$. Let us denote the class of $(0,i)$, for $i \in I$, by $0$. The alternating group $A_n$ acts canonically on $I$, hence on $[0,1] \times I$ and thus on $X$. Let $G = G(A_n \acts X)$ be the corresponding groupoid of germs. We have $G = ((0,1] \times (I \times I)) \cup (A_n \ltimes \gekl{0})$.

In this example, the Hausdorff cover is given by $\tiG = ([0,1] \times (I \times I)) \amalg (A_n \ltimes \gekl{0})$, with unit space $\tiX = ([0,1] \times I) \amalg \gekl{0}$. Here $(0,i) \in [0,1] \times I \subseteq \tiX$ corresponds to the subset $\menge{[\sigma,0]}{\sigma \in A_n, \, \sigma(i)=i} \subseteq G$. The set of dangerous points is given by $D = \gekl{0} \subseteq X$, and $\pi^{-1}(D) = \menge{(0,i)}{i \in I} \cup \gekl{0}$. It is easy to see that $\iota(D)$ becomes an isolated point in $\tiX$ and we have $\iota(D) = \inte(\pi^{-1}(D))$. The image of $C^*_r(\mfi)$ is given by the limit of the following diagram
\[
 \xymatrix{
 & ? \ar@{..>}[dl] \ar@{..>}[dr] & \\
 C([0,1]) \otimes M_n(\Cz) \ar[dr]_{\ev_0} & & C^*(A_n) \ar[dl]^{\pi} \\
 & M_n(\Cz) &
 }
\]
where $\ev_{0}$ denotes the evaluation map at $0$ and $\pi$ is induced by the canonical representation of $A_n$ in $M_n(\Cz)$ as permutation matrices. Here we make use of the canonical isomorphism $C^*_r([0,1] \times (I \times I)) \cong C([0,1]) \otimes M_n(\Cz)$.
\eex

\bex
\label{ex:NotKInj}
The following example shows that the canonical embedding $C^*_r(G) \into C^*_r(\tiG)$ is not injective in $K$-theory in general. Set $X_+ \coloneq \menge{r \e^{2 \pi i t} \in \Cz}{r \in [\frac{1}{2}, \frac{3}{2}], \, t \in [0,\frac{1}{2}]}$, $X_- \coloneq \menge{\e^{2 \pi i t} \in \Cz}{t \in [\frac{1}{2},1]}$ and $X \coloneq X_+ \cup X_-$. Consider the $(\Zz/2)$-action of $X$ determined by $r \e^{2 \pi i t} \ma (2-r) \e^{2 \pi i t}$ (so the action is trivial on $X_-$). Let $G$ be the corresponding groupoid of germs. The set of dangerous points is given by $D = \gekl{-1,1}$. The Hausdorff cover $\tiG$ is isomorphic to $(X_+ \rtimes (\Zz / 2 \Zz)) \amalg X_-$. Note that $\pi^{-1}(D)$ has empty interior. 
\setlength{\parindent}{0.5cm} \setlength{\parskip}{0cm}

Now it is straightforward to see that $K_1(C^*_r(\tiG)) \cong \gekl{0}$ but $K_1(C^*_r(G)) \not\cong \gekl{0}$.
\eex
\setlength{\parindent}{0cm} \setlength{\parskip}{0.5cm}

\bex
Let us discuss the example from \cite[\S~5.1]{CEPSS}. Let $\varphi \colon X \to X$ be a minimal homeomorphism of a compact Hausdorff space $X$ with no isolated points. Take a point $x_0 \in X$ and let $\cO(x_0)$ be the orbit of $x_0$ under $\varphi$, i.e., $\cO(x_0) = \menge{\varphi^n(x_0)}{n \in \Zz}$. The groupoid $G$ is given by $G = ((X \setminus \cO(x_0)) \cup \coprod_{\cO(x_0)} \Zz/2) \rtimes \Zz$. To describe the topology on $G$, let us denote the generator of the copy of $\Zz/2$ corresponding to $\varphi^n(x_0)$ by $a_n$. Then a subset of $(X \setminus \cO(x_0)) \cup \coprod_{\cO(x_0)} \Zz/2$ is a neighbourhood of $x \in X$ if and only if it contains an open subset of $X$ containing $x$, and it is a neighbourhood of $a_n$ if and only if it contains a subset of the form $(U \cup \gekl{a_n}) \setminus \gekl{\varphi^n(x_0)}$, where $U$ is an open subset of $X$ containing $\varphi^n(x_0)$. Then $(X \setminus \cO(x_0)) \cup \coprod_{\cO(x_0)} \Zz/2$ becomes an {\'e}tale groupoid with unit space $X$. The action of $\Zz$ on $(X \setminus \cO(x_0)) \cup \coprod_{\cO(x_0)} \Zz/2$ by groupoid automorphisms is induced by $\varphi$. By construction, the unit space of $G$ is given by $X$.

In this example, the unit space of the Hausdorff cover is given by $\tiX = X \cup \cO(x_0)$. Let us denote elements of $\cO(x_0) \subseteq \tiX$ by $\varphi^n(x_0)$ and elements of $X \subseteq \tiX$ which lie in the orbit of $x_0$ by $\widetilde{\varphi^n(x_0)}$. Note that such an element $\widetilde{\varphi^n(x_0)}$ corresponds to the subset $\gekl{\varphi^n(x_0), a_n} \subseteq G$. Let us explain the topology of $\tiX$. Points in $\cO(x_0)$ are isolated points. A sequence $(\varphi^{n_i}(x_0))$ in $\cO(x_0) \subseteq \tiX$ converges to $x \in X \subseteq \tiX$ if and only if $\varphi^{n_i}(x_0)$ converges to $x$ in $X$ originally. And a sequence $(x_i)$ in $X \subseteq \tiX$ converges to $x \in X \subseteq \tiX$ if and only if $x_i$ converges to $x$ in $X$ originally. The set of dangerous points is given by $\cO(x_0) \subseteq X$. The singular part $\tiX_\sing$ is given by $\cO(x_0) \subseteq \tiX$. Again, we have $\iota(D) = \inte(\pi^{-1}(D))$. The Hausdorff cover is given by $\tiG = (X \cup (\coprod_{\cO(x_0)} \Zz/2)) \rtimes \Zz$, and elements of $\tiG$ are denoted by $(g,k)$ with $g \in X \cup (\coprod_{\cO(x_0)} \Zz/2)$ and $k \in \Zz$.

In this example, Theorem~\ref{thm:CharCinC} tells us that $\tif \in C_c(\tiG)$ lies in the image of $\mfi \colon \Cc(G) \into C_c(\tiG)$ if and only if
\[
 \tif(\widetilde{\varphi^n(x_0)},k) = \tif(a_n,k) + \tif(\varphi^n(x_0),k)
\]
for all $n \in \Zz$ and $k \in \Zz$, where $a_n$ denotes the generator of the copy of $\Zz/2$ in $\tiG$ corresponding to $\varphi^n(x_0)$. 
\eex

\subsection{Examples from Thompson groups}
\label{ss:ExThompson}

In the following, we discuss groupoids of germs arising from the canonical action of Thompson's group $F$ on the unit interval $[0,1]$ and of Thompson's group $T$ on the circle $\Tz$. We refer the reader to \cite{CFP} for an introduction to Thompson groups.

\bex
\label{ex:F}
Consider the groupoid of germs $G = G(F \acts [0,1])$. The unit space of $G$ is given by $X = [0,1]$. The set of dangerous points is given by $D_F = (0,1) \cap \Zz[1/2]$, where $\Zz[1/2]$ denotes the set of dyadic rationals. The unit space $\tiX$ of the Hausdorff cover of $G$ is given by $\tiX \cong K_F \cup D_F$, where $K_F = ([0,1] \setminus D_F) \cup \menge{t_-,t_+}{t \in D_F}$. To describe the topology of $K_F$, let us describe convergence. A sequence $(x_n)$ in $K_F$, which is of the form $x_n = t_n$ for $t_n \in [0,1] \setminus D_F$ or $x_n = (t_n)_-$ or $x_n = (t_n)_+$ for $t_n \in D_F$, converges to $t \in [0,1] \setminus D_F$ if and only if $t_n$ converges to $t$ in $[0,1]$, converges to $t_-$ for $t \in D_F$ if and only if $(t_n)$ converges to $t$ in $[0,1]$ and $t_n < t$ or $x_n = t_-$ eventually, and converges to $t_+$ for $t \in D_F$ if and only if $t_n$ converges to $t$ in $[0,1]$ and $t_n > t$ or $x_n = t_+$ eventually. Note that $K_F$ is homeomorphic to the Cantor space. Let us further explain the topology on $K_F \cup D_F$. Points in $D_F$ are isolated. A sequence $(t_n)$ in $D_F$ converges to $t \in [0,1] \setminus D_F$ if and only if $t_n$ converges to $t$ in $[0,1]$, converges to $t_-$ for $t \in D_F$ if and only if $t_n$ converges to $t$ in $[0,1]$ and $t_n < t$ eventually, and converges to $t_+$ for $t \in D_F$ if and only if $t_n$ converges to $t$ in $[0,1]$ and $t_n > t$ eventually. An explicit homeomorphism $K_F \cup D_F \cong \tiX$ sends $t \in [0,1]  \setminus D_F \subseteq K_F$ to $t \in X$, $t \in D_F$ to $t \in X$, $t_-$ for $t \in D_F$ to the subset $\menge{[\sigma,t] \in G}{\sigma \in F, \, \sigma \equiv \id \text{ on } (t-\varepsilon,t] \text{ for some } \varepsilon > 0} = \lim_{s < t} s$, and $t_+$ for $t \in D_F$ to the subset $\menge{[\sigma,t] \in G}{\sigma \in F, \, \sigma \equiv \id \text{ on } [t, t+\varepsilon) \text{ for some } \varepsilon > 0} = \lim_{s > t} s$. Now consider the semidirect product $\Gamma \coloneq \spkl{2} \ltimes \Zz[1/2]$ with respect to the action of $\spkl{2} = \menge{2^n}{n \in \Zz} \subseteq \Qz\reg$ on $\Zz[1/2]$ by multiplication. $\Gamma$ acts partially on $K_F$ via $(a,b).t = at+b$, $(a,b).t_- = (at+b)_-$, $(a,b).t_+ = (at+b)_+$. Now consider the partial transformation groupoid
\[
 \tiG_{K_F} \coloneq \menge{(\gamma,x) \in \Gamma \times K_F}{\gamma.x \in K, \, \gamma.x = 0 \text{ if } x = 0, \, \gamma.x = 1 \text{ if } x = 1},
\]
equipped with the product topology ($\Gamma$ is viewed as a discrete group). Moreover, define
\[
 \tiG_{D_F} \coloneq \menge{(t,\gamma_-,\gamma_+)}{t \in D_F, \, \gamma_-, \gamma_+ \in \Gamma, \, \gamma_-.t = \gamma_+.t},
\]
viewed as a discrete groupoid. Note that $\tiG_{D_F}$ is Morita equivalent to the group $\Zz^2$. Now the Hausdorff cover is given by $\tiG = \tiG_{K_F} \cup \tiG_{D_F}$, and the singular part is given by $\tiG_\sing = \tiG_{D_F}$. In this case, Theorem~\ref{thm:CharCinC} tells us that $\tif \in C_c(\tiG)$ lies in the image of $\mfi \colon \Cc(G) \into C_c(\tiG)$ if and only if
$
 \sum_{\gamma_-} f(t,\gamma_-,\gamma_+) = f(\gamma_+,t_+)
$
for all $(\gamma_+,t_+) \in \tiG_{K_F}$, $t \in D_F$ and
$
 \sum_{\gamma_+} f(t,\gamma_-,\gamma_+) = f(\gamma_-,t_-)
$
for all $(\gamma_-,t_-) \in \tiG_{K_F}$, $t \in D_F$.
\eex

\bex
\label{ex:T}
Consider the groupoid of germs $G = G(T \acts \Tz)$. This case is similar to Example~\ref{ex:F} and hence we will be brief. The unit space of $G$ is given by $X = \Tz$, which we will identify with $\Rz / \Zz \cong [0,1] / { }_{0 \sim 1}$. The set of dangerous points is given by $D_T = \Tz \cap (\Zz[1/2] / \Zz)$. The unit space $\tiX$ of the Hausdorff cover of $G$ is given by $\tiX \cong K_T \cup D_T$, where $K_T = (\Tz \setminus D_T) \cup \menge{z_-,z_+}{z \in D_T}$. The topology of $K_T$ is analogous to the one on $K_F$. Again, $K_T$ is homeomorphic to the Cantor space. As above, $\Gamma = \spkl{2} \ltimes \Zz[1/2]$ acts partially on $K_T$ in an analogous way. Now consider the partial transformation groupoid
\[
 \tiG_{K_T} \coloneq \menge{(\gamma,x) \in \Gamma \times K_T}{\gamma.x \in K_T}.
\]
Note that $\tiG_{K_T}$ is the canonical groupoid model for the Cuntz algebra $\cO_2$. Moreover, define
\[
 \tiG_{D_T} \coloneq \menge{(z,\gamma_-,\gamma_+)}{z \in D_T, \, \gamma_-, \gamma_+ \in \Gamma, \, \gamma_-.z = \gamma_+.z},
\]
viewed as a discrete groupoid. Again, $\tiG_{D_T}$ is Morita equivalent to the group $\Zz^2$. Now the Hausdorff cover is given by $\tiG = \tiG_{K_T} \cup \tiG_{D_T}$, and the singular part is given by $\tiG_\sing = \tiG_{D_T}$.
\eex

\bremark
A similar analysis applies to other Thompson-like groups (acting on intervals or the circle) for which different slopes or break points are allowed.
\eremark

\subsection{Examples from self-similar groups}
\label{ss:ExSelfSim}

\begin{figure}[ht]
\centering 
\scalebox{0.8}{
\begin{tikzpicture}
\node[state] (gamma) {$\gamma$};
\node[state, right of=gamma] (gamma') {$\gamma'$};

\draw 
(gamma) edge[above] node{$x \, \vert \, x'$} (gamma');
\end{tikzpicture}
}
\caption{Arrow in Moore diagram}
\label{fig:ArrowMoore}
\end{figure}
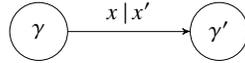

Let us discuss groupoids associated with self-similar groups. We start by briefly recalling the notion of self-similar groups and the associated groupoids (the reader may consult \cite{Nek} or \cite{StSz21,StSz23} and the references therein for details). Let $A$ be a finite alphabet and $T$ the tree given by finite words in $A$, i.e., $T = A^*$. A self-similar group is given by a subgroup $\Gamma$ of the automorphism group of $T$ such that, for all $\gamma \in \Gamma$ and $x \in T$, there exists $\gamma \vert_x \in \Gamma$ with the property that $\gamma(xw) = \gamma(x) (\gamma \vert_x) (w)$ for all $w \in T$. Note that it suffices to check this for $x \in A$. Moore diagrams are very convenient to describe self-similar groups and their actions on trees, as we will see in concrete examples below. At this point, we just recall that an arrow as in Figure~\ref{fig:ArrowMoore} in a Moore diagram means that $\gamma(x) = x'$ and $\gamma \vert_x = \gamma'$.

Given a self-similar group $\Gamma$ acting on $T$ as above, note that every $u \in T$ defines the partial bijection $T \cong uT, \, w \ma uw$ which we will denote by $u$ again. Its inverse will be denoted by $u^*$. Now let $X = \partial T$ be the boundary of the tree $T$. In other words, $X$ consists of the infinite paths in $T$, which we can also think of as infinite words in $A$. The groupoid $G$ attached to $\Gamma$ will be an ample groupoid with unit space $X$. It can be constructed in various ways. For example, we can take the monoid $M$ of partial bijections generated by $\menge{u \colon T \cong uT}{u \in T} \cup \menge{\gamma \colon T \cong T}{\gamma \in \Gamma}$. We have $M = \menge{u \gamma}{u \in T, \, \gamma \in \Gamma}$. The left inverse hull of $M$ is then the inverse semigroup $S$ whose non-zero elements (zero is the partial bijection which is nowhere defined) are given by $S\reg = \menge{u \gamma v^*}{u, v \in T, \, \gamma \in \Gamma}$. The groupoid $G$ attached to $\Gamma$ is then given by the tight groupoid of $S$. In other words, if $G(S)$ is the universal groupoid, then $G = G(S) \vert_{\partial T}$. Alternatively, to construct $G$, we can also extend the partial action (by partial bijections) $S \acts T$ to a partial action (by partial homeomorphisms) $S \acts \partial T$ and form the groupoid of germs. Explicitly, our groupoid $G$ is given by 
\[
 G = \menge{(u \gamma v^*, \omega) \in S \times \partial T}{u \gamma v^* \in S, \, \omega \in v \partial T}  / { }_{\sim},
\]
where for $u_1 \gamma_1 v_1^* \in S$ and $u_2 \gamma_2 v_2^* \in S$ with $\ell(v_2) \geq \ell(v_1)$, we have $(u_1 \gamma_1 v_1^*, \omega_1) \sim (u_2 \gamma_2 v_2^*, \omega_2)$ if and only if $\omega_1 = \omega_2$ and we have $v_2 = v_1 z$ for some $z \in T$ as well as $u_2 = u_1 \gamma_1(z)$ and $\gamma_1 \vert_z = \gamma_2$. Here and in the sequel, $\ell(\cdot)$ denotes the length of elements of $T$ (viewed as finite paths or finite words in $A$). We will write $[u \gamma v^*,\omega]$ for the class of $(u \gamma v^*,\omega)$ with respect to $\sim$.

In the following, we will focus on the case of contracting self-similar groups. This means that there exists a finite subset $N \subseteq \Gamma$ such that for all $\gamma \in \Gamma$ there exists $l \in \Nz$ such that for all $u \in T$ with $\ell(u) \geq l$, we have that $\gamma \vert_u \in N$. The smallest such $N$ is called the nucleus. Note that $n \vert_u \in N$ for all $n \in N$ and $u \in T$ if $N$ is the nucleus.

Given a contracting self-similar group $\Gamma$ with nucleus $N$, let us now describe a recipe to find the unit space $\tiX$ of the Hausdorff cover (as well as the canonical projection $\pi \colon \tiX \onto X$). As mentioned in Remark~\ref{rem:DesctiG}, this will then lead to a description of the Hausdorff cover $\tiG$ of $G$. We start by providing an approach to finding the closure $\overline{X}$ of $X$ in $G$, which will then lead to the set of dangerous points $D = s(\overline{X} \setminus X)$. In the following, we set $[u \gamma v^*, vX] \coloneq \menge{[u \gamma v^*, v \omega]}{\omega \in X}$.
\blemma
\label{lem:clX}
Let $\Gamma$ be a contracting self-similar group with nucleus $N$.
\setlength{\parindent}{0cm} \setlength{\parskip}{0cm}

\begin{enumerate}[label=(\roman*)]
    \item 
    \label{en:BGamma}
    For $\B_{\Gamma} \coloneq \menge{[u \gamma v^*, v X]}{u \gamma v^* \in S\reg}$, we have $\overline{X} = \bigcup_{U \in \B_{\Gamma}} \overline{U \cap X}^{U}$. Here $\overline{(\cdot)}^U$ denotes closure in $U$.
    \item 
    \label{en:Gamma-N}
    For every $\gamma \in \Gamma$ and $\omega \in X$, there exist $u, v \in T$ and $n \in N$ such that $[\gamma,\omega] = [u n v^*,\omega]$.
    \item
    \label{en:BN}
    For $\B_N \coloneq \menge{[u n v^*, v X]}{u n v^* \in S\reg, \, n \in N}$, we have $\overline{X} = \bigcup_{U \in \B_N} \overline{U \cap X}^{U}$.
    \item 
    \label{en:u=v}
    For $u, v \in T$ with $u \neq v$ and $n \in N$, we have $[u n v^*, vX] \cap X = \emptyset$.
    \item
    \label{en:BN'}
    For $\B_N' \coloneq \menge{[u n u^*, u X]}{u n u^* \in S\reg, \, n \in N}$, we have $\overline{X} = \bigcup_{U \in \B_N'} \overline{U \cap X}^{U}$.
    \item
    \label{en:u.u*}
    For $u \in T$ and $n \in N$, we have $[u n u^*, uX] \cap X = u ([n,X] \cap X) u^*$.
    \item
    \label{en:clX}
    We have
    \[
     \overline{X} = \bigcup_{u \in T, \, n \in N} u (\overline{[n,X] \cap X}^{[n,X]})u^*.
    \]
\end{enumerate}
\elemma
\setlength{\parindent}{0cm} \setlength{\parskip}{0cm}

\bproof
\ref{en:BGamma} is easy to see. For \ref{en:Gamma-N}, take $v \in T$ with $\ell(v) \geq l$ ($l$ is as in the definition of \an{contracting}) such that $\omega \in vX$. Then $[\gamma,\omega] = [\gamma vv^*, \omega] = [\gamma(v) \gamma \vert_v v^*, \omega]$. \ref{en:BN} follows from \ref{en:BGamma} and \ref{en:Gamma-N}. \ref{en:u=v} is straightforward to check, and \ref{en:BN'} follows from \ref{en:BN} and \ref{en:u=v}. \ref{en:u.u*} is easy to see, and \ref{en:clX} follows from \ref{en:BN'} and \ref{en:u.u*}.
\eproof
\setlength{\parindent}{0cm} \setlength{\parskip}{0.5cm}

Now, given $\gamma \in \Gamma$, set $\SF_{\gamma} = \menge{v \in T}{\gamma(v) = v, \, \gamma \vert_v = 1}$, where $1$ denotes the identity of $\Gamma$, and let $\TF_{\gamma} \coloneq \bigcup_{v \in \SF_{\gamma}} v X$. Furthermore, given $\omega \in X$ and $l \in \Nz$, let $\omega_l \in T$ be the finite word in $A$ consisting of the first $l$ letters of $\omega$. The following is easy to see and provides a way to work out $\overline{X}$.
\blemma
\label{lem:clgammaX}
For every $\gamma \in \Gamma$, we have $\overline{[\gamma,X] \cap X}^{[\gamma,X]} = \menge{[\gamma, \omega]}{\omega \in \partial \TF_{\gamma}}$, and $[\gamma,\omega]$ (for $\omega \in X$) lies in $\overline{[\gamma,X] \cap X}^{[\gamma,X]}$ if and only if for all $l \in \Nz$, $\omega_l \notin \SF_{\gamma}$ but there exists $z_l \in T$ such that $\omega_l z_l \in \SF_{\gamma}$.
\elemma

Given a contracting self-similar group $\Gamma$ with nucleus $N$, and given $\omega \in D = s(\overline{X} \setminus X)$, our goal now is to find $\menge{[u \gamma u^*, \omega] \in \overline{X}}{u \in T, \, \gamma \in \Gamma}$. First, let us define $\overline{X}_N(\omega) \coloneq \menge{[u n u^*, \omega] \in \overline{X}}{u \in T, \, n \in N}$. Moreover, given $l \in \Nz$, set $N(\omega,l) \coloneq \menge{n \in N}{[\omega_l n \omega_l^*, \omega] \in \overline{X}_N(\omega)}$ and $\overline{X}_N(\omega,l) \coloneq \menge{[\omega_l n \omega_l^*, \omega]}{n \in N(\omega,l)}$.
\blemma
\label{lem:pi-1}
Let $\Gamma$ be a contracting self-similar group with nucleus $N$, and take some $\omega \in D$.
\setlength{\parindent}{0cm} \setlength{\parskip}{0cm}

\begin{enumerate}[label=(\roman*)]
    \item 
    \label{en:clX=clXN}
    We have $\menge{[u \gamma u^*, \omega] \in \overline{X}}{u \in T, \, \gamma \in \Gamma} = \overline{X}_N(\omega) = \oX(\omega)$, where $\oX(\omega) = \oX \cap G_{\omega}^{\omega}$.
    \item 
    \label{en:clXlclXl+1}
    For every $l \in \Nz$, we have $\overline{X}_N(\omega,l) \subseteq \overline{X}_N(\omega,l+1)$.
    \item
    \label{en:clX=clXL}
    There exists $L \in \Nz$ such that 
    $\menge{[u \gamma u^*, \omega] \in \overline{X}}{u \in T, \, \gamma \in \Gamma} = \overline{X}_N(\omega,L)$.
    \item 
    \label{en:pi-1omega}
    Elements of $\pi^{-1}(\omega) \subseteq \tiX$ are subsets of the form $\menge{[\omega_L n \omega_L^*, \omega]}{n \in M}$ for some subset $M \subseteq N(\omega,L)$.
\end{enumerate}
\elemma
\setlength{\parindent}{0cm} \setlength{\parskip}{0cm}

\bproof
\ref{en:clX=clXN} follows from Lemma~\ref{lem:clX}~\ref{en:Gamma-N} and Lemma~\ref{lem:clX}~\ref{en:BN'}. To see \ref{en:clXlclXl+1}, observe that, if $\omega_{l+1} = \omega_l \omega(l+1)$, then $\omega_l n \omega_l^* \omega_{l+1} \omega_{l+1}^* = \omega_l n \omega(l+1) \omega_{l+1}^* = \omega_l n(\omega(l+1)) n \vert_{\omega(l+1)} \omega_{l+1}^* = \omega_{l+1} n \vert_{\omega(l+1)} \omega_{l+1}^*$. Now \ref{en:clX=clXL} follows from \ref{en:clX=clXN} and finiteness of $N$. Finally, \ref{en:pi-1omega} follows from \ref{en:clX=clXL}.
\eproof
\setlength{\parindent}{0cm} \setlength{\parskip}{0.5cm}

We record the following immediate consequences.
\bcor
\label{cor:SelfSimFinite}
Let $\Gamma$ be a contracting self-similar group with nucleus $N$. Then $\# \, \oX(\omega) \leq \# \, N$ for all $\omega \in X$, and $\# \, \pi^{-1}(\omega) \leq 2^{\# \, N}$ for all $\omega \in X$. 
\ecor
\setlength{\parindent}{0cm} \setlength{\parskip}{0cm}

Corollaries~\ref{cor:JC}, \ref{cor:JCc=0vsJ=0} and \ref{cor:SelfSimFinite} yield the following immediate consequence, which implies the result in \cite{GNSV} that for groupoids attached to contracting self-similar groups, the $C^*$-algebraic singular ideal $J$ vanishes if and only if the singular ideal $J_{\Cz}$ in the complex Steinberg algebra vanishes.
\bcor
\label{cor:SelfSimJvsJC}
Let $G$ be the ample groupoid attached to a contracting self-similar group. Then $J = \gekl{0}$ if and only if $J_{\Cz} = \gekl{0}$ if and only if $G$ does not satisfy condition ($\cS_0$).
\ecor
\setlength{\parindent}{0cm} \setlength{\parskip}{0.5cm}

\bremark
We arrive at the following recipe to determine $\pi^{-1}(\omega)$, for $\omega \in D$. First, find $L \in \Nz$ as in Lemma~\ref{lem:pi-1}~\ref{en:clX=clXL} as well as $N(\omega,L)$. Write $\omega = \omega_L \omega_{>L}$. Now find the set $\cM_{\omega}$ of all subsets $M \subseteq N(\omega,L)$ such that there is a sequence $(\omega_i)$ in $X$ converging to $\omega_{>L}$ such that $(\omega_i)$ converges in $G$ to all points in $\menge{[n,\omega_{>L}]}{n \in M}$ but to no points in $\menge{[n,\omega_{>L}]}{n \in N(\omega,L) \setminus M}$. Then $\omega_L \omega_i$ converges to the subset $\menge{[\omega_L n \omega_L^*, \omega]}{n \in M}$, and we have $\pi^{-1}(\omega) = \menge{\menge{[\omega_L n \omega_L^*, \omega]}{n \in M}}{M \in \cM_{\omega}}$.
\eremark

In the following, we discuss concrete examples. Our examples all belong to the class of multispinal groups (see \cite[\S~7]{StSz23}, for example).

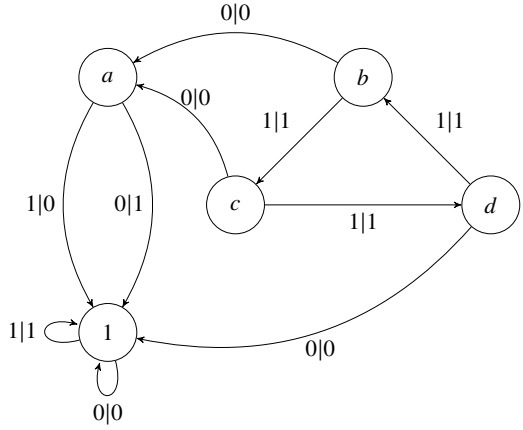
\begin{figure}[ht]
\centering 
\scalebox{0.8}{
\begin{tikzpicture}
\node[state] (b) {$b$};
\node[state, below left of=b] (c) {$c$};
\node[state, below right of=b] (d) {$d$};
\node[state, above left of=c] (a) {$a$};
\node[state, below left of=c] (1) {$1$};
\draw 
(a) edge[bend left, left] node{$0 \vert 1$} (1)
(a) edge[bend right, left] node{$1 \vert 0$} (1)
(b) edge[bend right, above] node{$0 \vert 0$} (a)
(c) edge[bend right, above] node{$0 \vert 0$} (a)
(d) edge[bend left, below] node{$0 \vert 0$} (1)
(b) edge[above left] node{$1 \vert 1$} (c)
(c) edge[below] node{$1 \vert 1$} (d)
(d) edge[above right] node{$1 \vert 1$} (b)
(1) edge[loop below] node{$0 \vert 0$} (1)
(1) edge[loop left] node{$1 \vert 1$} (1);
\end{tikzpicture}
}
\caption{Moore diagram for the Grigorchuk group}
\label{fig:Grigorchuk}
\end{figure}

\bex
\label{ex:Grig}
Let us first discuss the case of the Grigorchuk group (see Figure~\ref{fig:Grigorchuk} for its Moore diagram). The nucleus is given by $N = \gekl{1,a,b,c,d}$. 
\setlength{\parindent}{0.5cm} \setlength{\parskip}{0cm}

Following Lemmas~\ref{lem:clX} and \ref{lem:clgammaX}, we work out that the set of dangerous points is given by $D = \menge{u 1^{\infty}}{u \in T}$. Moreover, given $u \in T$, we work out $\pi^{-1}(u 1^{\infty})$ following Lemma~\ref{lem:pi-1}: Define $(u 1^{\infty})_{\bullet} \coloneq \gekl{u 1^{\infty}, [u \bullet u^*, u 1^{\infty}]}$, for $\bullet = b, c, d$. Then $\pi^{-1}(u 1^{\infty}) = \gekl{u 1^{\infty}, (u 1^{\infty})_b, (u 1^{\infty})_c, (u 1^{\infty})_d}$.

We conclude that $\tiX = (X \setminus D) \cup \menge{(u 1^{\infty})_{\bullet}}{\bullet = b, c, d; \, u \in T} \cup D$. Here points in $D \subseteq \tiX$ are isolated, and a sequence $(\bmx_i)$ in $\tiX$ converges to $\bmx \in X \setminus D \subseteq \tiX$ in $\tiX$ if and only if $\pi(\bmx_i)$ converges to $\pi(\bmx)$ in $X$, while $\bmx_i$ converges to $\bmx = (u 1^{\infty})_{\bullet} \in \tiX$ if and only if, eventually, $\bmx_i = \bmx$ or $\pi(\bmx_i) = u 1^{3n+*} 0 x_i$ for some $x_i \in X$, where $* = 0$ if $\bullet = d$, $* = 1$ if $\bullet = c$ and $* = 2$ if $\bullet = b$.

In this case we get $\iota(D) = \inte(\pi^{-1}(D))$. In particular, $\inte(\pi^{-1}(D)) \neq \emptyset$. But it is easy to see that the conditions in Corollary~\ref{cor:SuffCondSingVan} are satisfied, so that the singular ideal vanishes by Corollary~\ref{cor:SuffCondSingVan}, i.e., $J = \gekl{0}$, and $C^*_r(G) = C^*_\ess(G)$ (this is also shown in \cite[\S~5.6]{CEPSS}).
\eex
\setlength{\parindent}{0cm} \setlength{\parskip}{0.5cm}

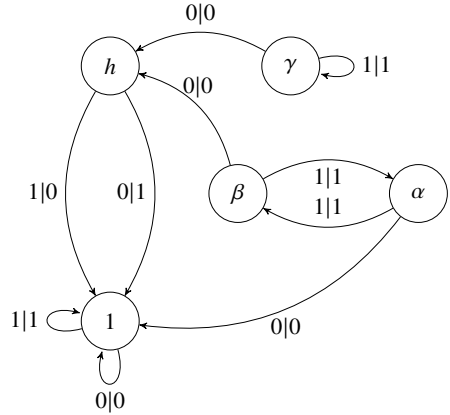
\begin{figure}[ht]
\centering 
\scalebox{0.8}{
\begin{tikzpicture}
\node[state] (h) {$h$};
\node[state, right of=h] (g) {$\gamma$};
\node[state, below right of=h] (b) {$\beta$};
\node[state, right of=b] (a) {$\alpha$};
\node[state, below left of=b] (1) {$1$};
\draw 
(h) edge[bend left, left] node{$0 \vert 1$} (1)
(h) edge[bend right, left] node{$1 \vert 0$} (1)
(g) edge[bend right, above] node{$0 \vert 0$} (h)
(g) edge[loop right] node{$1 \vert 1$} (g)
(b) edge[bend right, above] node{$0 \vert 0$} (h)
(a) edge[bend left, below] node{$0 \vert 0$} (1)
(a) edge[bend left, above] node{$1 \vert 1$} (b)
(b) edge[bend left, below] node{$1 \vert 1$} (a)
(1) edge[loop below] node{$0 \vert 0$} (1)
(1) edge[loop left] node{$1 \vert 1$} (1);
\end{tikzpicture}
}
\caption{Moore diagram for the Grigorchuk-Erschler group}
\label{fig:GrigorchukErschler}
\end{figure}

\bex
Next, let us discuss the case of the Grigorchuk-Erschler group (see Figure~\ref{fig:GrigorchukErschler} for its Moore diagram). The nucleus is given by $N = \gekl{1,h,\alpha,\beta,\gamma}$. 
\setlength{\parindent}{0.5cm} \setlength{\parskip}{0cm}

Following Lemmas~\ref{lem:clX} and \ref{lem:clgammaX}, we work out that the set of dangerous points is given by $D = \menge{u 1^{\infty}}{u \in T}$. Moreover, given $u \in T$, we work out $\pi^{-1}(u 1^{\infty})$ following Lemma~\ref{lem:pi-1}: Define $(u 1^{\infty})_{\bullet} \coloneq \gekl{u 1^{\infty}, [u \bullet u^*, u 1^{\infty}]}$, for $\bullet = \alpha, \beta$. Then $\pi^{-1}(u 1^{\infty}) = \gekl{u 1^{\infty}, (u 1^{\infty})_\alpha, (u 1^{\infty})_\beta}$.

We conclude that $\tiX = (X \setminus D) \cup \menge{(u 1^{\infty})_{\bullet}}{\bullet = \alpha, \beta; \, u \in T} \cup D$. Here points in $D \subseteq \tiX$ are isolated, and a sequence $(\bmx_i)$ in $\tiX$ converges to $\bmx \in X \setminus D \subseteq \tiX$ in $\tiX$ if and only if $\pi(\bmx_i)$ converges to $\pi(\bmx)$ in $X$, while $\bmx_i$ converges to $\bmx = (u 1^{\infty})_{\bullet} \in \tiX$ if and only if, eventually, $\bmx_i = \bmx$ or $\pi(\bmx_i) = u 1^{2n+*} 0 x_i$ for some $x_i \in X$, where $* = 0$ if $\bullet = \alpha$ and $* = 1$ if $\bullet = \beta$.

In this case we again get $\iota(D) = \inte(\pi^{-1}(D))$. In particular, $\inte(\pi^{-1}(D)) \neq \emptyset$. This time, the singular ideal is non-zero, i.e., $J \neq \gekl{0}$, and actually $J_R \neq \gekl{0}$ for any ring $R$, as shown by Nekrashevych (see also \cite[Example~7.7]{StSz23}).
\eex
\setlength{\parindent}{0cm} \setlength{\parskip}{0.5cm}

\bremark
Now let us briefly discuss Gupta-Sidki groups and GGS-groups (see \cite[\S~7.4]{StSz23}). In the language of \cite[\S~7.4]{StSz23}, if $\Phi(k)$ is an isomorphism for some $0 \leq k \leq m-2$ (which is in particular the case for Gupta-Sidki groups), then it is easy to see that $\inte(\pi^{-1}(D)) = \emptyset$. This explains why in this case, the singular ideal vanishes, i.e., $J = \gekl{0}$, and actually $J_R = \gekl{0}$ for any ring $R$ (as shown in \cite[\S~7.4]{StSz23}).
\eremark

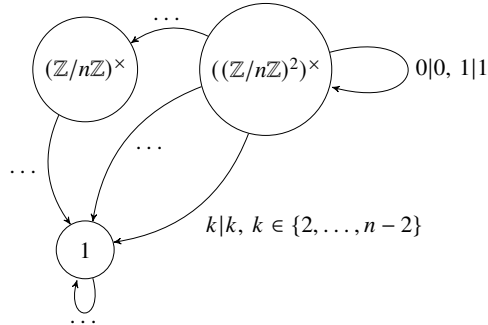
\begin{figure}[ht]
\centering 
\scalebox{0.8}{
\begin{tikzpicture}
\node[state] (Z) {$(\Zz / n \Zz)\reg$};
\node[state, right of=Z] (C) {$((\Zz / n \Zz)^2)\reg$};
\node[state, below of=Z] (1) {$1$};
\draw 
(Z) edge[bend right, left] node{$\dotso$} (1)
(C) edge[bend right, above] node{$\dotso$} (Z)
(C) edge[loop right] node{$0 \vert 0, \, 1 \vert 1$} (C)
(C) edge[bend left, below right] node{$k \vert k, \, k \in \gekl{2, \dotsc, n-2}$} (1)
(C) edge[bend right, below right] node{$\dotso$} (1)
(1) edge[loop below] node{$\dotso$} (1);
\end{tikzpicture}
}
\caption{Schematic Moore diagram for examples from \cite[\S~7.6]{StSz23}}
\label{fig:StSz7.6}
\end{figure}

\bex
Finally, we discuss the examples from \cite[\S~7.6]{StSz23}. Figure~\ref{fig:StSz7.6} shows their schematic Moore diagrams, where $n \in \Nz$ with $n > 2$, $(\Zz / n \Zz)\reg$ stands for $\Zz / n \Zz$ without identity element, $((\Zz / n \Zz)^2)\reg$ stands for $(\Zz / n \Zz)^2$ without identity element, and both $(\Zz / n \Zz)\reg$ and $((\Zz / n \Zz)^2)\reg$ have to be replaced by diagrams with $n-1$ states and $n^2-1$ states, respectively. 

Following Lemmas~\ref{lem:clX} and \ref{lem:clgammaX}, we work out that the set of dangerous points is given by 
\[
 D = \Big\{ u \omega \colon u \in T, \, \omega \in \prod_{i=1}^{\infty} \gekl{0,1} \Big\}.
\]
Moreover, given $u \in T$, it is possible to work out $\pi^{-1}(u 1^{\infty})$ following Lemma~\ref{lem:pi-1}. For example, we see that for every $u \in T$ and $\omega \in \prod_{i=1}^{\infty} \gekl{0,1}$, $\pi^{-1}(u \omega)$ always contains a subset which is isomorphic (as a group) to $(\Zz / n \Zz)^2$. Moreover, in this example, points in $\iota(D)$ are no longer isolated points in $\tiX$.
\eex


\begin{thebibliography}{99}

\bibitem{AR} C. \textsc{Anantharaman-Delaroche} and J. \textsc{Renault}, \emph{Amenable groupoids}, With a foreword by Georges Skandalis and Appendix B by E. Germain
Monogr. Enseign. Math., 36, L'Enseignement Math{\'e}matique, Geneva, 2000.

\bibitem{BHL} K. A. \textsc{Brix}, J. B. \textsc{Hume} and X. \textsc{Li}, \emph{Minimal covers with continuity-preserving transfer operators for topological dynamical systems}, preprint, arXiv:2408.11917.

\bibitem{BCFS} J. \textsc{Brown}, L. O. \textsc{Clark}, C. \textsc{Farthing} and A. \textsc{Sims}, \emph{Simplicity of algebras associated to {\'e}tale groupoids}, Semigroup Forum \emph{88} (2014), no. 2, 433--452.

\bibitem{BO} N. P. \textsc{Brown} and N. \textsc{Ozawa}, \emph{$C^*$-algebras and finite-dimensional approximations}, Grad. Stud. Math., 88, American Mathematical Society, Providence, RI, 2008.

\bibitem{BEM} A. \textsc{Buss}, R. \textsc{Exel} and R. \textsc{Meyer}, \emph{Inverse semigroup actions as groupoid actions}, Semigroup Forum \emph{85} (2012), no. 2, 227--243.

\bibitem{BM23} A. \textsc{Buss} and D. \textsc{Mart{\'i}nez}, \emph{Approximation properties of Fell bundles over inverse semigroups and non-Hausdorff groupoids}, Adv. Math. \emph{431} (2023), Paper No. 109251.

\bibitem{BM25} A. \textsc{Buss} and D. \textsc{Mart{\'i}nez}, \emph{Essential groupoid amenability and nuclearity of groupoid $C^*$-algebras}, preprint, arXiv:2501.01775.

\bibitem{CFP} J. W. \textsc{Cannon}, W. J. \textsc{Floyd} and W. R. \textsc{Parry}, \emph{Introductory notes on Richard Thompson's groups}, 
Enseign. Math. (2) \emph{42} (1996), no. 3-4, 215--256.

\bibitem{CEPSS} L. O. \textsc{Clark}, R. \textsc{Exel}, E. \textsc{Pardo}, A. \textsc{Sims} and C. \textsc{Starling}, \emph{Simplicity of algebras associated to non-Hausdorff groupoids}, Trans. Amer. Math. Soc. \emph{372} (2019), no.5, 3669--3712.

\bibitem{Con78} A. \textsc{Connes}, \emph{Sur la th{\'e}orie non commutative de l'int{\'e}gration}, Alg{\`e}bres d'op{\'e}rateurs (S{\'e}m., Les Plans-sur-Bex, 1978), pp. 19--143, 
Lecture Notes in Math., 725.

\bibitem{Con82} A. \textsc{Connes}, \emph{A survey of foliations and operator algebras}, Operator algebras and applications, Part 1 (Kingston, Ont., 1980), pp. 521--628, Proc. Sympos. Pure Math., 38, American Mathematical Society, Providence, RI, 1982.

\bibitem{Exel} R. \textsc{Exel}, \emph{Non-Hausdorff {\'e}tale groupoids}, Proc. Amer. Math. Soc. \emph{139} (2011), no.3, 897--907.

\bibitem{EP17} R. \textsc{Exel} and E. \textsc{Pardo}, \emph{Self-similar graphs, a unified treatment of Katsura and Nekrashevych $C^*$-algebras}, Adv. Math. \emph{306} (2017), 1046--1129.

\bibitem{EP22} R. \textsc{Exel} and D. R. \textsc{Pitts}, \emph{Characterizing groupoid $C^*$-algebras of non-Hausdorff {\'e}tale groupoids}, Lecture Notes in Math., 2306, Springer, Cham, 2022. 

\bibitem{Fell} J. M. G. \textsc{Fell}, \emph{A Hausdorff topology for the closed subsets of a locally compact non-Hausdorff space}, Proc. Amer. Math. Soc. \emph{13} (1962), 472--476.

\bibitem{Fol} G. B. \textsc{Folland}, \emph{Real Analysis}, Modern techniques and their applications, Second edition, Pure Appl. Math. (N. Y.), Wiley-Intersci. Publ., John Wiley \& Sons, Inc., New York, 1999.

\bibitem{GNSV} E. \textsc{Gardella}, V. \textsc{Nekrashevych}, B. \textsc{Steinberg} and A. \textsc{Vdovina}, \emph{Simplicity of $C^*$-algebras of contracting self-similar groups}, preprint, arXiv:2501.11482.

\bibitem{KKLRU} M. \textsc{Kennedy}, S. J. \textsc{Kim}, X. \textsc{Li}, S. \textsc{Raum} and D. \textsc{Ursu}, \emph{The ideal intersection property for essential groupoid $C^*$-algebras}, preprint, arXiv:2107.03980.

\bibitem{KS} M. \textsc{Khoshkam} and G. \textsc{Skandalis}, \emph{Regular representation of groupoid $C^*$-algebras and applications to inverse semigroups}, J. Reine Angew. Math. \emph{546} (2002), 47--72.

\bibitem{Kum} A. \textsc{Kumjian}, \emph{On $C^*$-diagonals}, Canad. J. Math. \emph{38} (1986), no. 4, 969--1008.

\bibitem{KM} B. K. \textsc{Kwa{\'s}niewski} and R. \textsc{Meyer}, \emph{Essential crossed products for inverse semigroup actions: simplicity and pure infiniteness}, Doc. Math. \emph{26} (2021), 271--335.

\bibitem{Li20} X. \textsc{Li}, \emph{Every classifiable simple $C^*$-algebra has a Cartan subalgebra}, Invent. Math. \emph{219} (2020), no. 2, 653--699.

\bibitem{MB} N. \textsc{Matte Bon}, \emph{Rigidity properties of full groups of pseudogroups over the Cantor set}, preprint, arXiv:1801.10133.


\bibitem{MW08} P. \textsc{Muhly} and D. P. \textsc{Williams}, \emph{Renault's equivalence theorem for groupoid crossed products}, NYJM Monogr., 3, State University of New York, University at Albany, Albany, NY, 2008.

\bibitem{Nek} V. \textsc{Nekrashevych}, \emph{Self-similar groups}, Math. Surveys Monogr., 117, American Mathematical Society, Providence, RI, 2005.

\bibitem{Nek18a} V. \textsc{Nekrashevych}, \emph{Finitely presented groups associated with expanding maps}, Geometric and cohomological group theory, 115--171, London Math. Soc. Lecture Note Ser., 444, Cambridge Univ. Press, Cambridge, 2018.

\bibitem{Nek18b} V. \textsc{Nekrashevych}, \emph{Palindromic subshifts and simple periodic groups of intermediate growth}, Ann. of Math. (2) \emph{187} (2018), no. 3, 667--719.

\bibitem{NS} S. \textsc{Neshveyev} and G. \textsc{Schwartz}, \emph{Non-Hausdorff {\'e}tale groupoids and $C^*$-algebras of left cancellative monoids}, M{\"u}nster J. Math. \emph{16} (2023), no.1, 147--175.

\bibitem{Ren} J. \textsc{Renault}, \emph{A groupoid approach to $C^*$-algebras}, Lecture Notes in Math., 793, Springer, Berlin, 1980.

\bibitem{Ren87} J. \textsc{Renault}, \emph{Repr{\'e}sentation des produits crois{\'e}s d'alg{\`e}bres de groupo{\"i}des}, J. Operator Theory \emph{18} (1987), no. 1, 67--97.

\bibitem{Ren91} J. \textsc{Renault}, \emph{The ideal structure of groupoid crossed product $C^*$-algebras}, With an appendix by Georges Skandalis, J. Operator Theory \emph{25} (1991), no. 1, 3--36.

\bibitem{Ren08} J. \textsc{Renault}, \emph{Cartan subalgebras in $C^*$-algebras}, Irish Math. Soc. Bulletin \emph{61} (2008), 29--63.

\bibitem{Ren15} J. \textsc{Renault}, \emph{Topological amenability is a Borel property}, Math. Scand. \emph{117} (2015), no. 1, 5--30.

\bibitem{StSz21} B. \textsc{Steinberg} and N. \textsc{Szak{\'a}cs}, \emph{Simplicity of inverse semigroup and {\'e}tale groupoid algebras}, Adv. Math. \emph{380} (2021), Paper No. 107611, 55 pp.

\bibitem{StSz23} B. \textsc{Steinberg} and N. \textsc{Szak{\'a}cs}, \emph{On the simplicity of Nekrashevych algebras of contracting self-similar groups}, Math. Ann. \emph{386} (2023), no. 3-4, 1391--1428.


\bibitem{Tim} T. \textsc{Timmermann}, \emph{The Fell compactification and non-Hausdorff groupoids}, Math. Z. \emph{269} (2011), no.3-4, 1105--1111.

\bibitem{Tu} J.-L. \textsc{Tu}, \emph{Non-Hausdorff groupoids, proper actions and $K$-theory}, Doc. Math. \emph{9} (2004), 565--597.

\bibitem{Yos} K. \textsc{Yoshida}, \emph{On the Simplicity of $C^*$-algebras Associated to Multispinal Groups}, preprint, arXiv:2102.02199.

\end{thebibliography}
\end{document}